\input epsf

\documentclass[11pt]{article}
\newcommand{\f}{\frac}

\usepackage{float}
\usepackage{graphicx}
\usepackage{latexsym,amsmath,amsfonts,amscd,,amssymb,subfigure}
\usepackage{graphics}
\usepackage{graphicx}
\usepackage{color}
\usepackage{epsfig}
\usepackage{changebar}
\usepackage{bm}
\numberwithin{equation}{section}

\begin{document}

\vspace{.5in}

\begin{center}

{\Large \bf A fifth-order absolutely convergent fixed-point fast sweeping hybrid alternative WENO scheme for steady state of hyperbolic conservation laws}

\end{center}

\vspace{.15in}

\centerline{
Liang Li\footnote{
School of Mathematics and Statistics, Huang Huai University, Zhumadian, Henan 463000, P.R. China. E-mail: liliangnuaa@163.com. Research was supported by the Science and Technology Project of Henan Province(252102220057) and National Natural Science Foundation of China(12501548).},
Jun Zhu\footnote{
College of Science and Key Laboratory of Mathematical Modelling and High Performance Computing of Air Vehicles (NUAA), MIIT, Nanjing University of Aeronautics and Astronautics, Nanjing, Jiangsu 210016, P.R.
China. E-mail: zhujun@nuaa.edu.cn. Research was supported by NSFC 12472292 and Science Challenge Project grant TZ2025007.},
Shanqin Chen\footnote{Department of Mathematical Sciences, Indiana University South Bend, South Bend, IN 46615, USA. 
E-mail: chen39@iu.edu}
and
Yong-Tao Zhang\footnote{Department of Applied and Computational
Mathematics and Statistics, University of Notre Dame, Notre Dame,
IN 46556, USA. E-mail: yzhang10@nd.edu. Research was partially supported by Simons Foundation MPS-TSM-00007854.}$\textsuperscript{,}
\renewcommand*{\thefootnote}{\fnsymbol{footnote}}
\setcounter{footnote}{0}\footnote{Corresponding author.}$
}
\baselineskip=1.7pc

\begin{abstract}
In this paper, we extend the previous work on absolutely convergent fixed-point fast sweeping WENO methods by Li et al. (J. Comput. Phys. 443: 110516, 2021) and design a 
fifth-order hybrid fast sweeping scheme for solving steady state problems of hyperbolic conservation laws. Unlike many other fast sweeping methods, the explicit property of fixed-point fast sweeping methods provides flexibility to apply the alternative weighted essentially non-oscillatory (AWENO) scheme with unequal-sized substencils as the local solver, which facilitates the usage of arbitrary monotone numerical fluxes. Furthermore, a novel hybrid technique is designed in the local solver to combine the nonlinear AWENO interpolation with the linear scheme for an additional improvement in efficiency of the high-order fast sweeping iterations. Numerical examples show that the developed fixed-point fast sweeping hybrid AWENO method with unequal-sized substencils can achieve absolute convergence (i.e., the residue of the fast sweeping iterations converges to machine zero / round off errors) more easily than the original AWENO method with equal-sized substencils, and is more efficient than the popular 
third-order  total variation diminishing (TVD) Runge-Kutta time-marching approach to converge to steady state solutions.    
\end{abstract}

\bigskip
{\bf Key Words:} alternative weighted essentially non-oscillatory scheme, fixed-point fast sweeping method, steady state, hyperbolic conservation laws, hybrid scheme, convergence.
\normalsize

\section{Introduction}
\label{sec1}
\setcounter{equation}{0}
\setcounter{figure}{0}
\setcounter{table}{0}

Steady state problems for hyperbolic partial differential equations (PDEs) are common
mathematical models in many applications such as fluid mechanics, wave problems, optimal
control, image processing, etc. An essential property of these boundary value problems is that the
solution information propagates along characteristics starting
from the boundary. The weighted essentially non-oscillatory (WENO) schemes are a class of effective numerical methods for spatial discretization of hyperbolic PDEs, which can achieve uniform high-order accuracy in smooth regions of the solution while maintaining sharp and essentially
non-oscillatory (ENO) \cite{HartenOsher, AWENO1} transitions of discontinuities. 
A third-order accurate
finite volume WENO scheme was first developed in \cite{LiuOsherChen}, and in \cite{JS} a general framework of constructing high-order finite difference WENO schemes was proposed. Later, WENO schemes were further developed and improved extensively in many aspects, for example, in the work on unstructured meshes to deal with problems defined on complex domain geometries \cite{HS,DK2,ZS,Y.Liu,ZQ3,Bals1}, improvement of accuracy at critical points \cite{HAPo,BorgesCarmona}, combinations with stiff system solvers for large time-step computations \cite{JZha,TanChShu,ZxuZha}, and on unequal-sized substencils \cite{JunZhu1,JunZhu2} for flexible choices of linear weights, etc. 

One of the key factors for efficiently solving steady state problems is to design fast iterative schemes for the highly coupled nonlinear systems resulted from high order nonlinear spatial discretizations such as WENO schemes. A class of efficient iterative schemes for solving hyperbolic steady state PDEs are the fast sweeping methods \cite{ZHAO}, which use alternating sweeping strategy to cover a family of characteristics in a certain direction simultaneously in each sweeping order. Coupled with the Gauss-Seidel iterations, these methods can achieve a fast convergence speed \cite{ZZQ,QZZ,WZ,ZMYZ}. Especially for steady state problems of the nonlinear hyperbolic conservation laws, the fixed-point fast sweeping WENO methods were developed \cite{S.Chen,WuLiang}, which were based on the fixed-point fast sweeping WENO methods for static Hamilton-Jacobi equations \cite{ZZC}. Unlike many other fast sweeping methods, this class of
fast sweeping methods has good properties that are suitable for complex nonlinear problems such as hyperbolic conservation laws. For example, they have explicit forms and do not involve inverse operation of nonlinear local systems, and they can be applied in solving general nonlinear hyperbolic PDEs with any monotone numerical fluxes and high order approximations (e.g. high order WENO approximations) easily. 

It was observed in \cite{WuLiang} that iteration residues of the high order fixed-point fast sweeping WENO methods do not fully converge to values of round-off error level for difficult steady state problems of hyperbolic conservation laws. The issue makes it difficult to determine the convergence criterion for the fast sweeping methods
when different problems are solved,
and challenging to apply the methods to complex problems in real applications. Furthermore, residue errors may dominate the simulation errors and affect the convergence of the numerical solution to the desired solution of the PDEs.
With the adoption of WENO schemes with unequal-sized substencils (US-WENO) \cite{JunZhu1,JunZhu2} as the local solver, fully convergent fixed-point fast sweeping methods (called absolutely convergent fixed-point fast sweeping WENO methods) were developed in \cite{LZZ,LZZ2} on both structured and unstructured meshes to solve hyperbolic conservation laws. The iteration residues of absolutely convergent fixed-point fast sweeping WENO methods fully converge to values of machine zero / round-off error level for difficult steady state problems. The methods were combined with
the inverse Lax-Wendroff procedure to handle complex numerical boundary conditions in \cite{LILW}. 

In this paper, we extend the absolutely convergent fixed-point fast sweeping finite difference WENO methods \cite{LZZ} and further improve their flexibility, accuracy, and efficiency. Instead of using the Lax-Friedrichs flux splitting in the WENO local solver \cite{LZZ}, we apply the alternative WENO (AWENO) scheme with unequal-sized substencils as the local solver, which facilitates the usage of arbitrary monotone numerical fluxes and preserves absolute convergence in the proposed fast sweeping methods. Smaller numerical errors are obtained for the fast sweeping methods here, using monotone numerical fluxes with lower dissipation. The alternative formulation for finite difference ENO schemes was first designed in \cite{AWENO1}, and was applied to WENO schemes in \cite{AWENO} for the development of AWENO schemes.
Rather than using flux splitting and
performing WENO reconstructions for flux functions as in \cite{JS}, the AWENO schemes carry out WENO interpolations for the conservative variables directly to compute numerical fluxes.
To further enhance the computational efficiency of the high-order fast sweeping iterations, a hybrid technique is used in the local solver to combine the nonlinear AWENO interpolation in stencils where the solution is identified as ``non-smooth'', with the high-order linear upwind scheme in stencils where the solution is identified as ``smooth''. Since the costly nonlinear AWENO interpolation is not applied everywhere, more efficient computations are achieved.
To detect places where the numerical solution may lose regularity and oscillations may occur with approximations by high-order polynomials,
an interesting and effective approach in \cite{ZQH,USA} is to check whether an extreme point of the reconstructed polynomial on the large stencil lies within the stencil. However, such an approach needs to identify the extreme points of a high-degree polynomial and could be complicated. Here we design a simpler method to hybrid the AWENO scheme and the higher-order linear scheme, and avoid finding the extreme points of a high-degree reconstruction polynomial. As in \cite{LILW}, the inverse Lax-Wendroff (ILW) numerical boundary procedure \cite{ILW3,S.Tan} is applied to effectively deal with complex geometries of computational domains.

The rest of the paper is organized as follows. The numerical methods are described in Section 2. In Section 3, we provide extensive numerical experiments to test and evaluate the proposed fast sweeping hybrid AWENO method, and to compare it with other methods. Numerical results show its flexibility in using different monotone numerical fluxes, absolute convergence in solving difficult benchmark problems, improved efficiency with the hybrid technique, and much faster speed than that of the popular third-order  TVD Runge-Kutta time-marching approach to converge to steady state solutions. Concluding remarks are given in Section 4.

\section {Description of the numerical methods}
\label{secfs}
\setcounter{equation}{0}
\setcounter{figure}{0}
\setcounter{table}{0}

In this section, we present the numerical schemes in detail. 
Consider the steady state problems of two dimensional (2D) hyperbolic conservation laws
\begin{equation}
f(u)_x+g(u)_y=R,  \qquad   (x,y)\in\Omega,   \\
\label{(e2.1)}
\end{equation}
on a bounded domain $\Omega$ with appropriate boundary conditions imposed on $\partial\Omega$. Here $u$ is the vector of the unknown conservative variables, $f(u)$ and $g(u)$ are the vectors of flux functions, and $R$ is the source term. For example, the steady Euler system of equations in compressible fluid dynamics has
that the conservative variables $U = (\rho, \rho u, \rho v, E)^T$, the fluxes $f= F(U) = (\rho u, \rho u^2 + p, \rho u v, u (E + p))^T$, and $g= G(U) = (\rho v, \rho u v, \rho v^2 +p, v(E + p))^T$.
Here $\rho$ is the density of the fluid, $(u,v)^T $ is the velocity vector, $p$ is the pressure, and $E=\frac{p}{\gamma'-1}+\frac{1}{2}\rho(u^2+v^2)$ is the total energy where the constant $\gamma'=1.4$ for the case of air. 

\subsection{The local solver: hybrid AWENO scheme}

We describe a hybrid AWENO scheme with unequal-sized substencils to discretize (\ref{(e2.1)}), which is the local solver of the proposed fast sweeping method.
The computational domain is partitioned by a uniform grid $\{(x_i, y_j)\}$ with the grid sizes $\Delta{x}$ and $\Delta{y}$ in the $x$-direction and the $y$-direction respectively. Denote $x_{i+1/2}=\frac{1}{2}(x_{i}+x_{i+1})$, and $y_{j+1/2}=\frac{1}{2}(y_{j}+y_{j+1})$; a computational cell $I_{i}=[x_{i-1/2},x_{i+1/2}]$ of the $x$-direction, and $J_{j}=[y_{j-1/2},y_{j+1/2}]$ of the $y$-direction. The conservation finite difference scheme is used to approximate the flux terms in (\ref{(e2.1)})
\begin{equation}
(f(u)_{x}+g(u)_{y})|_{x=x_{i},y=y_{j}}\approx\frac{1}{\Delta{x}}(\hat{f}_{i+\frac{1}{2},j}-\hat{f}_{i-\frac{1}{2},j})+\frac{1}{\Delta{y}}(\hat{g}_{i,j+\frac{1}{2}}-\hat{g}_{i,j-\frac{1}{2}}),
\end{equation}
where $\hat{f}_{i+\frac{1}{2},j}$ and $\hat{g}_{i,j+\frac{1}{2}}$ are the numerical fluxes. An alternative formulation for computing the numerical flux was designed in \cite{AWENO1}. Here we describe the procedure for computing $\hat{f}_{i+\frac{1}{2},j}$, while it is similar for $\hat{g}_{i,j+\frac{1}{2}}$. For a $(2k+1)$th-orer scheme, the desired accuracy is achieved by using
\begin{equation}
\hat{f}_{i+\frac{1}{2},j}={f}_{i+\frac{1}{2},j}+\sum_{l=1}^{k}a_{2l}\Delta{x}^{2l}(\frac{\partial^{2l}f}{\partial{x}^{2l}})_{i+\frac{1}{2},j}.
\end{equation}
The coefficients $a_{2l}, l=1,\cdots,k$ are obtained using the Taylor-expansion technique.
See \cite{AWENO1} for details. For $k=2$ (i.e., the fifth-order scheme), the numerical flux has the form
\begin{equation}
\hat{f}_{i+\frac{1}{2},j}={f}_{i+\frac{1}{2},j}-\frac{1}{24}(\Delta{x})^{2}f_{xx}|_{i+\frac{1}{2},j}+\frac{7}{5760}(\Delta{x})^{4}f_{xxxx}|_{i+\frac{1}{2},j},
\label{eq22}
\end{equation}
where the first term is approximated by a monotone numerical flux of $f(u)$, i.e.
\begin{equation}
{f}_{i+\frac{1}{2},j}=h({u}_{i+\frac{1}{2},j}^{-},{u}_{i+\frac{1}{2},j}^{+}),
\label{monflux}
\end{equation}
for linear stability of solving hyperbolic PDEs. $h(\cdot,\cdot)$ is the monotone numerical flux function of $f(u)$ obtained from an exact or approximate Riemann solver. The values of $u_{i+\frac{1}{2},j}^{-}$ and $u_{i+\frac{1}{2},j}^{+}$ are obtained by high-order approximations with nonlinear stability, e.g. the WENO approximations \cite{AWENO}.
Note that unlike the traditional finite difference WENO reconstruction process \cite{JS} in which a robust flux-splitting procedure is needed, in the alternative formulation, the finite difference WENO interpolation procedure is applied directly to numerical values of the conservative variables to obtain the cell-interface values ${u}_{i+\frac{1}{2},j}^{-}$ and ${u}_{i+\frac{1}{2},j}^{+}$. In this way, arbitrary monotone numerical fluxes can be used as in a finite volume method. 
The following subsections describe the fifth-order WENO interpolation on unequal-sized substencils, the proposed simple hybrid technique, approximations for the high-order derivatives, and the monotone numerical fluxes used here to compute the flux $\hat{f}_{i+\frac{1}{2},j}$ in (\ref{eq22}).  

\subsubsection{WENO interpolation on unequal-sized substencils}
We first briefly review the fifth-order WENO interpolation on equal-sized substencils in the AWENO scheme \cite{AWENO}, which is based on the idea of combining lower-order approximations on equal-sized small stencils to form a higher-order approximation on the whole stencil. Due to the nice dimension-by-dimension property of finite difference schemes, we only need to describe the procedure in one spatial direction, e.g. the $x-$direction here.  The following three equal-sized substencils are used for interpolation on the target cell $I_i$, which are
$T_{1}=\{I_{i-2},I_{i-1},I_{i}\},T_{2}=\{I_{i-1},I_{i},I_{i+1}\},T_{3}=\{I_{i},I_{i+1},I_{i+2}\}.$ For simplicity of presentation without causing confusion, we omit the subscript $j$ for numerical values of $u$ along the $x-$direction grid line $y=y_j$, and denote $h=\Delta x$. The corresponding Lagrange interpolation polynomials using numerical values on grid points of these three substencils (note that the middle point of a computational cell in a substencil is a grid point) are
\begin{equation}
\begin{split}
&q_1(x)=(\frac{u_{i-2}-2u_{i-1}+u_{i}}{2})(\frac{x-x_{i}}{h})^2 + (\frac{u_{i-2}-4u_{i-1}+3u_{i}}{2})(\frac{x-x_{i}}{h}) + u_{i},\\
&q_2(x)=(\frac{u_{i-1}-2u_{i}+u_{i+1}}{2})(\frac{x-x_{i}}{h})^2 + (\frac{u_{i+1}-u_{i-1}}{2})(\frac{x-x_{i}}{h}) + u_{i},\\
&q_3(x)=(\frac{u_{i}-2u_{i+1}+u_{i+2}}{2})(\frac{x-x_{i}}{h})^2 + (\frac{-3u_{i-2}+4u_{i-1}-u_{i}}{2})(\frac{x-x_{i}}{h}) + u_{i}.\\
\end{split}\label{js}
\end{equation}
The nonlinear weights $\omega_{k}$ are computed as
\begin{equation}
\omega_{k}=\frac{\bar{\omega}_{k}}{\sum_{s=1}^{3}\bar{\omega}_{s}}, \quad \bar{\omega}_{k}=\frac{\gamma_{k}}{(\varepsilon+\beta_{k})^{2}}, \qquad k=1,2,3,
\end{equation}
where $\{\gamma_{1}=\frac{1}{16},\gamma_{2}=\frac{5}{8},\gamma_{3}=\frac{5}{16}\}$ are the linear weights and $\varepsilon$ is a small positive constant to avoid a division by zero. The local smoothness indicators
$\beta_{k}= \sum_{\alpha=1}^{2} \int_{x_{i-1/2}}^{x_{i+1/2}}{h}^{2\alpha-1}(\frac{d^{\alpha}q_{k}(x)}{dx^{\alpha}})^2dx$, which measures the smoothness of the interpolation polynomial $q_k(x)$ on the substencil $T_{k}, k=1,2,3$. Then, we obtain the WENO interpolation of $u(x,y)$ at the point $(x_{i+\frac{1}{2}},y_{j})$  on the target cell $I_i$ as
\begin{equation}
{u}_{i+1/2,j}^{-}=\omega_{1}q_{1}(x_{i+1/2})+\omega_{2}q_{2}(x_{i+1/2})+\omega_{3}q_{3}(x_{i+1/2}).
\end{equation}
Numerical experiments in \cite{LZZ}
show that the fast sweeping methods using WENO local solver with unequal-sized substencils have a better convergence property than those with equal-sized substencils. This will also be verified for the proposed fast sweeping methods with the hybrid AWENO local solver in this paper.  
In the following we describe the fifth-order WENO interpolation on unequal-sized substencils  
for computing the approximation $u^{-}_{i+\frac{1}{2},j}$ on the target cell $I_i$.
The fifth-order WENO interpolation on unequal-sized substencils is a convex combination of approximations by a fourth-degree polynomial and two linear polynomials. The linear weights for the combination can be taken as any positive constant with the only requirement that their sum be one \cite{JunZhu1}. Three interpolation polynomials based on three substencils $S_{1}=\{I_{i-2},I_{i-1},I_{i},I_{i+1},I_{i+2}\},S_{2}=\{I_{i-1},I_{i}\},S_{3}=\{I_{i},I_{i+1}\}$ are
\begin{equation}
\begin{split}
p_1(x)&=(\frac{u_{i-2}-4u_{i-1}+6u_{i}-4u_{i+1}+u_{i+2}}{24})(\frac{x-x_{i}}{h})^4 \\
&+ (\frac{-u_{i-2}+2u_{i-1}-2u_{i+1}+u_{i+2}}{12})(\frac{x-x_{i}}{h})^3 \\
&+ (\frac{-u_{i-2}+16u_{i-1}-30u_{i}+16u_{i+1}-u_{i+2}}{24})(\frac{x-x_{i}}{h})^2\\
 &+ (\frac{u_{i-2}-8u_{i-1}+8u_{i+1}-u_{i+2}}{12})(\frac{x-x_{i}}{h}) + u_i,\\
p_2(x)&=\frac{u_{i}-u_{i-1}}{h}(x-x_{i})+u_{i},\\
p_3(x)&=\frac{u_{i+1}-u_{i}}{h}(x-x_{i})+u_{i}.\\
\end{split}
\label{p1S1p2p3}
\end{equation}
The local smoothness indicators
\begin{equation}
\beta_{k}= \sum_{\alpha=1}^{r_k} \int_{x_{i-1/2}}^{x_{i+1/2}}{h}^{2\alpha-1}(\frac{d^{\alpha}p_{k}(x)}{dx^{\alpha}})^2dx,\qquad k=1,2,3,
\end{equation}
measure the smoothness of $p_{k}(x)$ in each substencil $S_{k}$, where $r_1=4$ and $r_2=r_3=1$. The final expressions for evaluating $\beta_{k}$ are:
\begin{equation}
\begin{split}
\beta_{1}=&\frac{1}{60480}\bigg[1228889u_{i-1}^2 + (-3495756u_{i} - 601771u_{i-2} + 2100862u_{i+1} - 461113u_{i+2})u_{i-1}\\
  & + 82364u_{i-2}^2 + (799977u_{i} - 461113u_{i+1} + 98179u_{i+2})u_{i-2} +1228889u_{i+1}^2 \\
  & + (-3495756u_{i} - 601771u_{i+2})u_{i+1} + 2695779u_{i}^2 + 799977u_{i}u_{i+2}+82364u_{i+2}^2\bigg]\\
\beta_{2}=&(u_{i-1}-u_i)^2,\\
\beta_{3}=&(u_{i}-u_{i+1})^2.\\
\end{split}
\end{equation}
The nonlinear weights are computed as in \cite{BorgesCarmona}:
\begin{equation}
\omega_{l_1}=\frac{\bar{\omega}_{l_1}}{\sum_{l_2=1}^{3}\bar{\omega}_{l_2}}, \quad \bar{\omega}_{l_1}=\gamma_{l_1}(1+\frac{\tau}{\varepsilon+\beta_{l_1}}), \qquad l_1=1,2,3.
\end{equation}
Here the small positive constant $\varepsilon$ is taken as  $10^{-6}$ to avoid the case that the denominator becomes zero. Based on a balance between the sharp and essentially non-oscillatory shock transitions in nonsmooth regions and high-order accuracy in smooth regions, following the practice in \cite{Y.Liu,JunZhu1}, we choose the linear weights as  $\gamma_{1}=0.98$, $\gamma_{2}=0.01$, $\gamma_{3}=0.01$. $\tau$ is a quantity that depends on the absolute differences between the smoothness indicators: $\tau=(\frac{\sum_{l_1=2}^{3}|\beta_{1}-\beta_{l_1}|}{2})^2$.
The final fifth-order WENO interpolation on unequal-sized substencils for the approximation ${u}_{i+1/2,j}^{-}$ on the target cell $I_i$ is 
\begin{equation}
{u}_{i+1/2,j}^{-}=\omega_{1}\big[\frac{1}{\gamma_{1}}p_{1}(x_{i+1/2})-\frac{\gamma_{2}}{\gamma_{1}}p_{2}(x_{i+1/2})-\frac{\gamma_{3}}{\gamma_{1}}p_{3}(x_{i+1/2})\big]+\omega_{2}p_{2}(x_{i+1/2})+\omega_{3}p_{3}(x_{i+1/2}).
\label{weno5formu}
\end{equation}
The computation of approximation ${u}_{i+1/2,j}^{+}$ on the target cell $I_{i+1}$, which is needed in the numerical flux (\ref{monflux}), is mirror-symmetric with respect to $x_{i+1/2}$.  

\subsubsection{A novel hybrid interpolation method}
To further improve the computational efficiency and accuracy of numerical solutions of WENO schemes, an effective approach is to hybridize the nonlinear WENO approximations with the corresponding high-order approximations of linear upwind schemes. The costly nonlinear WENO approximations are only applied to cells which are identified as ``troubled-cells'', i.e., where the numerical solution may lose regularity.  An interesting and effective approach to identify troubled-cells, which was proposed in \cite{ZQH},  is to determine whether the reconstructed high-order polynomial (i.e., $p_{1}(x)$ here) has an extreme point in the large stencil (i.e., $S_{1}$ here). If at least one extreme point of the reconstructed polynomial on the large stencil lies within the stencil, then either oscillations may have occurred with such an approximation by this high-order polynomial, or there are smooth local extrema in the numerical solution here. For both cases, a nonlinear WENO approximation is applied to avoid possible oscillations or nonlinear instability. Otherwise, if there is no extreme point of the reconstructed high-order polynomial in the stencil, the target cell is not identified as a 
troubled-cell, and the high-order polynomial is directly used to compute the numerical values at the cell interfaces. Although this approach may mis-flag smooth local extrema in the numerical solution as troubled-cells, it is simple and free of parameters, which is a nice property compared with many other indicators for troubled-cells \cite{QSiam}. However, the procedure to identify extreme points of the reconstructed high-order polynomial (e.g., the quartic polynomial $p_{1}(x)$ here) often involves relatively complicated algebraic operations and may not be efficient.  This motivates us to design a simplified and more efficient approach here. Note the fact that if the solution is not smooth which results in nonlinear instability / oscillation (i.e., extreme point appears in the approximation polynomial) in the large stencil $S_{1}$, then at least one of these three small substencils $T_k,k=1,2,3$ must contain the non-smooth region, and the associated quadratic approximation polynomial must have extreme point in that substencil. So instead of identifying extreme points for a high-order polynomial, we can just identify extreme points for these simple quadratic polynomials $q_k(x), k=1, 2, 3$. If $q_k(x)$ has no extreme points within the corresponding  substencil $T_k$ for all $k=1,2,3$, the target cell $I_i$ is not identified as a 
troubled-cell, and the high-order interpolation polynomial $p_1(x)$ is directly used to compute the numerical values at the cell interfaces of $I_i$ to provide the values needed in the monotone flux. The extreme points $root_{k}$ corresponding to the quadratic polynomial $q_{k}(x)$ for $k=1,2,3$, which are the roots of $q'_k(x)$, are the following:
\begin{equation}\begin{split}
&root_{1}=-\frac{u_{i-2}-4u_{i-1}+3u_{i}}{2(u_{i-2}-2u_{i-1}+u_{i})}h+x_{i},\\
&root_{2}=-\frac{u_{i+1}-u_{i-1}}{2u_{i-1}-4u_{i}+2u_{i+1}}h+x_{i},\\
&root_{3}=-\frac{-3u_{i-2}+4u_{i-1}-u_{i}}{2(u_{i}-2u_{i+1}+u_{i+2})}h+x_{i}.\\
\end{split}\end{equation}
Hence, if these three extreme points satisfy the condition $root_{k} \notin T_{k}$, for all $k=1,2,3$, then the quartic polynomial $p_1(x)$ in (\ref{p1S1p2p3}), which is obtained by interpolation on the large stencil $S_{1}$, is used for the fifth-order approximation
\begin{equation}{u}_{i+1/2,j}^{-}=p_{1}(x_{i+1/2})\label{equ}\end{equation}
as in the linear upwind scheme; otherwise, the fifth-order WENO interpolation on unequal-sized substencils  (\ref{weno5formu}) is used. 
This new method avoids the complex procedure to find extreme points of high-order polynomials in \cite{ZQH,USA} and greatly simplifies the computation in the hybrid method. In particular, if it is found that the extreme point of one of the quadratic polynomials falls within the corresponding substencil, we mark the target cell as a troubled-cell and stop computing the other extreme point(s), to save simulation time. Numerical experiments in the next section show that this hybrid interpolation scheme preserves the absolute convergence for the proposed fast sweeping AWENO method. 

\subsubsection{Completion of the scheme}

For the other high-order derivative terms in (\ref{eq22}), as pointed out in \cite{AWENO}, lower-order
approximations for them are sufficient to guarantee the desired fifth-order accuracy of the scheme and they contribute much less to spurious oscillations due to at least
$\Delta{x}^2$ in the coefficients. Hence the second and fourth derivative terms here are approximated by the central finite difference schemes as
\begin{equation}
f_{xx}|_{i+\frac{1}{2},j}=\frac{1}{48\Delta{x}^2}(-5f_{i-2,j}+39f_{i-1,j}-34f_{i,j}-34f_{i+1,j}+39f_{i+2,j}-5f_{i+3,j}),
\end{equation}
\begin{equation}
f_{xxxx}|_{i+\frac{1}{2},j}=\frac{1}{2\Delta{x}^4}(f_{i-2,j}-3f_{i-1,j}+2f_{i,j}+2f_{i+1,j}-3f_{i+2,j}+f_{i+3,j}).
\end{equation}
To show the flexibility of using different monotone numerical fluxes in the  proposed fast sweeping hybrid AWENO method, we apply the HLLC flux \cite{HLLC3,HLLC2,HLIU} which has low dissipation and the local Lax-Friedrichs (LLF) flux \cite{LLF} which is simple to handle complex problems.
The HLLC flux for two-dimensional Euler equations in the $x-$direction is 
\begin{equation}
{H}^{^{HLLC}}(U_{l},U_{r})=\left\{
             \begin{array}{ll}
             F_{l},                                  & if \ S_{L}>0,  \\
             F_{l}+S_{L}(U_{l}^{*}-U_{l}), & if \  S_{L}\leq{0}<S_{M},\\
             F_{r}+S_{R}(U_{r}^{*}-U_{r}), & if \  S_{M}\leq{0}<S_{R},\\
             F_{r},                                  & if \  S_{R}\leq0,
             \end{array}
\right.
\end{equation}
where $F_{l}=F(U_{l})$ and $F_{r}=F(U_{r})$, $U_{l}$ and $U_{r}$ are the left and right states. Note that here a quantity with subscript $l$ denotes the numerical value from the left direction, and that with subscript $r$ denotes the numerical value from the right direction.   $U_{l}^{*}$ and $U_{r}^{*}$ are the averaged states between the two acoustic-wave speeds $S_{L}$, $S_{R}$.  The acoustic-wave speeds are computed as follows:
$$S_{L}=\min[v_{l}-c_{l},{v}^{*}-{c}^{*}], \quad S_{R}=\max[v_{r}+c_{r},{v}^{*}+{c}^{*}],$$
where
\begin{equation}
\left\{
             \begin{array}{l}
                          {H}^{*}=\frac{H_{l}+H_{r}R}{1+R},\\
             R=\sqrt{\frac{\rho_{r}}{\rho_{l}}}, \\
             {v}^{*}=\frac{v_{l}+v_{r}R}{1+R},                                   \\
             {c}^{*}=\sqrt{(\gamma-1)({H}^{*}-\frac{1}{2}{{v}^{*2}}) }.
             \end{array}
\right.
\end{equation}
Here, $H=\frac{E+p}{\rho}$ is the enthalpy, and $c=\sqrt{\frac{\gamma p}{\rho}}$ is the sound speed. The left average state $U_{l}^{*}$ is computed as follows
$$
\begin{gathered}
U_{l}^{*}=\rho_{l}(\frac{S_{L}-v_{l}}{S_{L}-S_{M}})\begin{pmatrix} 1 \\ u_{l}\\ S_{M}\\ \frac{E_{l}}{\rho_{l}}+(S_{M}-v_{l})[S_{M}+\frac{p_{l}}{\rho_{l}(S_{L}-v_{l})}] \end{pmatrix},
\end{gathered}
$$
where
$$S_{M}=\frac{\rho_{r}v_{r}(S_{R}-v_{r})-\rho_{l}v_{l}(S_{L}-v_{l})+p_{l}-p_{r}}{\rho_{r}(S_{R}-v_{r})-\rho_{l}(S_{L}-v_{l})}.$$
The right average state $U_{r}^{*}$ is obtained symmetrically.
The flux formula in the $y-$direction is similar. See e.g. \cite{Torobook}.
For the LLF flux, in the $x-$direction it is defined by
\begin{equation}
H^{^{LLF}}(U_{l},U_{r})=\frac{1}{2}[F(U_{l})+F(U_{r})-\lambda(U_{r}-U_{l})],
\end{equation}
where $\lambda$ is the upper bound between $U_l$ and $U_r$ for the absolute values of eigenvalues of the Jacobian matrix $F'(U)$. Similarly for the $y-$direction. 

\bigskip
\noindent{\bf Remark.} Note that an alternative WENO scheme with unequal-sized substencils was also adopted in \cite{USA}. However, we would like to point out that a different hybrid scheme is designed in this paper, which is simpler and shows a better absolutely convergence property in the proposed fast sweeping methods. This point will also be emphasized in the section for numerical experiments. 

\subsection{Absolutely convergent fixed-point fast sweeping scheme}

Discretization of the PDE (\ref{(e2.1)}) by the fifth-order hybrid AWENO scheme with unequal-sized substencils, described in Section 2.1, results in a
nonlinear algebraic system
$$0=-(\hat{f}_{i+1/2,j}-\hat{f}_{i-1/2,j})/  \Delta x-(\hat{g}_{i,j+1/2}-\hat{g}_{i,j-1/2})/\Delta y+R(u_{ij},x_i,y_j)$$
\begin{equation}
\hspace{2in}i=1,\cdots,N; \quad j=1,\cdots, M,
\label{eq3}
\end{equation}
where $N$ and $M$ are the number of grid points in the $x$ and $y$ directions, respectively.
The right-hand-side (RHS) of (\ref{eq3}) is a nonlinear function of the
numerical values at the grid points of the AWENO scheme's stencil, which is denoted by $L$ in the following. 
A popular method to compute steady state solution of hyperbolic conservation laws is to perform time-marching iterations using the third-order TVD Runge-Kutta method \cite{AWENO1}. Another class of efficient iterative methods for steady state hyperbolic conservation laws is the 
fixed-point fast sweeping scheme \cite{WuLiang,LZZ}, which has improved efficiency over the TVD Runge-Kutta method. Here we apply the fixed-point fast sweeping scheme to the nonlinear system (\ref{eq3}) obtained by the fifth-order hybrid AWENO scheme with unequal-sized substencils, which has the following form:
$$u_{i,j}^{n+1}=u_{i,j}^{n}+\frac{\gamma}{\alpha_{x}/\Delta x +\alpha_{y}/\Delta  y}L(u_{i-3,j}^{\ast},\cdots,u_{i+3,j}^{\ast};u_{i,j}^{n};u_{i,j-3}^{\ast},\cdots, u_{i,j+3}^{\ast}),$$
\begin{equation}
\hspace{2in} i=i_{1},\cdots,i_{N}; \quad j=j_{1},\cdots,j_{M}.
\label{eq11}
\end{equation}
$u_{i,j}^{n}$ is the $n$th step iteration value at the grid point $(x_i, y_j)$.  $\gamma$ is the Courant-Friedrichs-Lewy (CFL) number and $\frac{\gamma}{\alpha_{x}/\Delta{x}+\alpha_{y}/\Delta{y}}$ actually corresponds to the time-step size $\Delta{t_{n}}$ in time-marching iterations.
$\alpha_{x}=\max_{u}{|f'(u)|}$ and $\alpha_{y}=\max_{u}{|g'(u)|}$ represent the maximum characteristic speeds in each spatial direction. The nonlinear function $L$, which denotes the fifth-order hybrid AWENO discretization with unequal-sized substencils, depends on $13$ numerical values of the AWENO scheme's stencil. In the fast sweeping method, the Gauss-Seidel iterations and alternating direction sweepings are used.
The Gauss-Seidel philosophy requires that the newest
numerical values of $u$ are used in the fifth-order hybrid AWENO scheme's stencil whenever they are available.
The iterations do not just proceed in only one direction $i = 1:N, j = 1:M$ as in a usual time-marching iteration, but in
the following four alternating directions repeatedly:
$$\mbox{(1) } i=1:N, j=1:M;$$
$$\mbox{(2) } i=N:1, j=1:M;$$
$$\mbox{(3) } i=N:1, j=M:1;$$
$$\mbox{(4) } i=1:N, j=M:1.$$
This is expressed as the notation ``$i=i_{1},\cdots,i_{N};j=j_{1},\cdots,j_{M}$'' in the scheme (\ref{eq11})
to show the alternating sweeping directions in iterations.
Via using alternating direction sweepings, the characteristics property of hyperbolic PDEs is utilized, which is one key factor in fast sweeping methods \cite{ZHAO}. By combining it with the Gauss-Seidel philosophy, a fast convergence to steady state numerical solution of hyperbolic conservation laws is obtained. This will be verified for the proposed hybrid AWENO fast sweeping methods in the following numerical experiments. In the  Gauss-Seidel iterations, the newest numerical values on
the stencil of the hybrid AWENO scheme are used if they are available. This important component of the method is demonstrated by
the notation $u^{\ast}$ in the scheme (\ref{eq11}) to represent the numerical values in the hybrid AWENO stencil, with the understanding that $u_{k,l}^{\ast}$ could be $u_{k,l}^{n}$ or $u_{k,l}^{n+1}$, depending on the current sweeping direction.

\section{Numerical experiments}
\setcounter{equation}{0}
\setcounter{figure}{0}
\setcounter{table}{0}

In this section we perform numerical experiments to test the proposed fifth-order fast sweeping hybrid AWENO method with unequal-sized substencils. To demonstrate the improved absolute convergence of the fast sweeping method, we compare its numerical results with those using the AWENO scheme with equal-sized substencils \cite{AWENO}, for the regular shock reflection problem which is a typical benchmark problem to test the convergence of high order methods in solving steady
flow problems. To verify the improved efficiency of the new fast sweeping method, we compare its results with those obtained by the third-order TVD Runge-Kutta (RK) time-marching method and the numerical results without the hybrid technique. For the convenience of presentation, we abbreviate these schemes as the following: FS-HAUSWENO for the developed fast sweeping hybrid AWENO method with unequal-sized substencils, FS-AUSWENO for the fast sweeping AWENO method with unequal-sized substencils (the hybrid technique is not applied), RK-AWENO for the TVD-RK time-marching using the AWENO scheme with equal-sized substencils, and RK-HAUSWENO for the TVD-RK time-marching using the hybrid AWENO method with unequal-sized substencils.  The convergence of the iterations is measured by the average residue which is defined as
\begin{equation}
ResA=\frac{1}{N\cdot M}\sum_{j=1}^{{M}}\sum_{i=1}^{{N}}|R_{i,j}|,
\end{equation}
where the local residue at the grid point $(x_i,y_j)$ is
$$R_{i,j}=\frac{u_{i,j}^{n+1}-u_{i,j}^{n}}{\Delta t_{n}}.$$
$n$ is the iteration step, and $\Delta t_n=\f{\gamma}{\alpha_x/\Delta x+\alpha_y/\Delta y}$.
For the 2D Euler system of equations, $ResA=\sum_{j=1}^{{M}}\sum_{i=1}^{N}\frac{|R1_{i,j}|+|R2_{i,j}|+|R3_{i,j}|+|R4_{i,j}|}{4\times N\times M}$, where $R\ast_{i,j}$'s are local residuals of these four conservative variables, i.e., $R1_{i,j}=\frac{\rho_{i,j}^{n+1}-\rho_{i,j}^{n}}{\Delta t_{n}}, R2_{i,j}=\frac{(\rho u)_{i,j}^{n+1}-(\rho u)_{i,j}^{n}}{\Delta t_{n}}, R3_{i,j}=\frac{(\rho v)_{i,j}^{n+1}-(\rho v)_{i,j}^{n}}{\Delta t_{n}}, R4_{i,j}=\frac{E_{i,j}^{n+1}-E_{i,j}^{n}}{\Delta t_{n}}$. If not specified otherwise, the convergence criterion is set to $ResA<10^{-12}$. Note that the number of iterations reported in every table here counts a complete update of numerical values in all grid points once as one iteration.

\bigskip
\noindent{\bf Example 1. A Euler system of equations with a smooth solution}

Consider the following two-dimensional Euler system of equations
\[\setlength{\abovedisplayskip}{3pt}\setlength{\belowdisplayskip}{3pt}
\frac{\partial}{\partial{t}}\left(\begin{array}{cccc}
 \rho\\ \rho u\\ \rho v\\E
\end{array}\right)+\frac{\partial}{\partial{x}}
\left(\begin{array}{cccc}
 \rho u\\ \rho u^{2}+p\\ \rho uv\\u(E+p)
\end{array}\right)
+\frac{\partial}{\partial{y}}
\left(\begin{array}{cccc}
 \rho v\\ \rho uv\\ \rho v^{2}+p\\v(E+p)
\end{array}\right)=
0,
\]
with the exact steady state solution $\rho(x,y,\infty)=1+0.2\sin(x-y), u(x,y,\infty)=1,
v(x,y,\infty)=1,$ and $p(x,y,\infty)=1$. The computational domain is $(x,y) \in [0,2\pi] \times [0,2\pi],$ and the exact steady state solution is applied as the boundary conditions in both directions. This problem with a smooth solution is used here to test the accuracy and verify the fifth-order convergence rate of the proposed new fast sweeping method. We take the exact steady state solution as the numerical initial condition to start the iterations. Since
the exact steady state solution does not satisfy the numerical schemes, it will be driven by the iterative schemes and march to their numerical steady states. 
The $L^{1}$ and $L^{\infty}$ numerical errors and orders of accuracy for the density, and the CPU costs of different iterative schemes to converge to their steady states are listed in Table $\ref{1.1}$. The numerical results show that all schemes achieve the desired fifth-order accuracy when they converge to steady state solutions. However, the numerical errors of the schemes with the HLLC flux are much smaller (more than $50\%$ smaller) than those with the LLF flux on a same-size grid. Also, regardless of the choice between the HLLC flux and the LLF flux, the schemes with unequal-sized substencils (FS-HAUSWENO, FS-AUSWENO, RK-HAUSWENO) exhibit much smaller numerical errors than the scheme with equal-sized substencils (RK-AWENO). Note that numerical errors of the hybrid method (FS-HAUSWENO scheme) and the method without hybrid technique (FS-AUSWENO scheme) are similar, which indicates a high accuracy of the AWENO scheme on unequal-sized substencils. In fact, we directly use the fifth-order linear scheme for this problem, and the numerical errors of the linear scheme are also similar to those of both the FS-HAUSWENO and the FS-AUSWENO schemes (the data are not shown to save space here). Furthermore, the CPU costs reported in Table $\ref{1.1}$ show that the FS-HAUSWENO scheme is the most efficient among all four schemes. Comparison of CPU times among the FS-HAUSWENO scheme, the RK-HAUSWENO scheme and the FS-AUSWENO scheme shows that the fast sweeping approach is more efficient than the TVD-RK time-marching, and the hybrid technique further improves efficiency of the fast sweeping iterations. Note that the computational times are relatively short in this example. More obvious CPU-time differences can be observed in the following examples. 
\begin{table}
		\centering
\begin{tabular}{|c|c|c|c|c|c|c|c|c|c|c|}\hline
	\multicolumn{1}{|c|}{}&	\multicolumn{5}{|c|}{RK-AWENO-LLF, $\gamma$=1.0 }&\multicolumn{5}{|c|}{RK-AWENO-HLLC, $\gamma$=1.0 }\\
\hline
 $N\times N$&$L^{1}$ error&order&$L^{\infty}$ error&order&time&$L^{1}$ error&order&$L^{\infty}$ error&order&time\\
 \hline
20$\times$20 & 2.77E-04 &   & 1.38E-03 &  & 0.17 & 1.29E-04 &   & 8.54E-04 &   & 0.25 \\
    \hline
30$\times$30 & 3.88E-05 & 4.85  & 2.24E-04 & 4.49 & 0.42&1.77E-05 & 4.90  & 1.22E-04 & 4.80 & 0.58\\
    \hline
40$\times$40 & 9.32E-06 & 4.96  & 5.70E-05 & 4.76 & 0.89 & 4.24E-06 & 4.96  & 2.81E-05 & 5.10 & 1.06\\
    \hline
50$\times$50 & 3.08E-06 & 4.96  & 1.93E-05 & 4.85 &  1.55 & 1.40E-06 & 4.98  & 9.51E-06 & 4.86 &  1.64\\
    \hline
60$\times$60 & 1.23E-06 & 4.97  & 7.80E-06 & 4.97 &  2.58   & 5.65E-07 & 4.97  & 3.90E-06 & 4.89 &  2.52 \\
    \hline
70$\times$70 & 5.77E-07 & 4.98  & 3.64E-06 & 4.95 & 3.97 & 2.62E-07 & 4.98  & 1.82E-06 & 4.94 & 3.76\\
    \hline
80$\times$80 & 2.97E-07 & 4.96  & 1.88E-06 & 4.94 & 5.67&  1.35E-07 & 4.98  & 9.24E-07 & 5.08 &5.31\\
    \hline
    	\multicolumn{1}{|c|}{}&	\multicolumn{5}{|c|}{RK-HAUSWENO-LLF, $\gamma$=1.0 }&\multicolumn{5}{|c|}{RK-HAUSWENO-HLLC, $\gamma$=1.0 }\\
\hline
 $N\times N$&$L^{1}$ error&order&$L^{\infty}$ error&order&time&$L^{1}$ error&order&$L^{\infty}$ error&order&time\\
 \hline
20$\times$20  & 3.55E-05 &   & 1.93E-04 & &0.13 &1.59E-05 &   & 1.01E-04 & &0.19\\
    \hline
30$\times$30  & 4.94E-06 & 4.86  & 2.82E-05 & 4.74&0.33 &2.22E-06 & 4.86  & 1.43E-05 &4.82 &0.39\\
    \hline
40$\times$40  & 1.20E-06 & 4.92  & 6.42E-06 & 5.14 &0.67& 5.40E-07 & 4.91  & 3.13E-06 & 5.28 &0.77\\
    \hline
50$\times$50  & 3.98E-07 & 4.94  & 2.13E-06 & 4.94&1.19& 1.80E-07 & 4.93  & 1.03E-06 & 4.98&1.33\\
    \hline
60$\times$60   & 1.610E-07 & 4.96  & 8.60E-07 & 4.98&1.94& 7.29E-08 & 4.95  & 4.13E-07 & 5.01&2.13\\
    \hline
70$\times$70  & 7.50E-08 & 4.96  & 3.99E-07 & 4.98 &2.83 & 3.40E-08 & 4.95  & 1.91E-07 & 4.99 &3.16\\
    \hline
80$\times$80 & 3.86E-08 & 4.97  & 2.05E-07 & 4.98&4.40 & 1.75E-08 & 4.96  & 9.83E-08 & 4.99&4.63\\
    \hline
    	\multicolumn{1}{|c|}{}&	\multicolumn{5}{|c|}{FS-AUSWENO-LLF, $\gamma$=1.0 }&\multicolumn{5}{|c|}{FS-AUSWENO-HLLC, $\gamma$=1.0 }\\
\hline
 $N \times N$&$L^{1}$ error&order&$L^{\infty}$ error&order&time&$L^{1}$ error&order&$L^{\infty}$ error&order&time\\
 \hline
20$\times$20  & 3.45E-05 &   & 1.99E-04 & &0.14  &1.54E-05 &   & 1.01E-04 & &0.13\\
    \hline
30$\times$30  & 4.84E-06 & 4.85  & 2.88E-05 & 4.77&0.34 &2.17E-06 & 4.83  & 1.43E-05 &4.82 &0.34\\
    \hline
40$\times$40  & 1.18E-06 & 4.90  & 6.57E-06 & 5.13 &0.69& 5.32E-07 & 4.89  & 3.13E-06 & 5.28 &0.66\\
    \hline
50$\times$50  & 3.93E-07 & 4.93  & 2.18E-06 & 4.95&1.27& 1.78E-07 & 4.92  & 1.03E-06 & 4.98&1.16\\
    \hline
60$\times$60   & 1.60E-07 & 4.94  & 8.79E-07 & 4.98&2.02& 7.22E-08 & 4.93  & 4.13E-07 & 5.01&1.94\\
    \hline
70$\times$70  & 7.45E-08 & 4.95  & 4.08E-07 & 4.98 &3.13 & 3.37E-08 & 4.94  & 1.91E-07 & 4.99 &2.97\\
    \hline
80$\times$80 & 3.84E-08 & 4.96  & 2.10E-07 & 4.98&4.63 & 1.74E-08 & 4.95  & 9.83E-08 & 4.99&4.67\\
    \hline
	\multicolumn{1}{|c|}{}&	\multicolumn{5}{|c|}{FS-HAUSWENO-LLF, $\gamma$=1.0 }&\multicolumn{5}{|c|}{FS-HAUSWENO-HLLC, $\gamma$=1.0 }\\
\hline
 $N\times N$&$L^{1}$ error&order&$L^{\infty}$ error&order&time&$L^{1}$ error&order&$L^{\infty}$ error&order&time\\
 \hline
20$\times$20  & 3.45E-05 &   & 1.99E-04 & &0.11  &1.54E-05 &   & 1.01E-04 & &0.13\\
    \hline
30$\times$30  & 4.84E-06 & 4.85  & 2.87E-05 & 4.77&0.28 &2.17E-06 & 4.83  & 1.43E-05 &4.81 &0.27\\
    \hline
40$\times$40  & 1.18E-06 & 4.90  & 6.56E-06 & 5.13 &0.59& 5.32E-07 & 4.89  & 3.13E-06 & 5.28 &0.52\\
    \hline
50$\times$50  & 3.93E-07 & 4.93  & 2.17E-06 & 4.95&1.00& 1.78E-07 & 4.92  & 1.03E-06 & 4.98&0.92\\
    \hline
60$\times$60   & 1.60E-07 & 4.94  & 8.78E-07 & 4.97&1.70& 7.22E-08 & 4.93  & 4.13E-07 & 5.01&1.53\\
    \hline
70$\times$70  & 7.45E-08 & 4.95  & 4.08E-07 & 4.98 &2.41 & 3.37E-08 & 4.94  & 1.91E-07 & 4.99 &2.25\\
    \hline
80$\times$80 & 3.84E-08 & 4.96  & 2.10E-07 & 4.98&3.50 & 1.74E-08 & 4.95  & 9.84E-08 & 4.98&3.33\\
    \hline
		\end{tabular}\caption{\label{1.1}Example 1. Euler system of equations with a smooth solution. $L^{1}$ and $L^{\infty}$ numerical errors, orders of accuracy, and CPU times of different iterative schemes. CPU time unit: second.}
	\end{table}

\bigskip
\noindent{\bf Example 2. Regular shock reflection}

\begin{figure}
\centering
\subfigure[RK-AWENO-LLF]{
\begin{minipage}[t]{0.5\linewidth}
\centering
\includegraphics[width=3.1in]{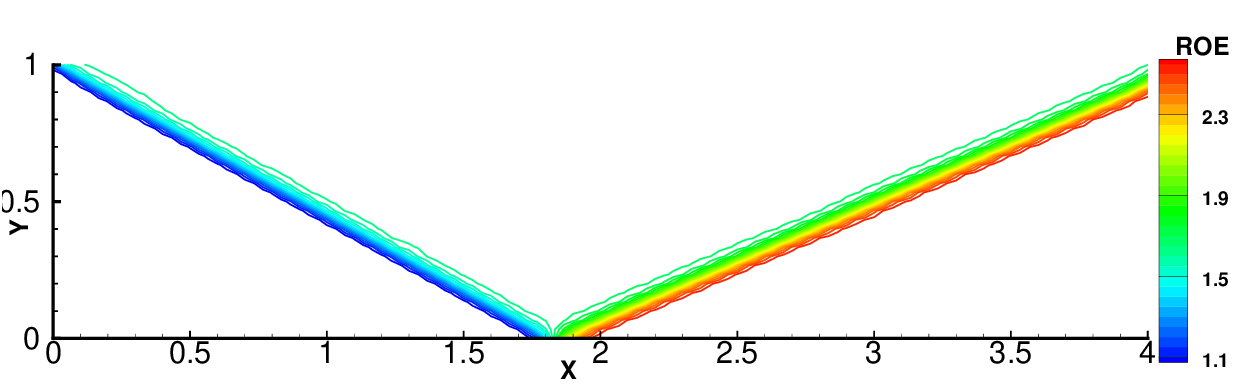}
\end{minipage}%
}%
\subfigure[RK-AWENO-HLLC]{
\begin{minipage}[t]{0.5\linewidth}
\centering
\includegraphics[width=3.1in]{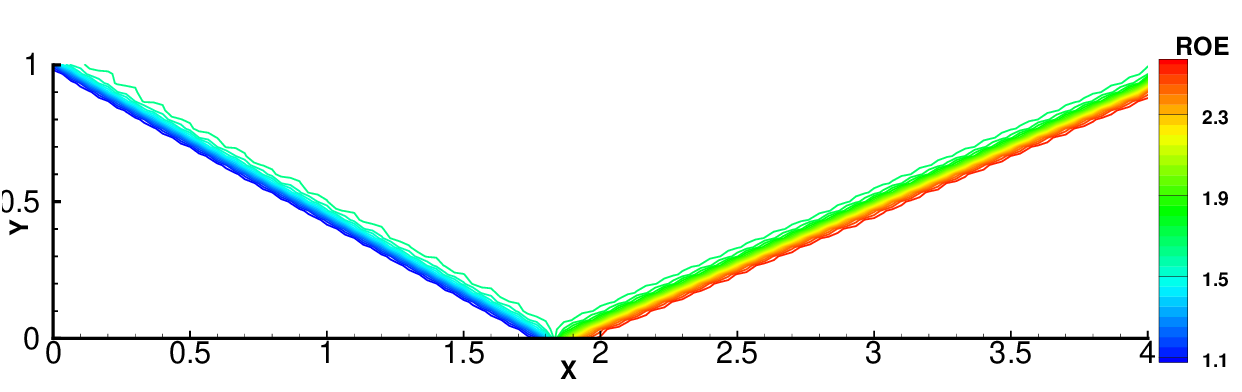}
\end{minipage}%
}%

\subfigure[RK-HAUSWENO-LLF]{
\begin{minipage}[t]{0.5\linewidth}
\centering
\includegraphics[width=3.1in]{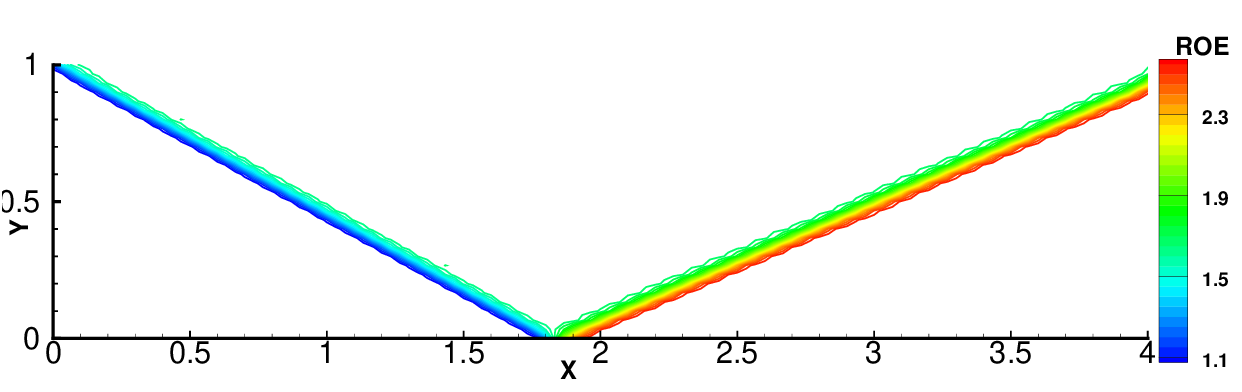}
\end{minipage}%
}%
\subfigure[RK-HAUSWENO-HLLC]{
\begin{minipage}[t]{0.5\linewidth}
\centering
\includegraphics[width=3.1in]{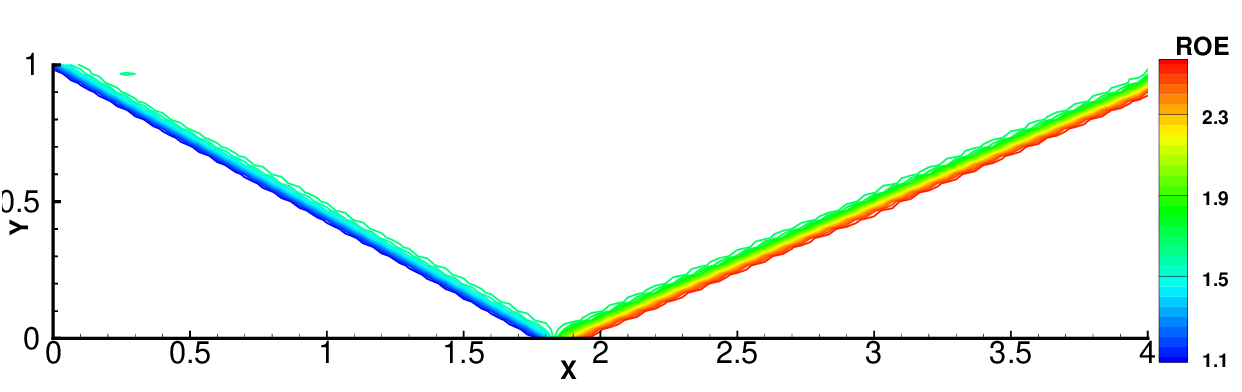}
\end{minipage}%
}%

\subfigure[FS-HAUSWENO-LLF]{
\begin{minipage}[t]{0.5\linewidth}
\centering
\includegraphics[width=3.1in]{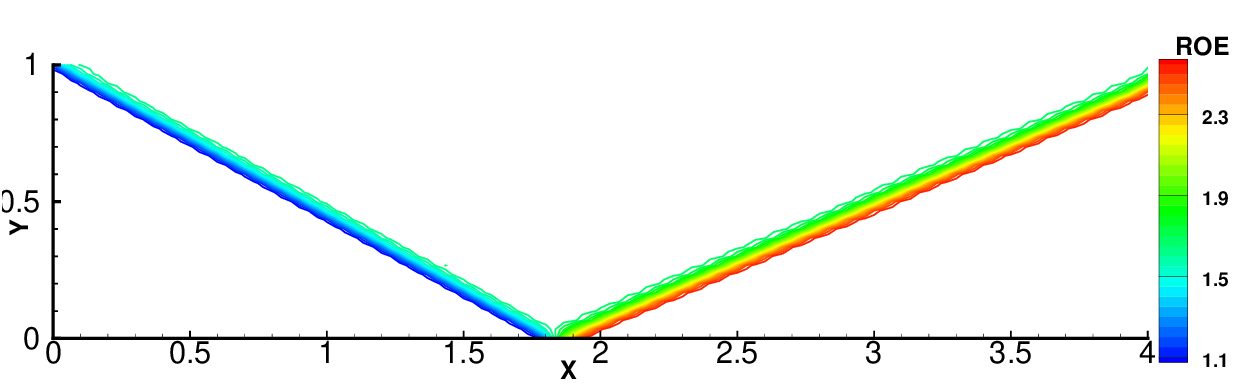}
\end{minipage}%
}
\subfigure[FS-HAUSWENO-HLLC]{
\begin{minipage}[t]{0.5\linewidth}
\centering
\includegraphics[width=3.1in]{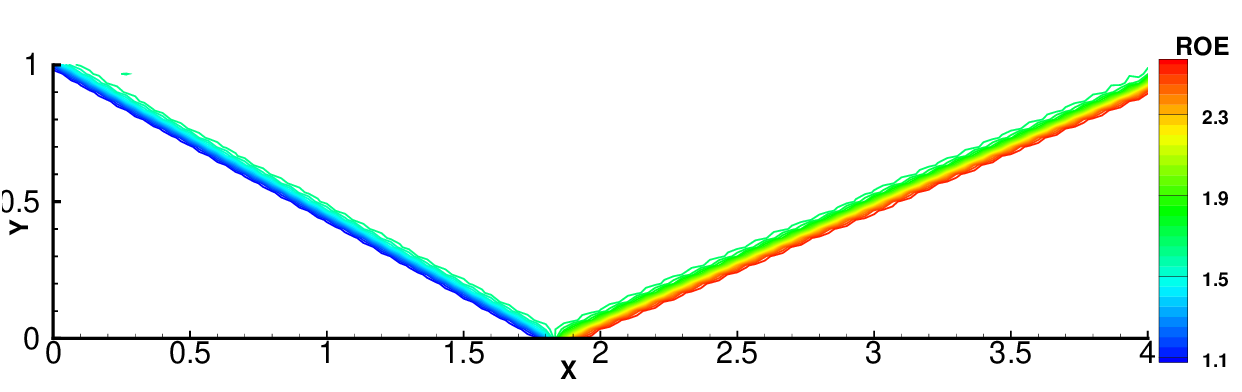}
\end{minipage}%
}%

\subfigure[FS-HAUSWENO-LLF]{
\begin{minipage}[t]{0.5\linewidth}
\centering
\includegraphics[width=3.1in]{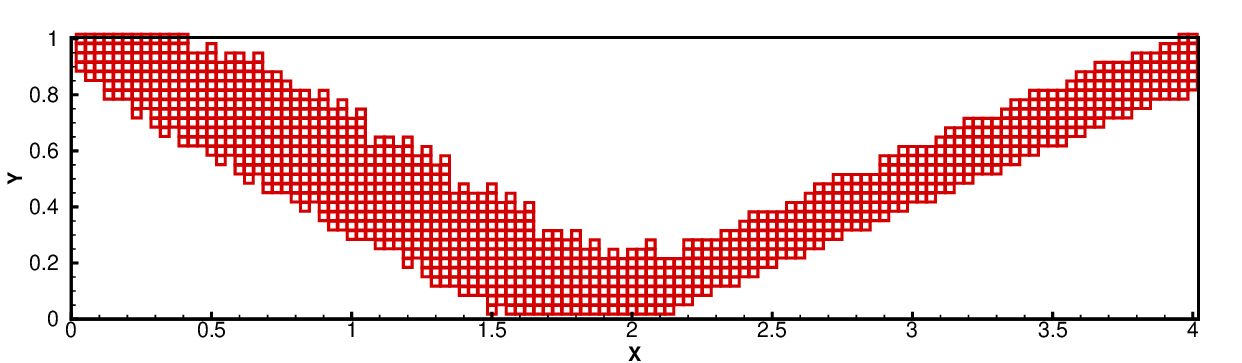}
\end{minipage}%
}
\subfigure[FS-HAUSWENO-HLLC]{
\begin{minipage}[t]{0.5\linewidth}
\centering
\includegraphics[width=3.1in]{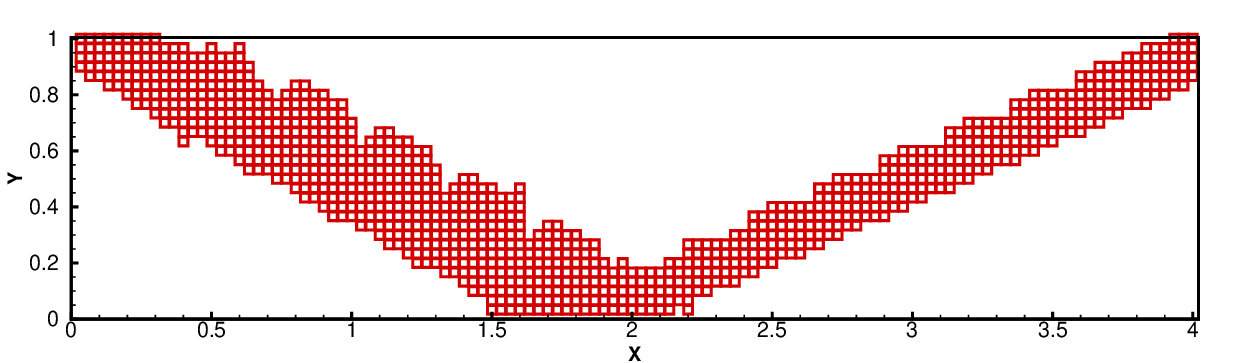}
\end{minipage}%
}%
\centering
\caption{\label{2.2}Example 2: Regular shock reflection. 30 equally spaced density contours from 1.1 to 2.6 of the steady states of numerical solutions by different iterative schemes, and the identified troubled-cells where WENO interpolation is used in the FS-HAUSWENO scheme.}
\end{figure}

\begin{figure}
\centering
\subfigure[RK-AWENO-LLF]{
\begin{minipage}[t]{0.4\linewidth}
\centering
\includegraphics[width=2.8in]{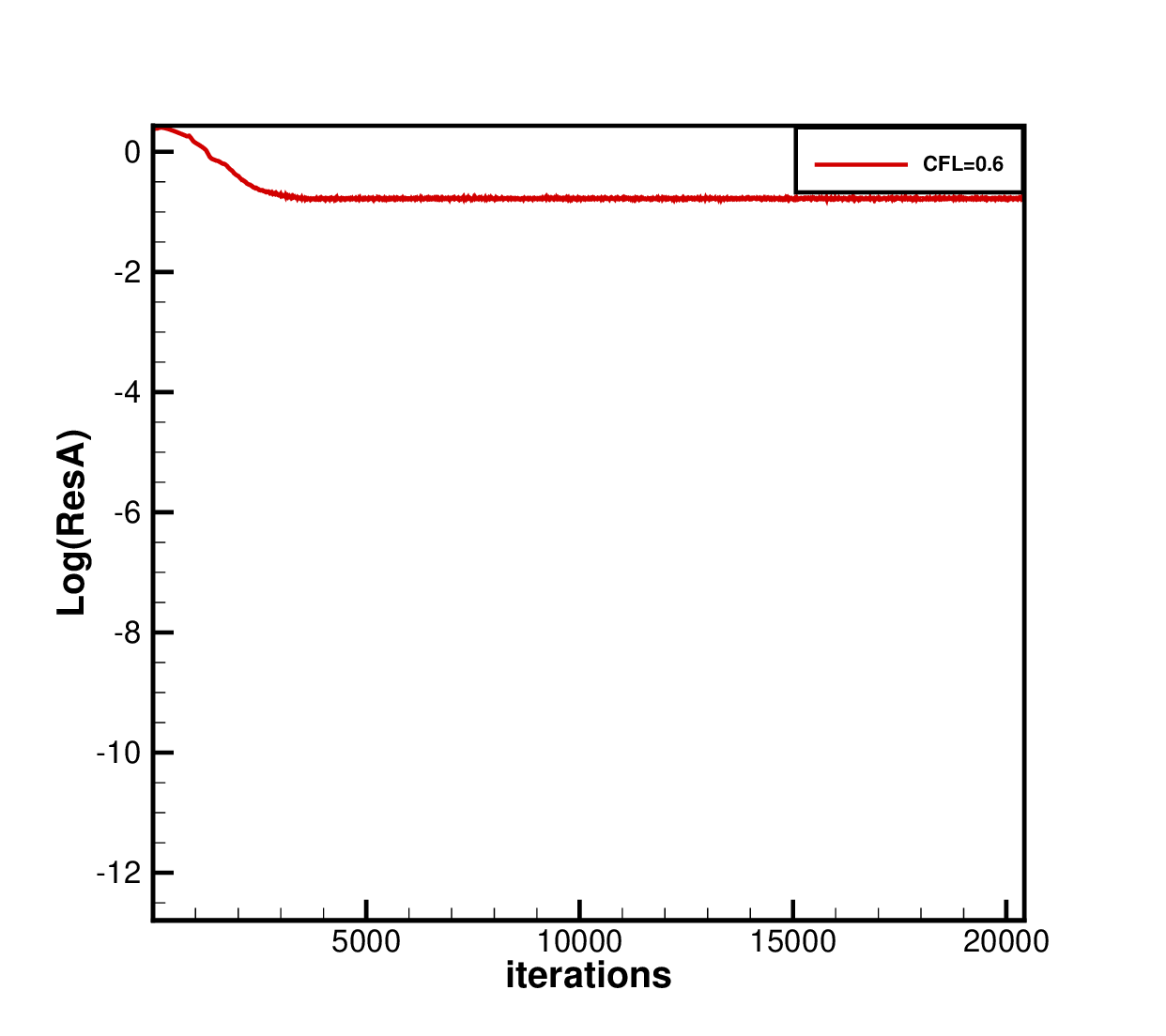}
\end{minipage}%
}%
\subfigure[RK-AWENO-HLLC]{
\begin{minipage}[t]{0.4\linewidth}
\centering
\includegraphics[width=2.8in]{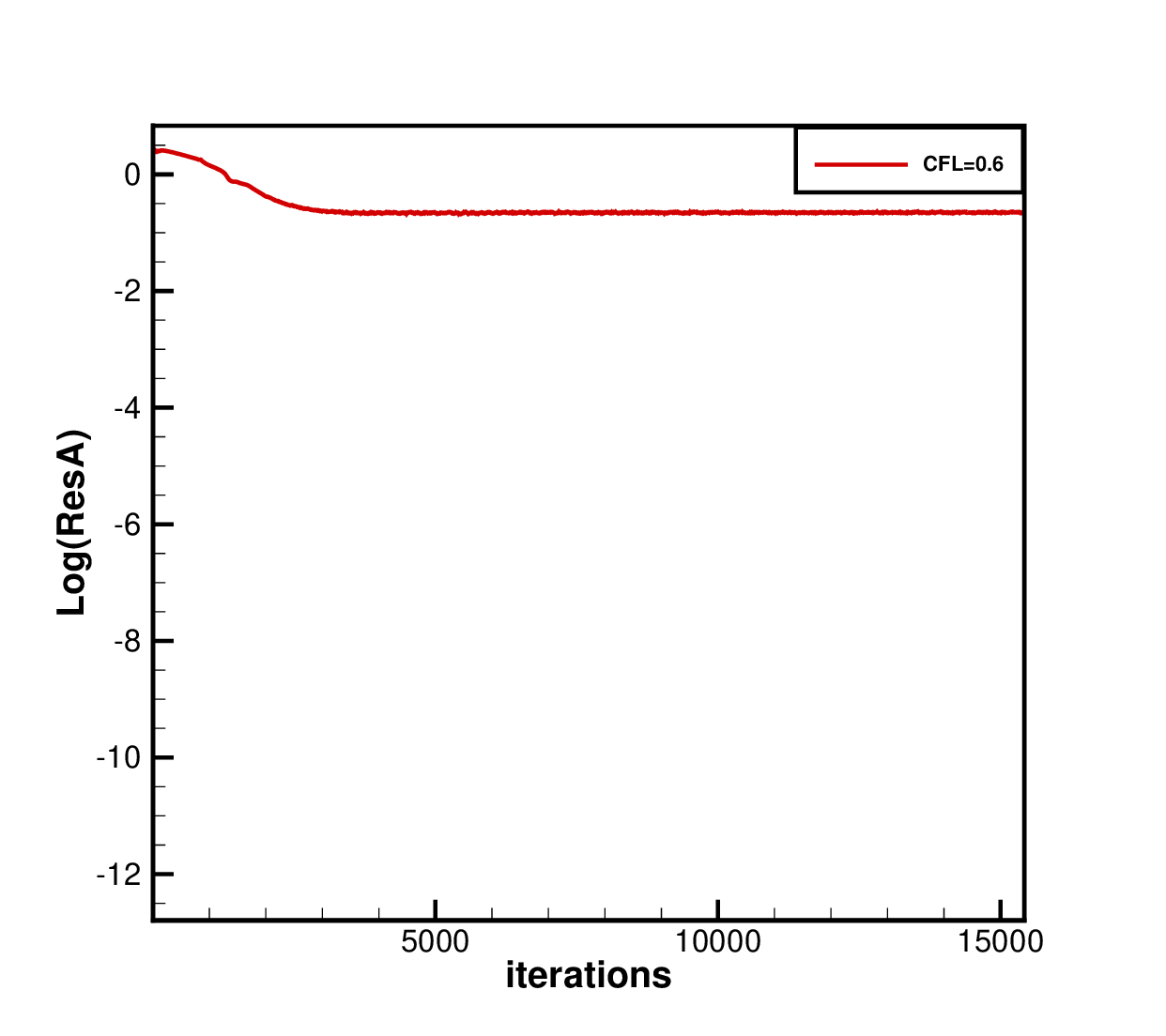}
\end{minipage}%
}%

\subfigure[RK-HAUSWENO-LLF]{
\begin{minipage}[t]{0.4\linewidth}
\centering
\includegraphics[width=2.8in]{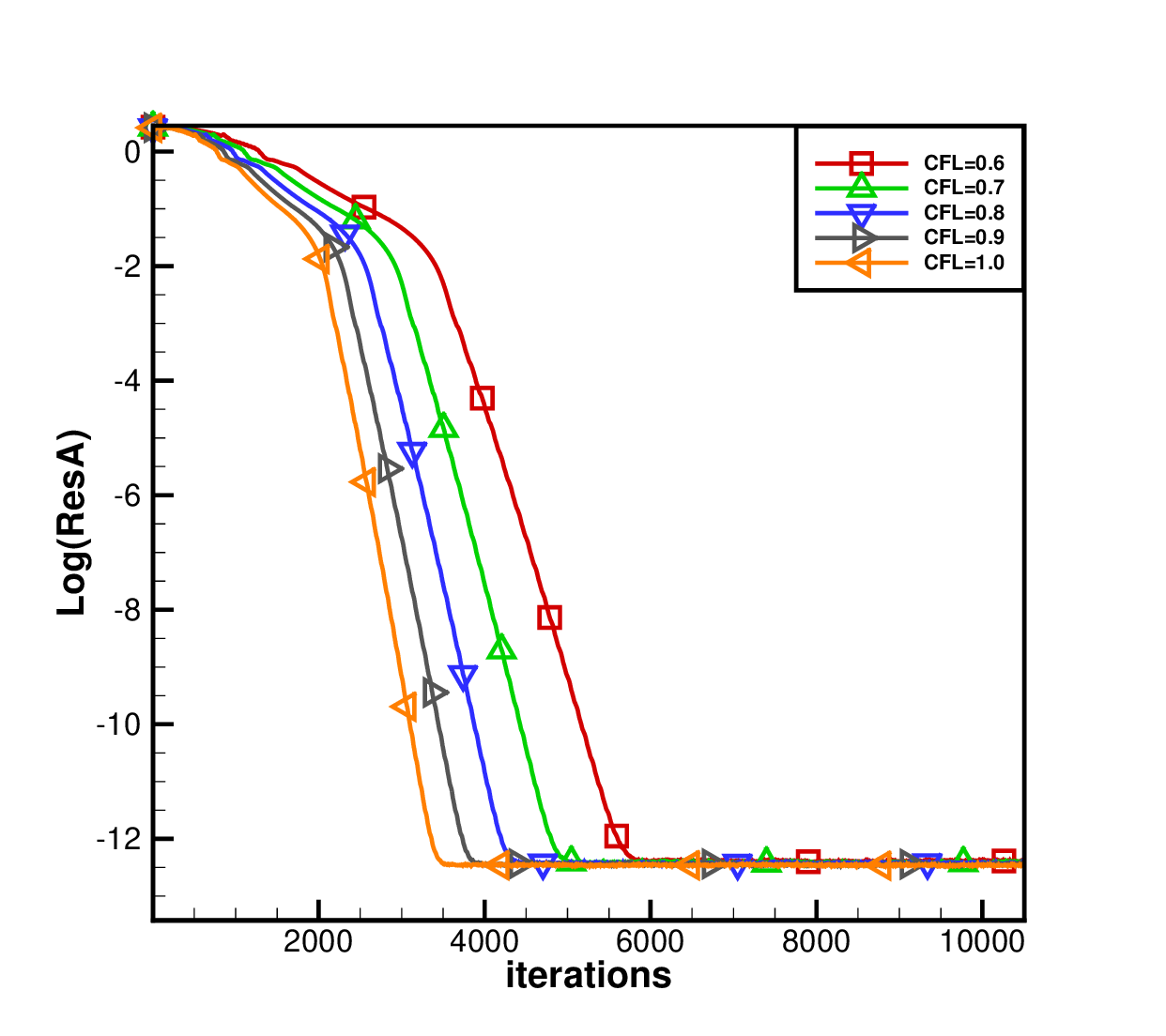}
\end{minipage}%
}%
\subfigure[RK-HAUSWENO-HLLC]{
\begin{minipage}[t]{0.4\linewidth}
\centering
\includegraphics[width=2.8in]{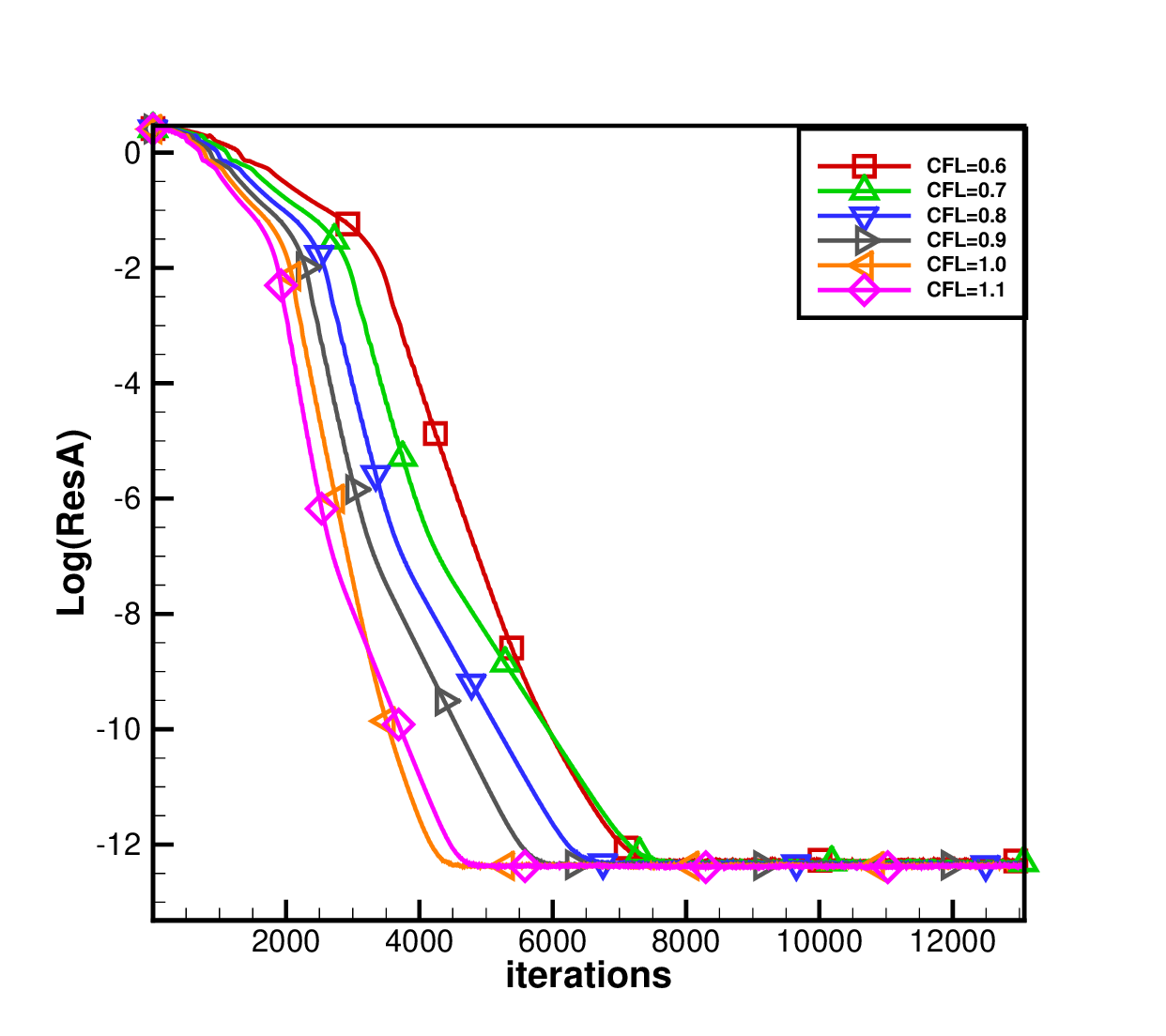}
\end{minipage}%
}%

\subfigure[FS-HAUSWENO-LLF]{
\begin{minipage}[t]{0.4\linewidth}
\centering
\includegraphics[width=2.8in]{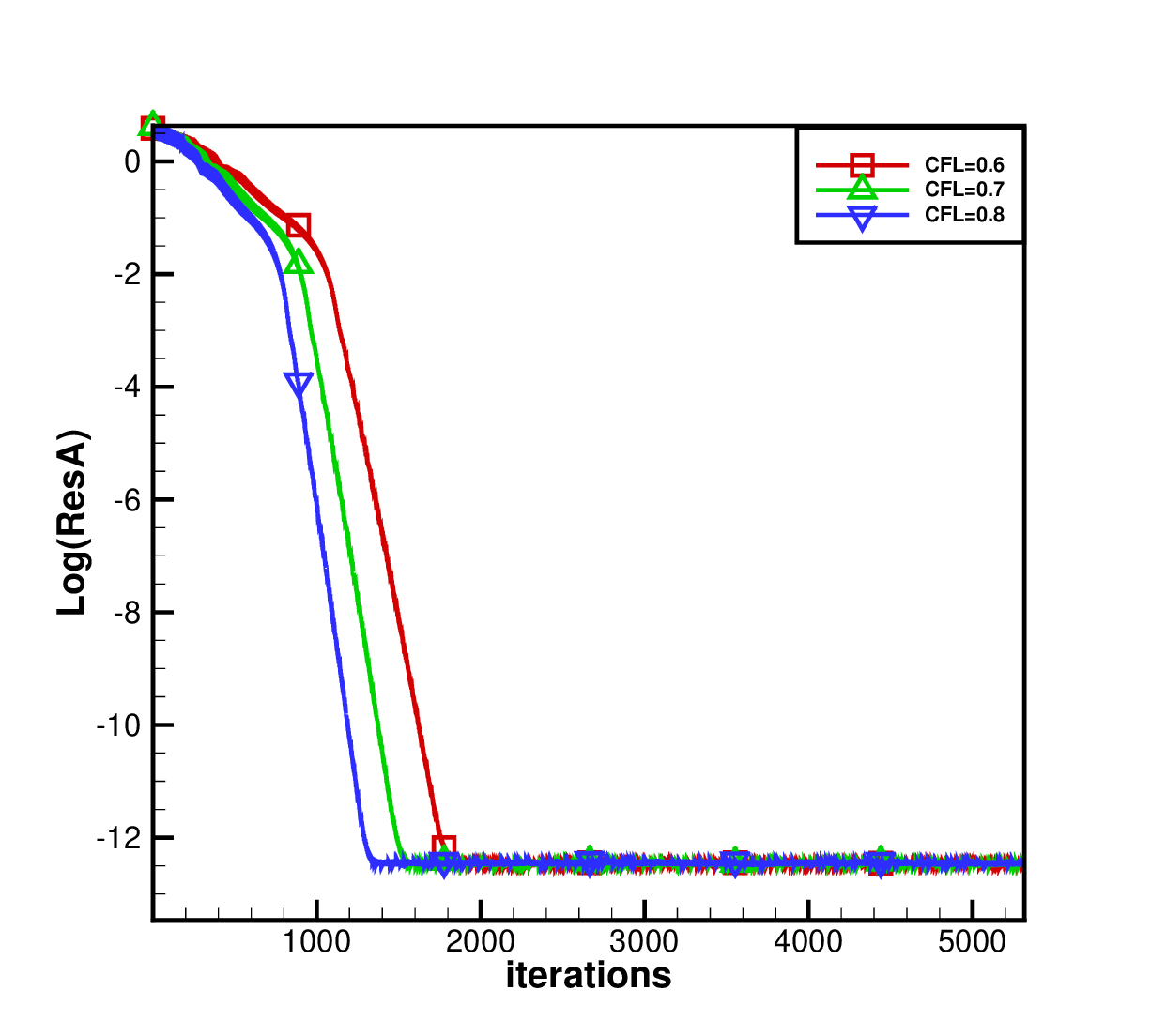}
\end{minipage}%
}
\subfigure[FS-HAUSWENO-HLLC]{
\begin{minipage}[t]{0.4\linewidth}
\centering
\includegraphics[width=2.8in]{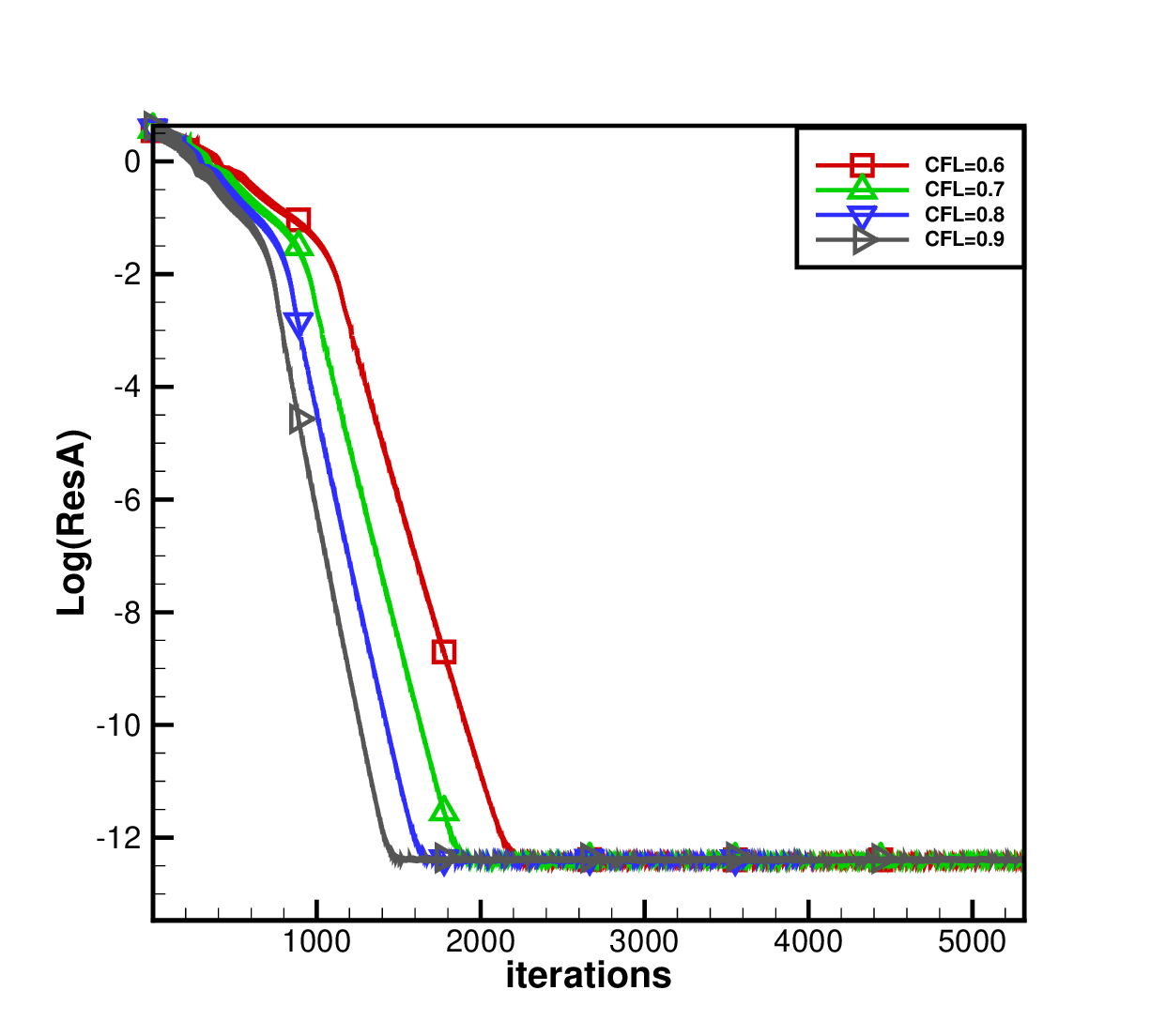}
\end{minipage}%
}%
\centering
\caption{\label{2.3}Example 2: Regular shock reflection. The convergence history of the residue as a function of number of iterations for different schemes with various CFL numbers.}
\end{figure}

In this example, the regular shock reflection problem is solved to test the performance of the FS-HAUSWENO scheme. This is a typical benchmark problem to test the convergence of high order numerical methods for solving 2D steady-flow problems. The computational domain is $[0,4]\times[0,1]$. The boundary conditions are that of a reflection condition along the bottom boundary, supersonic outflow along the right boundary, and Dirichlet conditions on the other two sides:
\begin{equation*}
(\rho,u,v,p)^{T}=
\begin{cases}
(1.0,2.9,0,5/7)^{T}\mid_{(0,y,t)^{T}},\\
(1.69997,2.61934,-0.50632,1.52819)^{T}\mid_{(x,1,t)^{T}}.\\
\end{cases}
\end{equation*}
Initially, we set the solution for the entire domain to be the values at the left boundary to start iterations of the schemes. For this steady-flow problem, because it is challenging for the residues of many high-order schemes including the WENO schemes \cite{WuLiang,SSCW,SCW} to converge to machine zero or round off errors, it is used here to verify the absolute convergence of the proposed hybrid fast sweeping method. The computational grid is $120\times30$, and these different iterative schemes are applied to this problem for comparison. The contour plots of the density variable of the numerical steady states for all of these schemes are comparable and are presented in 
Fig.~\ref{2.2}. Here we omit the similar plots of the FS-AUSWENO schemes in this example and also the following examples, to save space. All schemes accurately capture shock waves and are stable. The cells where the WENO interpolation is used in the FS-HAUSWENO scheme are also shown in  
Fig.~\ref{2.2}. It is observed that the hybrid technique identifies the region of troubled-cells well. However, the residue of the RK-AWENO scheme with both the LLF and the HLLC fluxes fails to converge to round off errors, as shown in Fig. $\ref{2.3}$ (a) and (b).  Hence, the RK-AWENO scheme does not have the property of absolute convergence. Fig. $\ref{2.3}$ shows that the other schemes based on unequal-sized substencils, including the proposed FS-HAUSWENO scheme, achieve absolute convergence regardless of the LLF flux or the HLLC flux (the plots for the FS-AUSWENO scheme are omitted to save space). We would like to emphasize that the novel hybrid technique proposed in this paper preserves the absolute convergence for the alternative WENO scheme with unequal-sized substencils, as shown in Fig. $\ref{2.3}$, while the hybrid scheme directly based on the high-order polynomial $p_{1}(x)$ does not preserve the absolute convergence for this typical example (see Fig. 19 in \cite{USA}). This comparison suggests that the hybrid scheme here provides a robust and smooth transition between the high-order interpolation of linear scheme and WENO interpolation. To evaluate computational efficiency,  we perform tests with various CFL numbers to compare these iterative schemes. The CFL numbers are gradually increased to the maximum value that permits the convergence of iterations. The number of iterations, the final time, and the total CPU time when the schemes converge are reported for the TVD-RK time-marching method and the fast sweeping method in Table \ref{2.1}. Note that when the CFL number is increased to a value that makes the iterations not converge, we denote it by ``not conv'' in the table. It is observed that the largest CFL numbers permitted in these schemes for convergence are comparable. However, the fast sweeping method takes fewer iterations to converge to steady state than the TVD-RK method, and saves approximately $40\%-60\%$ of CPU time cost with a comparison of FS-HAUSWENO and RK-HAUSWENO under the same CFL number. Also, it is verified in Table \ref{2.1} that the hybrid fast sweeping method (FS-HAUSWENO) is more efficient than the one without the hybrid technique (FS-AUSWENO). It is interesting to observe that the savings of computational costs via using the hybrid technique here is relatively less than those in the other hybrid finite difference WENO schemes based on the finite volume reconstruction approach and flux splitting (e.g. \cite{ZQH}), due to the relatively smaller percentages of the computation times which the nonlinear weights take in the AWENO schemes. This is also consistent with the results on the computational times of the AWENO schemes in \cite{USA}. 
\begin{table}
		\centering
\begin{tabular}{|c|c|c|c|c|c|c|}\hline
		&\multicolumn{3}{|c|}{RK-HAUSWENO-LLF}&\multicolumn{3}{|c|}{RK-HAUSWENO-HLLC}\\\hline
            CFL  & iteration number & final time & CPU time & iteration number & final time & CPU time \\\hline
0.6	& 5607	&4.73	&21.33 	& 7038	&5.96	&30.81                   \\                           \hline
0.7	&4806	&4.73	&18.58 	&7140	&7.06	&30.83                   \\                           \hline
0.8	&4203	&4.73	&16.02 	&6225	&7.04	&28.33                  \\                           \hline
0.9	&3735	&4.73	&14.05 	&5514	&7.01	&22.63                  \\                           \hline
1.0	&3363	&4.73	&13.05 	&4179	&5.90	&17.41                 \\                           \hline
1.1	&not conv	&	&          &4464	&6.94	&18.61                         \\                                \hline
1.2	&not conv	&	&     	    &not conv	&	&                    \\
		\end{tabular}
\begin{tabular}{|c|c|c|c|c|c|c|}\hline
		&\multicolumn{3}{|c|}{FS-AUSWENO-LLF}&\multicolumn{3}{|c|}{FS-AUSWENO-HLLC}\\\hline
            CFL  & iteration number & final time & CPU time & iteration number & final time & CPU time \\\hline
0.6	& 1755	&4.44	&16.92 	& 2135	&5.41	&18.08                  \\                           \hline
0.7	&1487	&4.39	&13.34 	&1827	&5.40	&15.39                   \\                           \hline
0.8	&1290	&4.35	&11.47 	&1598	&5.40	&13.45                  \\                           \hline
0.9	&not conv	&	&          &1414	&5.37	&12.19                         \\ \hline
1.0	&not conv	&	&          &not conv	&	&                         \\
		\end{tabular}
\begin{tabular}{|c|c|c|c|c|c|c|}\hline
		&\multicolumn{3}{|c|}{FS-HAUSWENO-LLF}&\multicolumn{3}{|c|}{FS-HAUSWENO-HLLC}\\\hline
            CFL  & iteration number & final time & CPU time & iteration number & final time & CPU time \\\hline
0.6	&1755	&4.42	&13.25 	&2135	&5.41	&15.09                  \\                           \hline
0.7	&1487	&4.39	&11.33 	&1827	&5.40	&12.88                   \\                           \hline
0.8	&1290	&4.35	&9.86 	&1598	&5.40	&11.36                  \\                           \hline
0.9	&not conv	&	&          &1414	&5.37	&10.42                         \\                                \hline
1.0	&not conv	&	&          &not conv	&	&                         \\                                \hline

		\end{tabular}
\caption{\label{2.1}Example 2: Regular shock reflection. Number of iterations, the final time and total CPU time when convergence is obtained. Convergence criterion is $10^{-12}$. CPU time unit: second.}
\end{table}

\bigskip
\noindent{\bf Example 3. Supersonic flow past a plate with an attack angle}

\begin{figure}
\centering
\subfigure[RK-HAUSWENO-LLF]{
\begin{minipage}[t]{0.45\linewidth}
\centering
\includegraphics[width=2.8in]{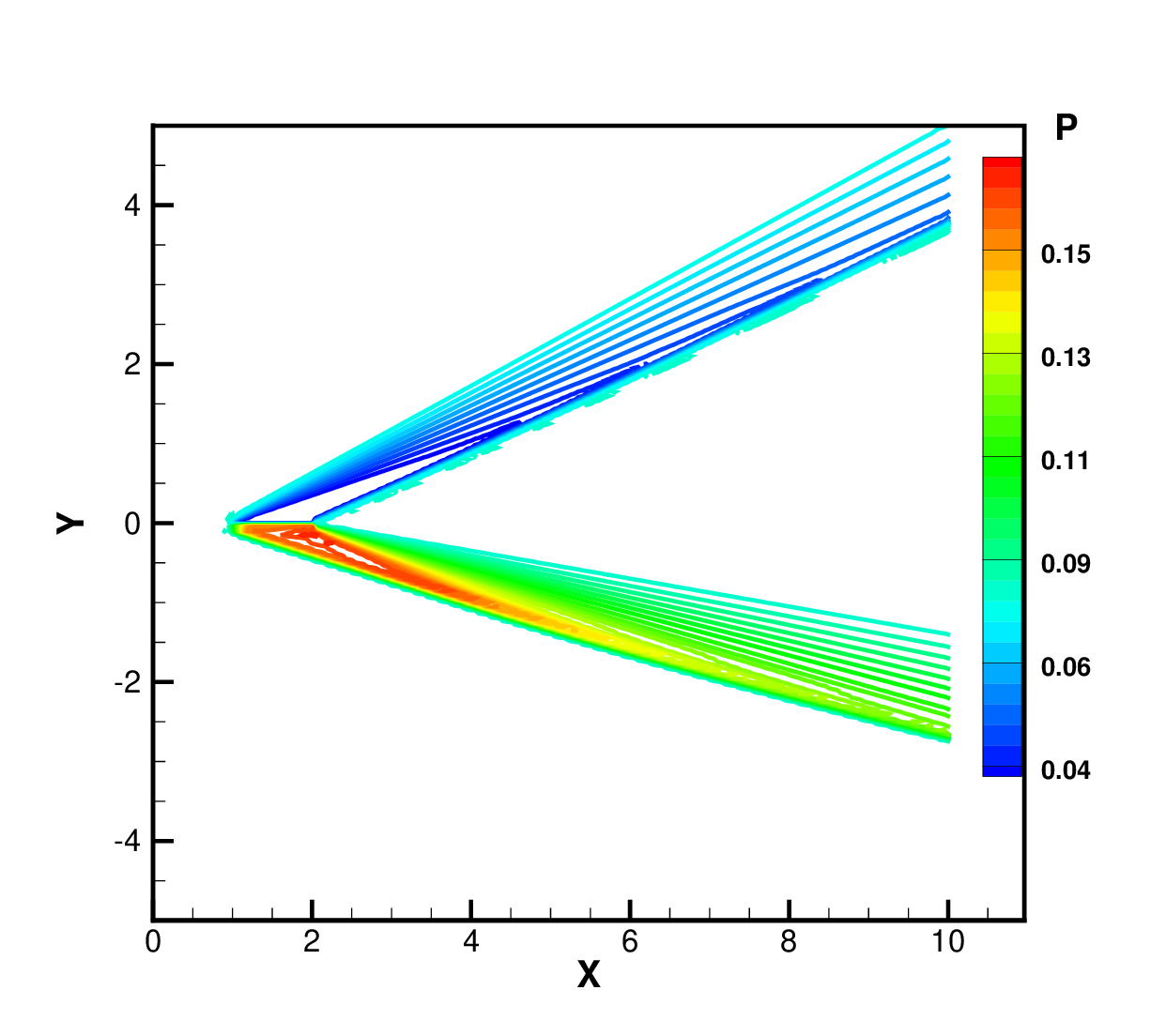}
\end{minipage}%
}%
\subfigure[RK-HAUSWENO-HLLC]{
\begin{minipage}[t]{0.45\linewidth}
\centering
\includegraphics[width=2.8in]{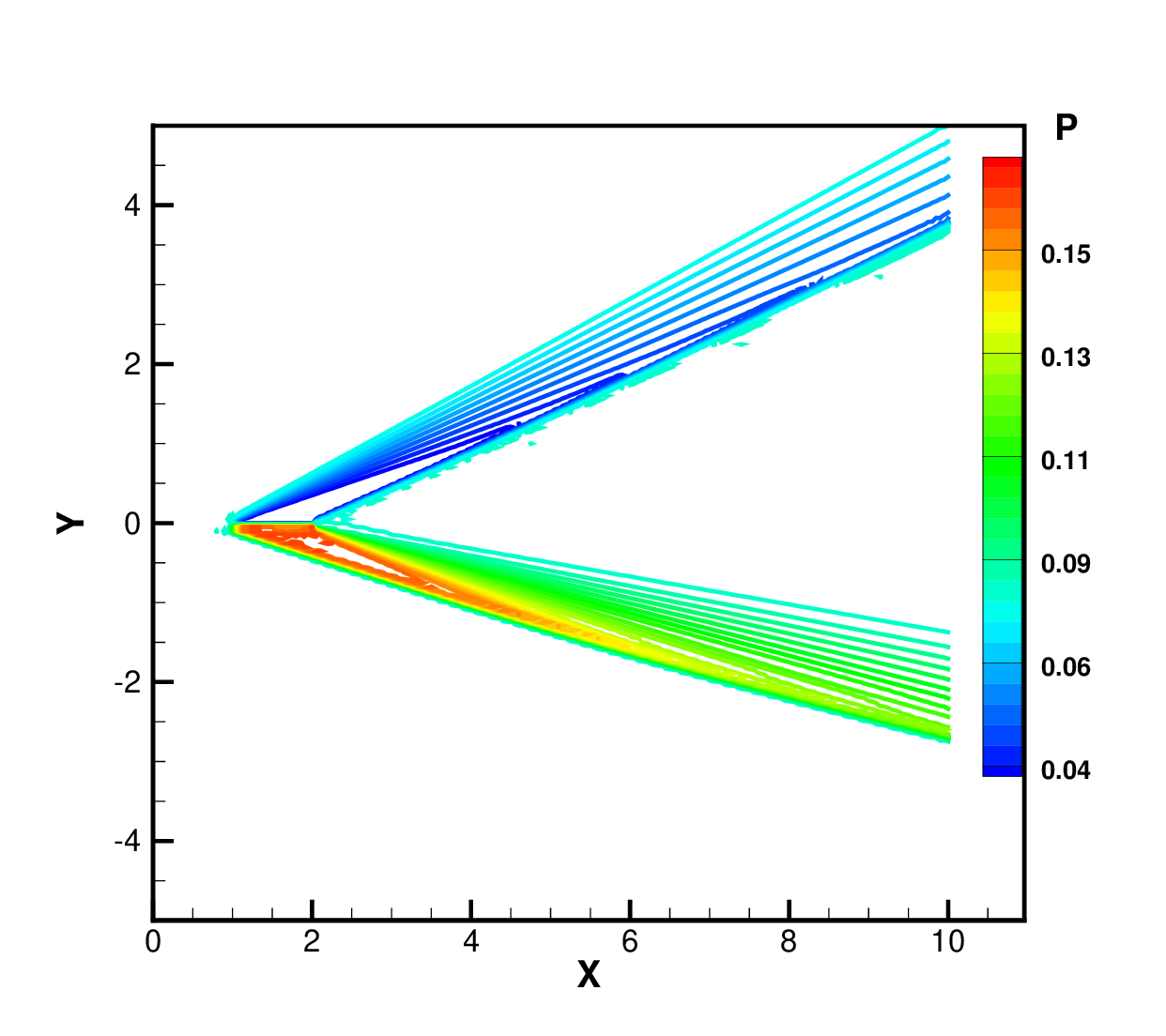}
\end{minipage}%
}%

\subfigure[FS-HAUSWENO-LLF]{
\begin{minipage}[t]{0.45\linewidth}
\centering
\includegraphics[width=2.8in]{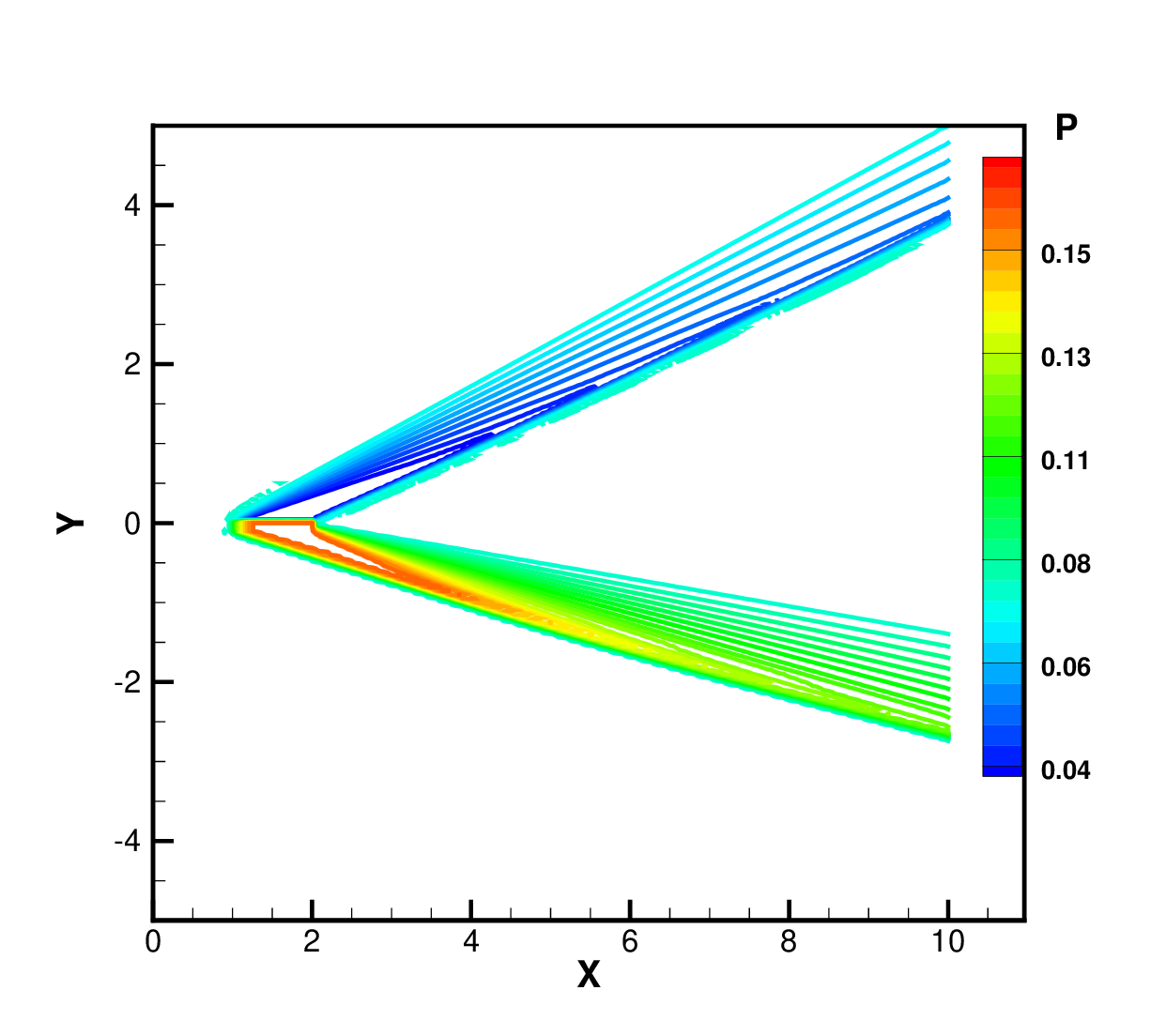}
\end{minipage}%
}
\subfigure[FS-HAUSWENO-HLLC]{
\begin{minipage}[t]{0.42\linewidth}
\centering
\includegraphics[width=2.8in]{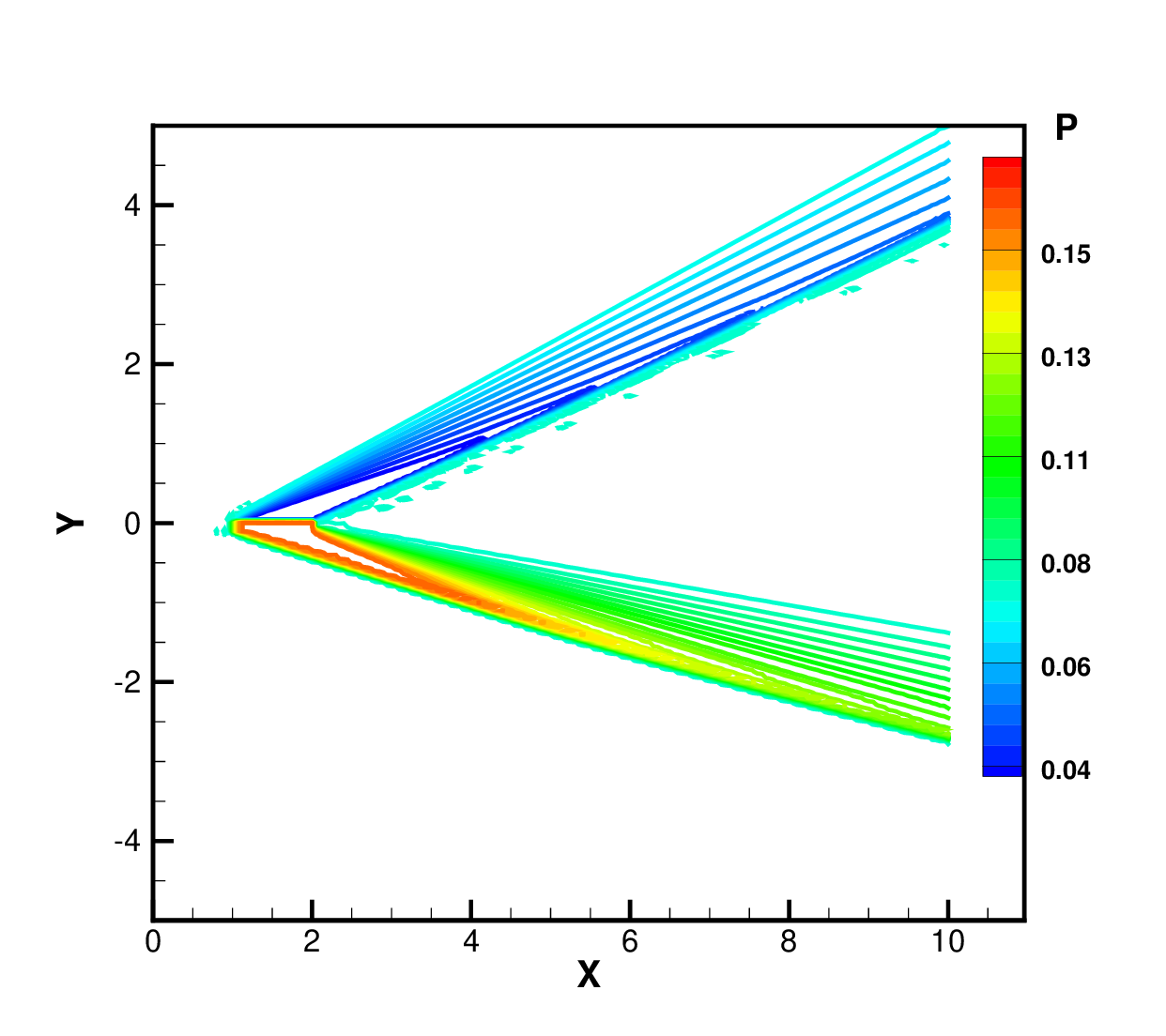}
\end{minipage}%
}%

\subfigure[FS-HAUSWENO-LLF]{
\begin{minipage}[t]{0.42\linewidth}
\centering
\includegraphics[width=2.8in]{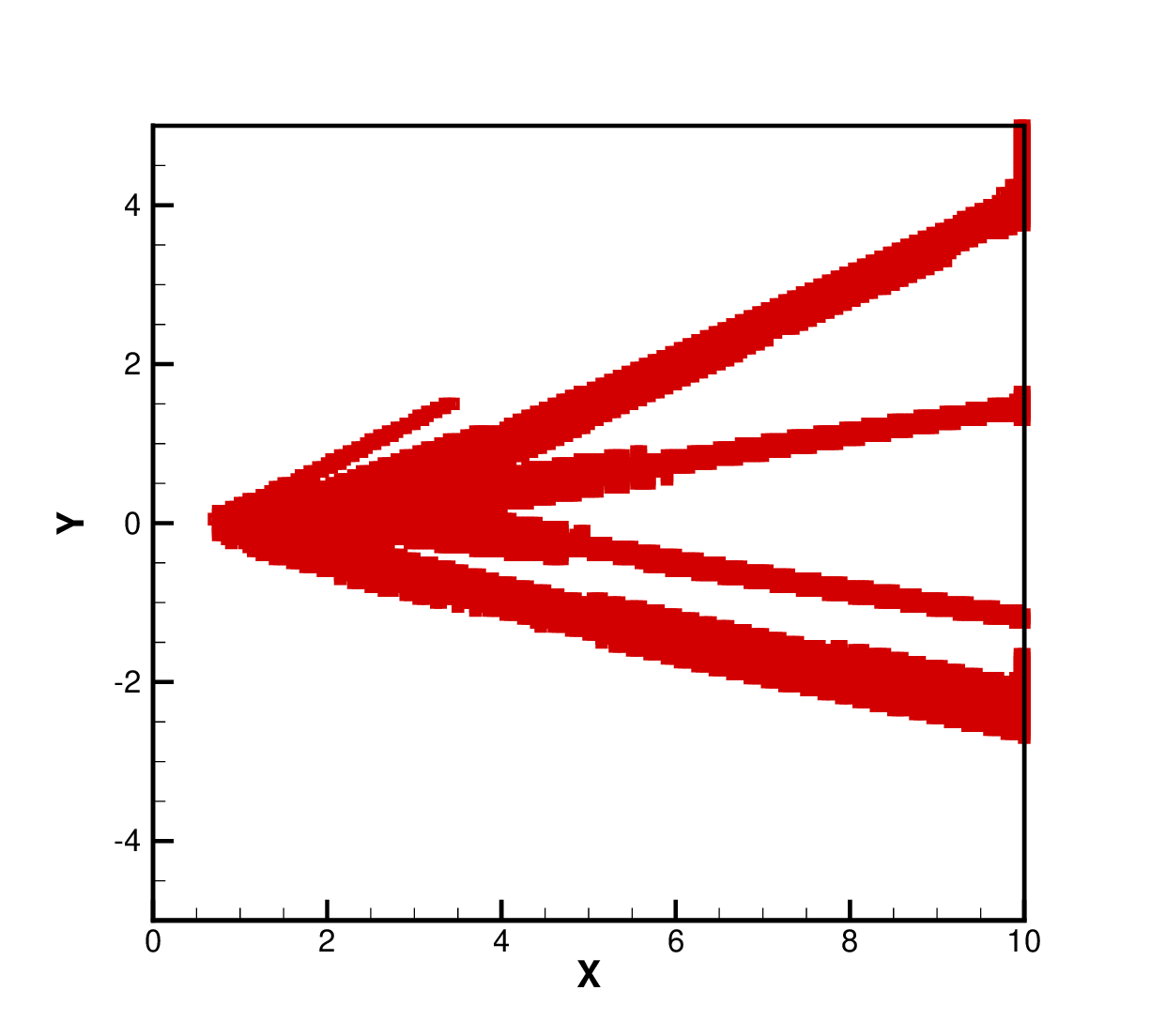}
\end{minipage}%
}
\subfigure[FS-HAUSWENO-HLLC]{
\begin{minipage}[t]{0.45\linewidth}
\centering
\includegraphics[width=2.8in]{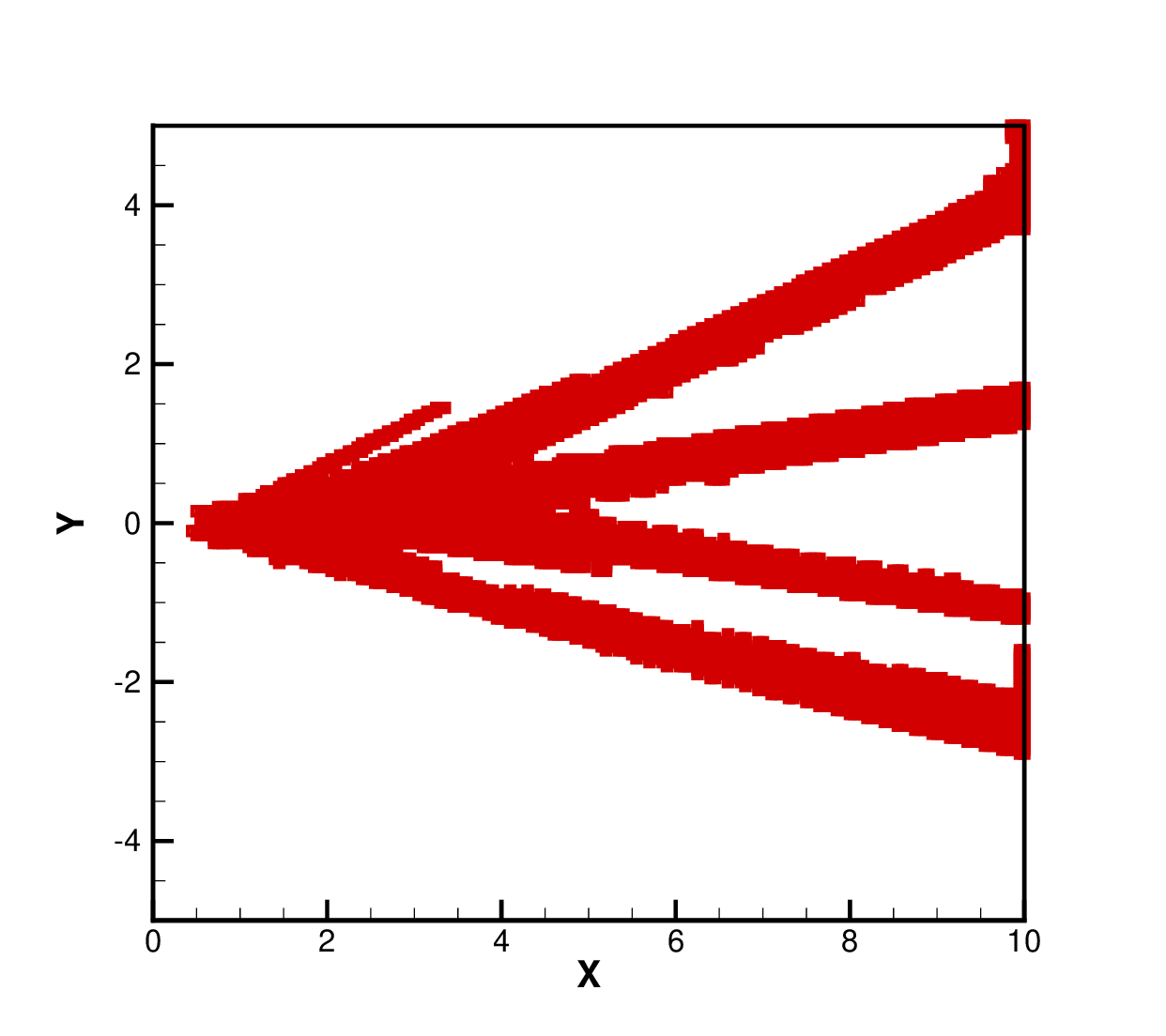}
\end{minipage}%
}%
\centering
\caption{\label{3.2}Example 3: Supersonic flow past a plate with an attack angle. 30 equally spaced pressure contour from 0.04 to 0.17 of the converged steady states of numerical solutions by different iterative schemes, and the identified troubled-cells where the WENO interpolation is used in the FS-HAUSWENO scheme.}
\end{figure}

\begin{figure}
\centering
\subfigure[RK-HAUSWENO-LLF]{
\begin{minipage}[t]{0.45\linewidth}
\centering
\includegraphics[width=3.0in]{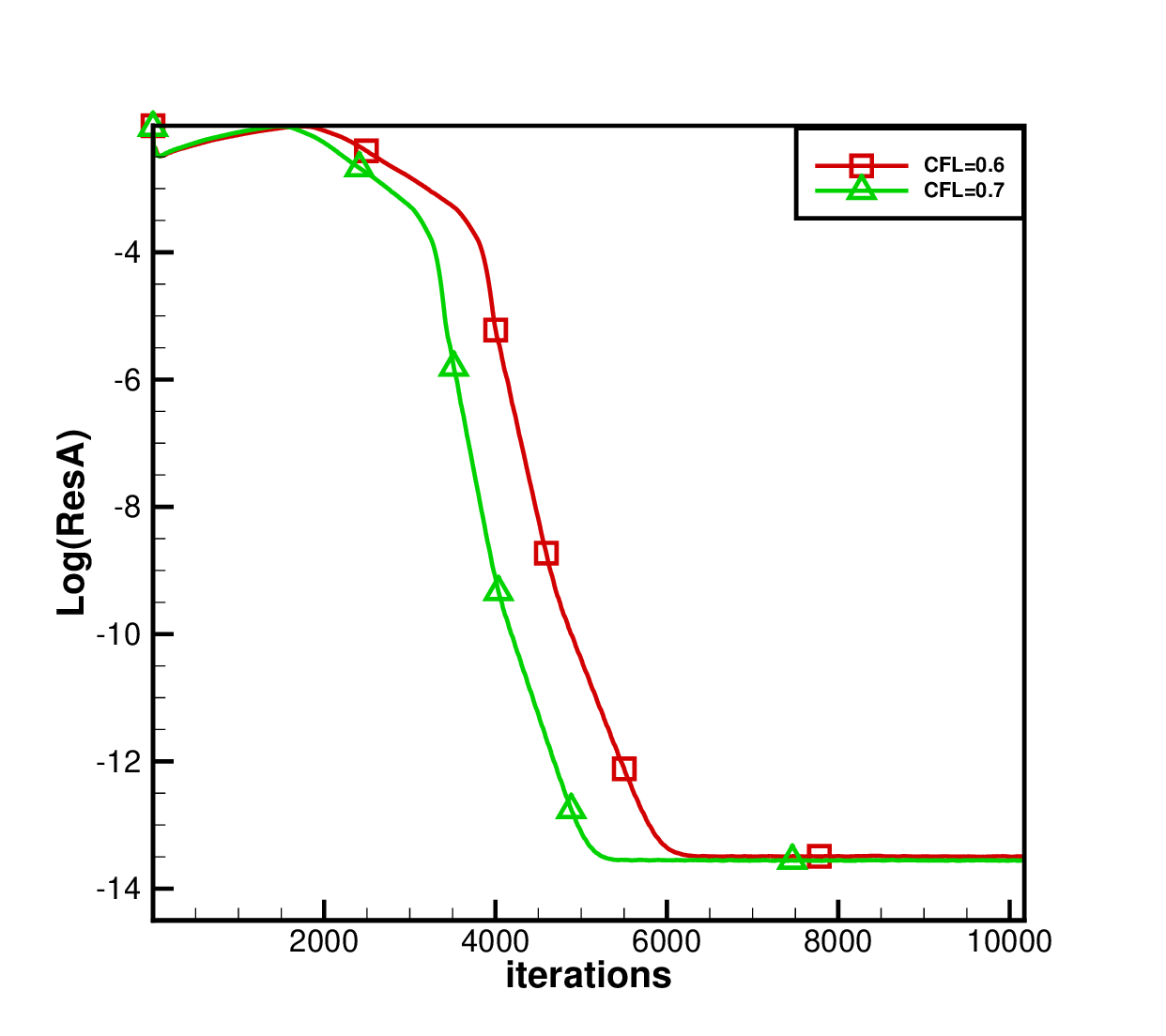}
\end{minipage}%
}%
\subfigure[RK-HAUSWENO-HLLC]{
\begin{minipage}[t]{0.45\linewidth}
\centering
\includegraphics[width=3.0in]{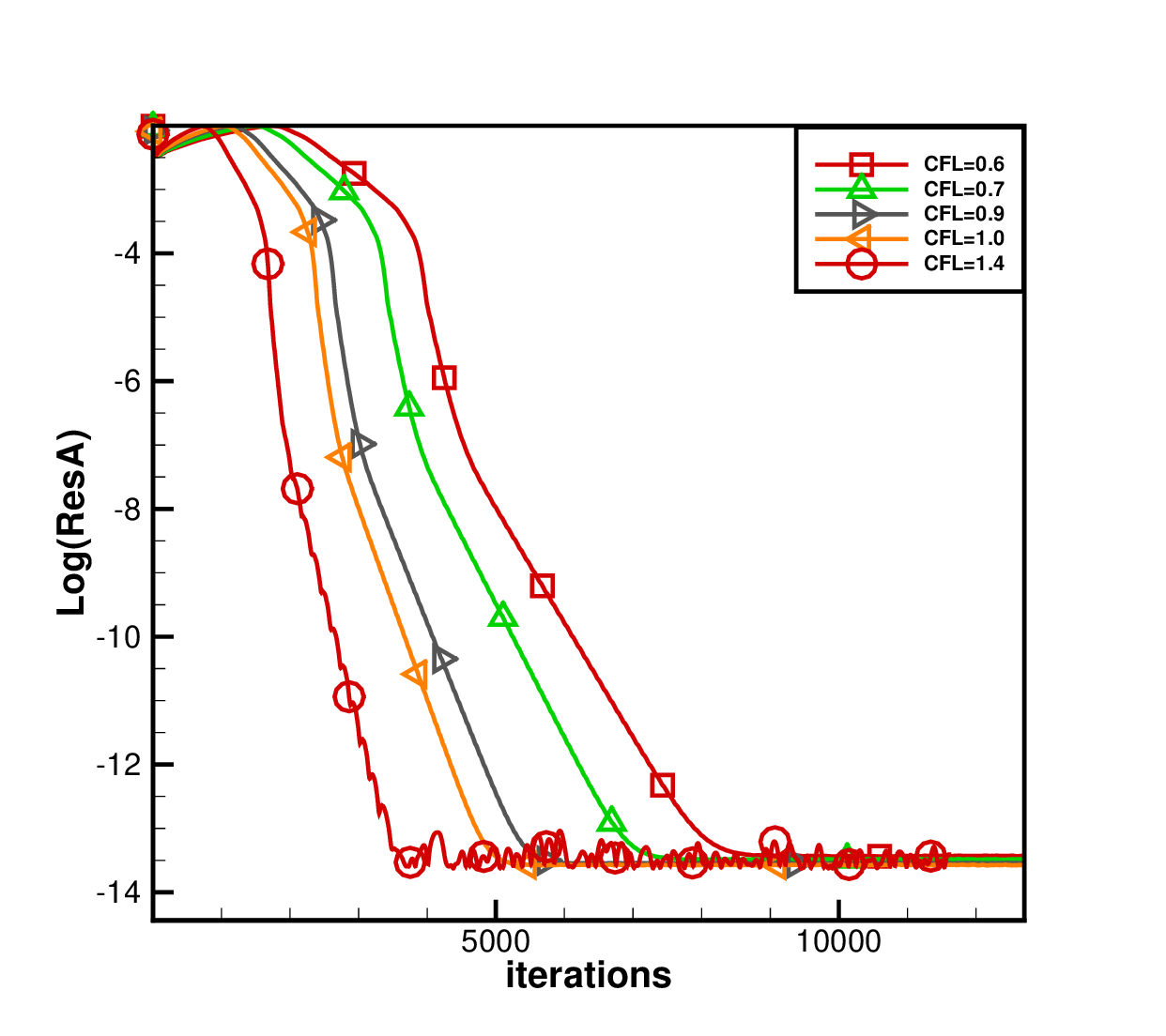}
\end{minipage}%
}%

\subfigure[FS-HAUSWENO-LLF]{
\begin{minipage}[t]{0.45\linewidth}
\centering
\includegraphics[width=3.0in]{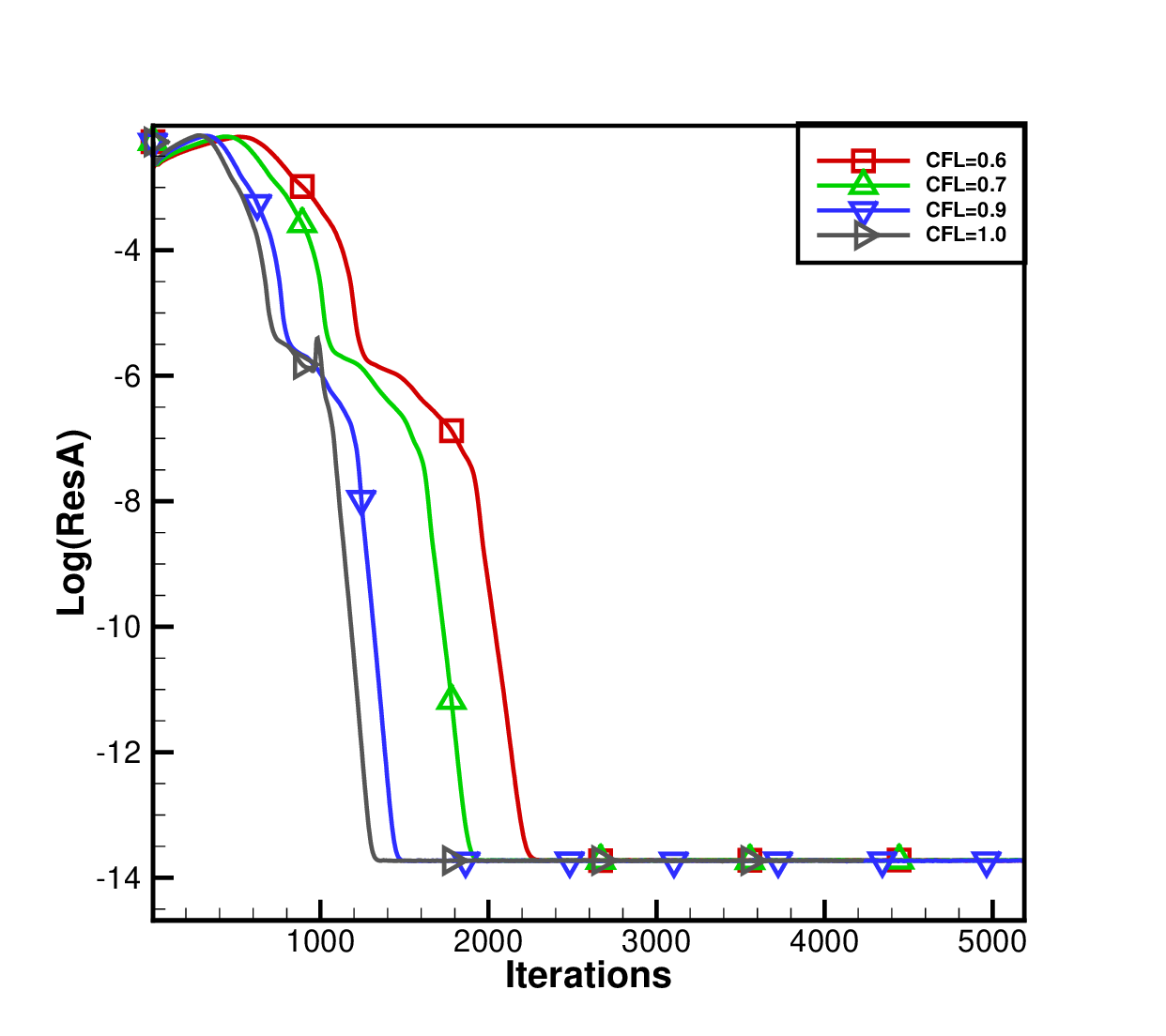}
\end{minipage}%
}
\subfigure[FS-HAUSWENO-HLLC]{
\begin{minipage}[t]{0.45\linewidth}
\centering
\includegraphics[width=3.0in]{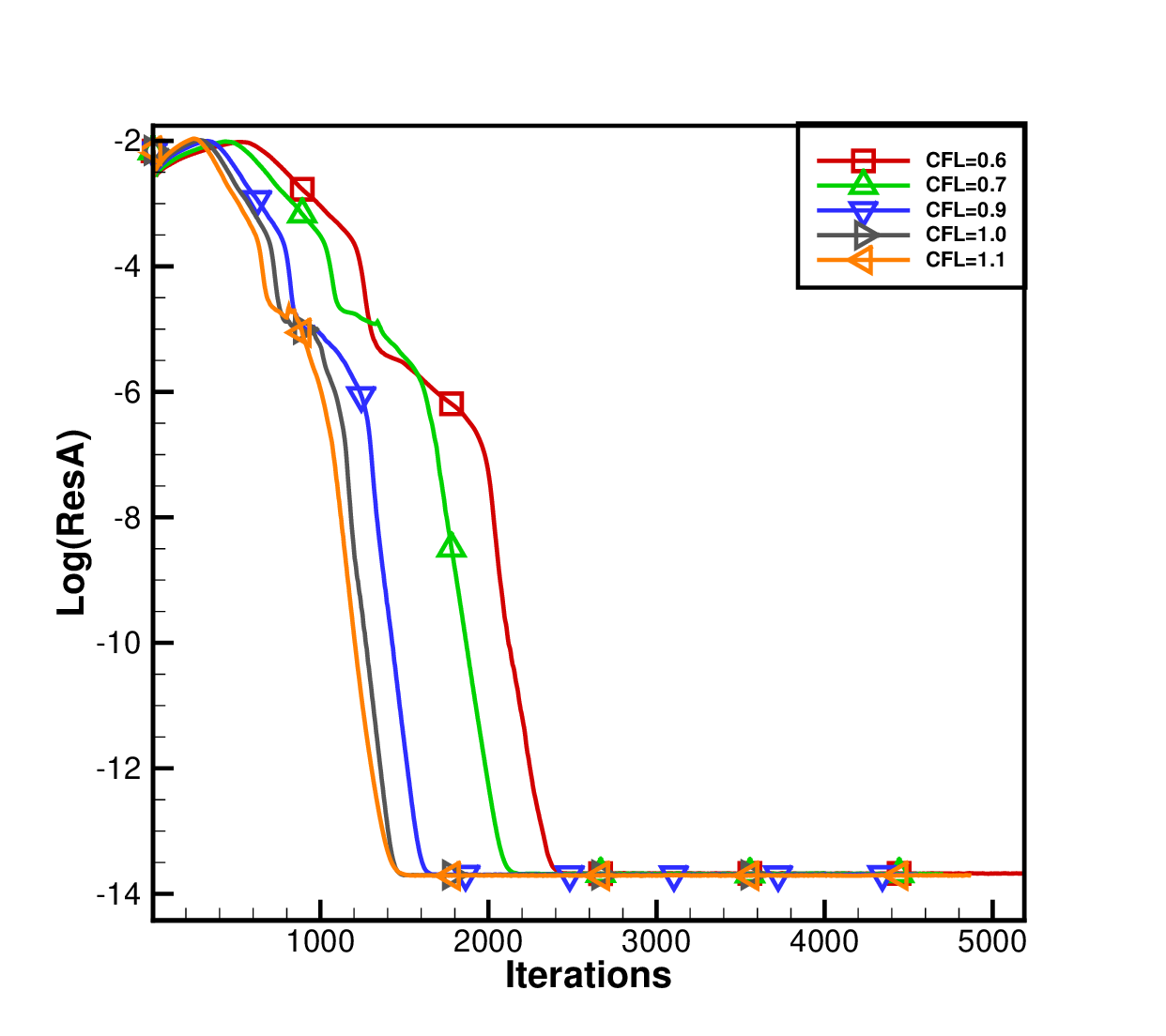}
\end{minipage}%
}%
\centering
\caption{\label{3.3}Example 3: Supersonic flow past a plate with an attack angle. The convergence history of the residue as a function of number of iterations for different schemes with various CFL numbers.}
\end{figure}

In this example, we consider a supersonic flow that traverses a plate characterized by an attack angle of $\alpha=10 ^{\circ}$, which is also a challenging steady-state problem for the residue of high order numerical schemes to converge to small values at the round off error level \cite{LZZ}. The free stream Mach number is $M_{\infty}=3$. The ideal gas flows in a direction from the left towards the plate, which is located at the interval $x \in [1,2]$ and $y=0$. The initial conditions are $p=\frac{1}{\gamma' M_{\infty}^{2}}$, $\rho$=1, $u=\cos(\alpha)$ and $v=\sin(\alpha)$ where $\gamma'=1.4$. The computational domain is $[0,10]\times [-5,5]$, and the computational grid is $200\times 200$. The slip boundary condition is imposed on the plate. The physical values of the inflow and the outflow boundary conditions are applied in different directions. The threshold value of the convergence criterion is set to $10^{-12}$. Fig.~\ref{3.2} shows the plots of the pressure contours of the numerical steady-state solutions obtained using these two hybrid schemes (FS-HAUSWENO and RK-HAUSWENO), and the cells where the WENO interpolation is used in the FS-HAUSWENO scheme. The comparable results are observed for both hybrid schemes. Again, it is observed that the hybrid technique identifies the region of troubled-cells pretty well in general, except some cells around the line $y=0$. Fig.~\ref{3.3} shows residue history in terms of iterations for the FS-HAUSWENO and the RK-HAUSWENO schemes, with the LLF flux and the HLLC flux under various CFL numbers. It is observed that the residue of iterations settles down to values at the level of round off errors for all cases, and the absolute convergence of the proposed hybrid fast sweeping method is achieved.
Table \ref{3.1} reports the number of iterations, the final time, and the total CPU time for the RK-HAUSWENO scheme, the FS-AUSWENO scheme and the FS-HAUSWENO scheme with various CFL numbers to achieve absolute convergence. It is verified that the fast sweeping method takes fewer iterations and is more efficient to converge to steady state than the TVD-RK method, and saves approximately up to $50\%$ CPU time cost. Comparison of the costs for the FS-AUSWENO scheme and the FS-HAUSWENO scheme shows a saving of up to $30\%$ CPU time using the hybrid technique in the fast sweeping AWENO method for this example.

\begin{table}
		\centering
\begin{tabular}{|c|c|c|c|}\hline
			\multicolumn{4}{|c|}{RK-HAUSWENO-LLF }\\\hline
            $\gamma:$ CFL number & iteration number & final time & CPU time \\\hline
0.6	&5466	&20.06	&231.11\\\hline
0.7	&4683	&20.05	&206.55\\\hline
0.8	&not conv	&	&\\
		\end{tabular}
\begin{tabular}{|c|c|c|c|}\hline
			\multicolumn{4}{|c|}{RK-HAUSWENO-HLLC }\\\hline
            $\gamma:$ CFL number & iteration number & final time & CPU time \\\hline
0.6	&7242	&26.62	&294.03\\\hline
0.7	&6204	&26.60	&245.52\\\hline
0.9	&4824	&26.60	&208.33\\\hline
1.0	&4341	&26.59	&170.58\\\hline
1.4	&3138	&26.91	&122.34\\\hline
1.5	&not conv	&	&\\
		\end{tabular}
\begin{tabular}{|c|c|c|c|}\hline
			\multicolumn{4}{|c|}{FS-AUSWENO-LLF }\\\hline
            $\gamma:$ CFL number & iteration number & final time & CPU time \\\hline
0.6	&2140	&23.85	&244.83\\\hline
0.7	&1816	&23.62	&208.31\\\hline
0.9	&1388	&23.22	&159.66\\\hline
1.0	&1248	&23.20	&140.92\\\hline
1.1	&not conv	&	&\\
		\end{tabular}
\begin{tabular}{|c|c|c|c|}\hline
			\multicolumn{4}{|c|}{FS-AUSWENO-HLLC}\\\hline
            $\gamma:$ CFL number & iteration number & final time & CPU time \\\hline
0.6	&2248	&25.02	&228.80\\\hline
0.7	&1920	&24.93	&193.09\\\hline
0.9	&1516	&25.31	&152.73\\\hline
1.0	&1360	&25.23	&137.09\\\hline
1.1	&1304	&26.62	&130.45\\\hline
1.2	&not conv	&	&\\
		\end{tabular}
\begin{tabular}{|c|c|c|c|}\hline
			\multicolumn{4}{|c|}{FS-HAUSWENO-LLF }\\\hline
            $\gamma:$ CFL number & iteration number & final time & CPU time \\\hline
0.6	&2140	&23.85	&185.81\\\hline
0.7	&1816	&23.62	&155.61\\\hline
0.9	&1388	&23.21	&120.16\\\hline
1.0	&1248	&23.20	&106.31\\\hline
1.1	&not conv	&	&\\
		\end{tabular}
\begin{tabular}{|c|c|c|c|}\hline
			\multicolumn{4}{|c|}{FS-HAUSWENO-HLLC}\\\hline
            $\gamma:$ CFL number & iteration number & final time & CPU time \\\hline
0.6	&2248	&25.02	&168.41\\\hline
0.7	&1920	&24.93	&140.06\\\hline
0.9	&1516	&25.31	&111.31\\\hline
1.0	&1360	&25.23	&100.86\\\hline
1.1	&1308	&26.70	&91.22\\\hline
1.2	&not conv	&	&\\\hline
		\end{tabular}
		\caption{\label{3.1}Example 3: Supersonic flow past a plate with an attack angle. Number of iterations, the final time and total CPU time when convergence is obtained. Convergence criterion threshold value is $10^{-12}$. CPU time unit: second.}
	\end{table}

\bigskip
\noindent{\bf Example 4: Flow around a circular cylinder}

\begin{figure}
\centering
\subfigure[RK-HAUSWENO-LLF]{
\begin{minipage}[t]{0.3\linewidth}
\centering
\includegraphics[width=2.0in]{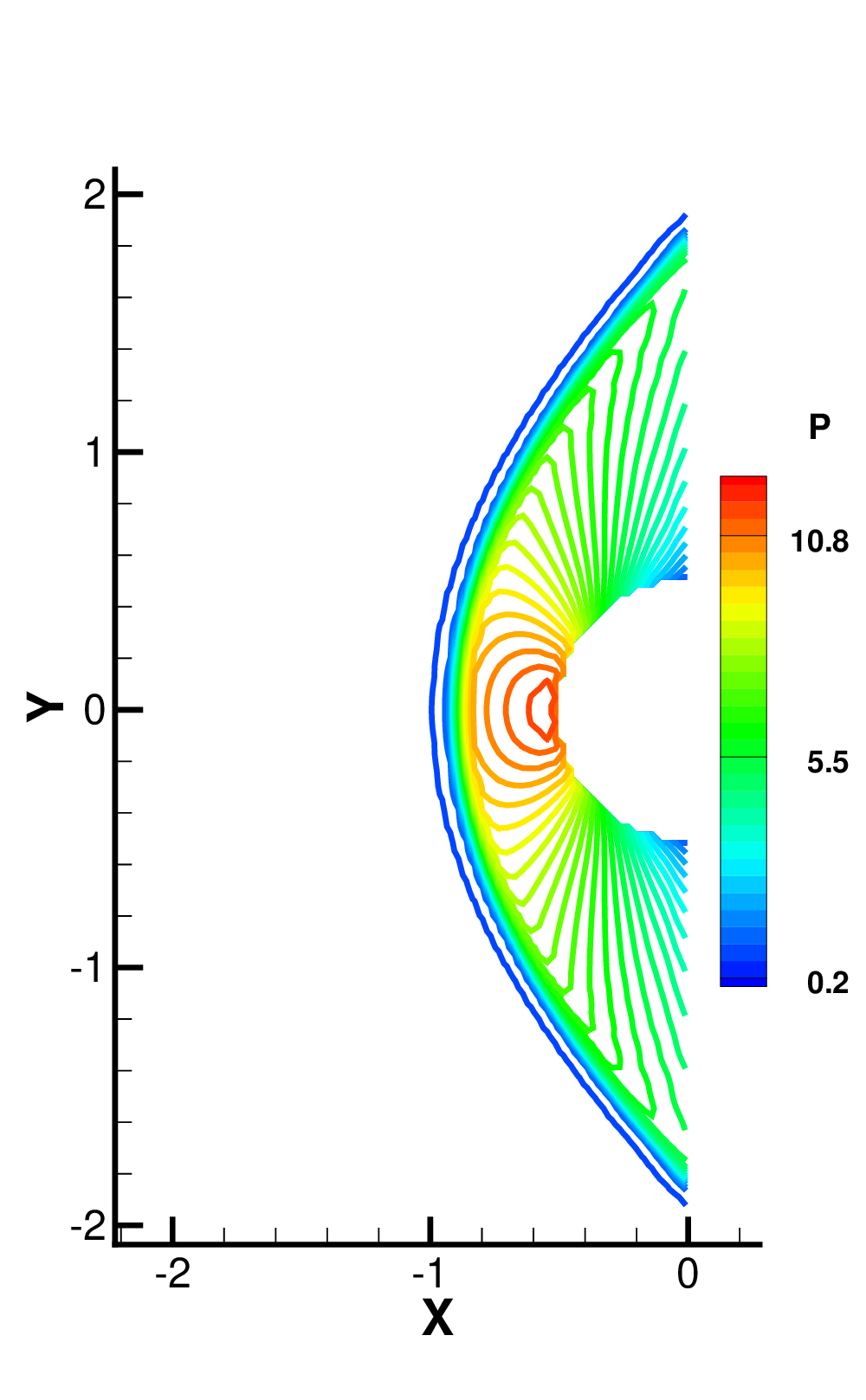}
\end{minipage}%
}
\subfigure[FS-HAUSWENO-LLF]{
\begin{minipage}[t]{0.3\linewidth}
\centering
\includegraphics[width=2.0in]{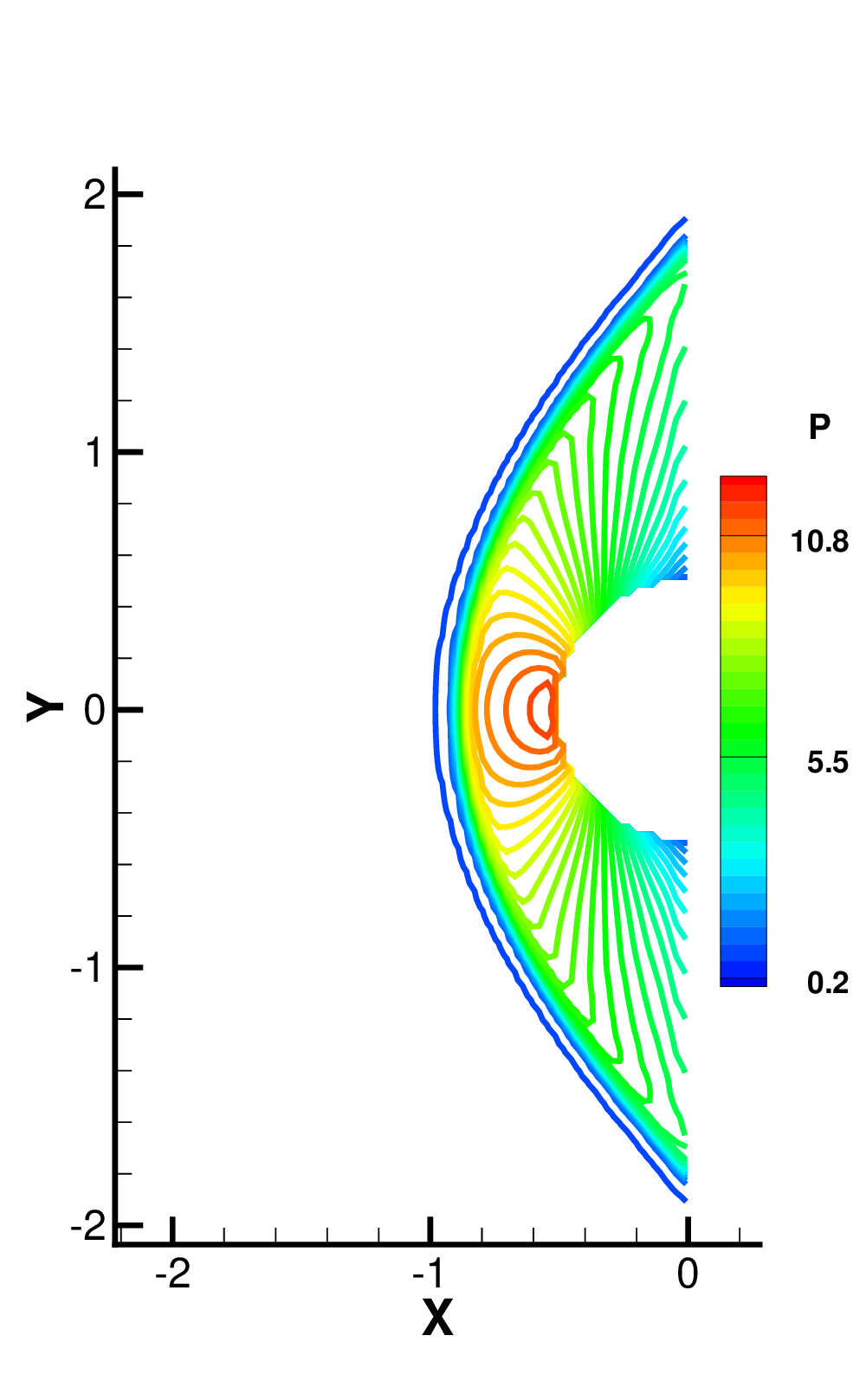}
\end{minipage}%
}%
\subfigure[FS-HAUSWENO-LLF]{
\begin{minipage}[t]{0.3\linewidth}
\centering
\includegraphics[width=2.0in]{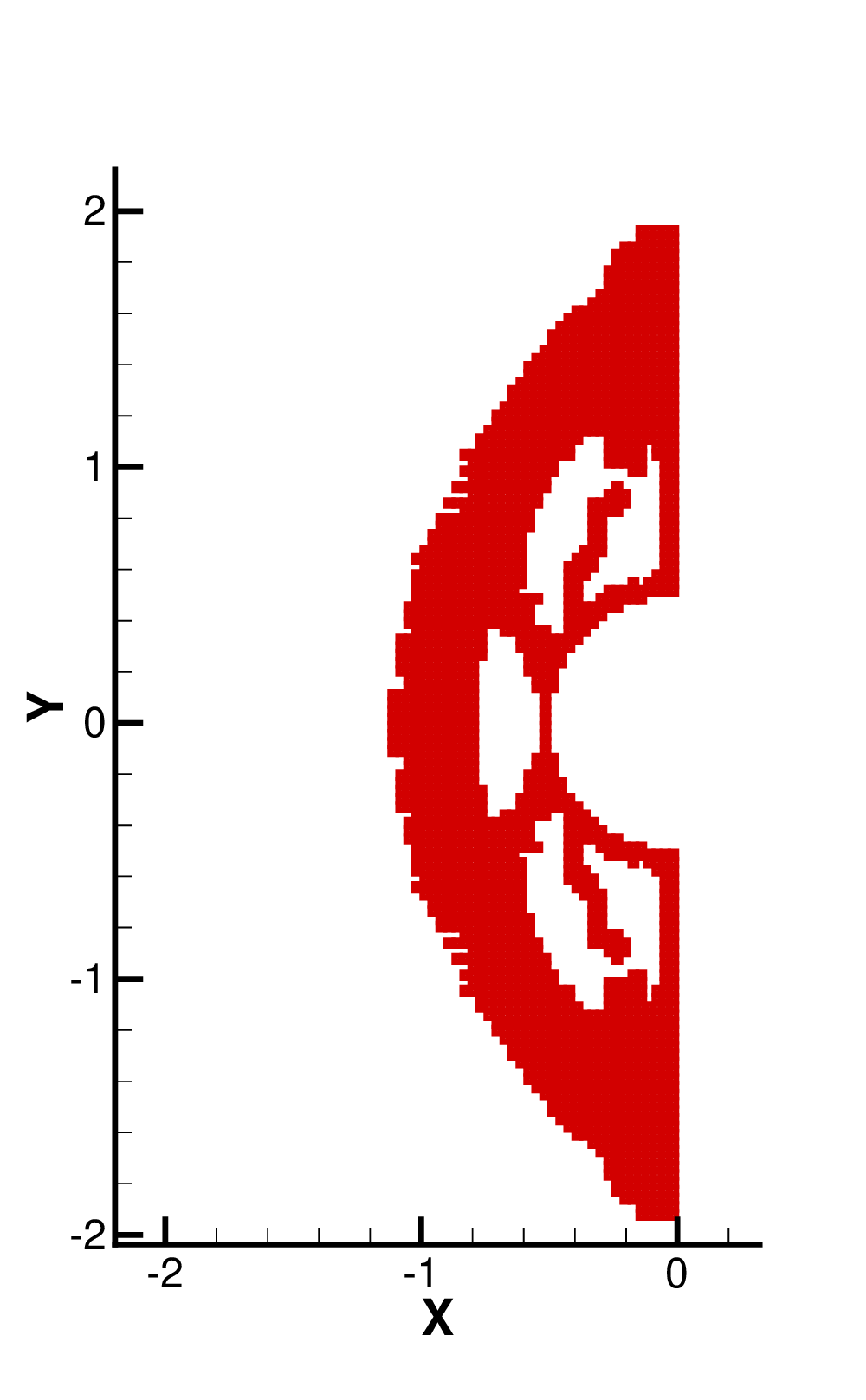}
\end{minipage}%
}
\caption{\label{4.2}Example 4: Flow around a circular cylinder. 30 equally spaced pressure contour from 0.2 to 12 of the converged steady states of numerical solutions by different iterative schemes, and the identified troubled-cells where the WENO interpolation is used in the FS-HAUSWENO scheme.}
\end{figure}

\begin{figure}
\centering
\subfigure[RK-HAUSWENO-LLF]{
\begin{minipage}[t]{0.4\linewidth}
\centering
\includegraphics[width=3.0in]{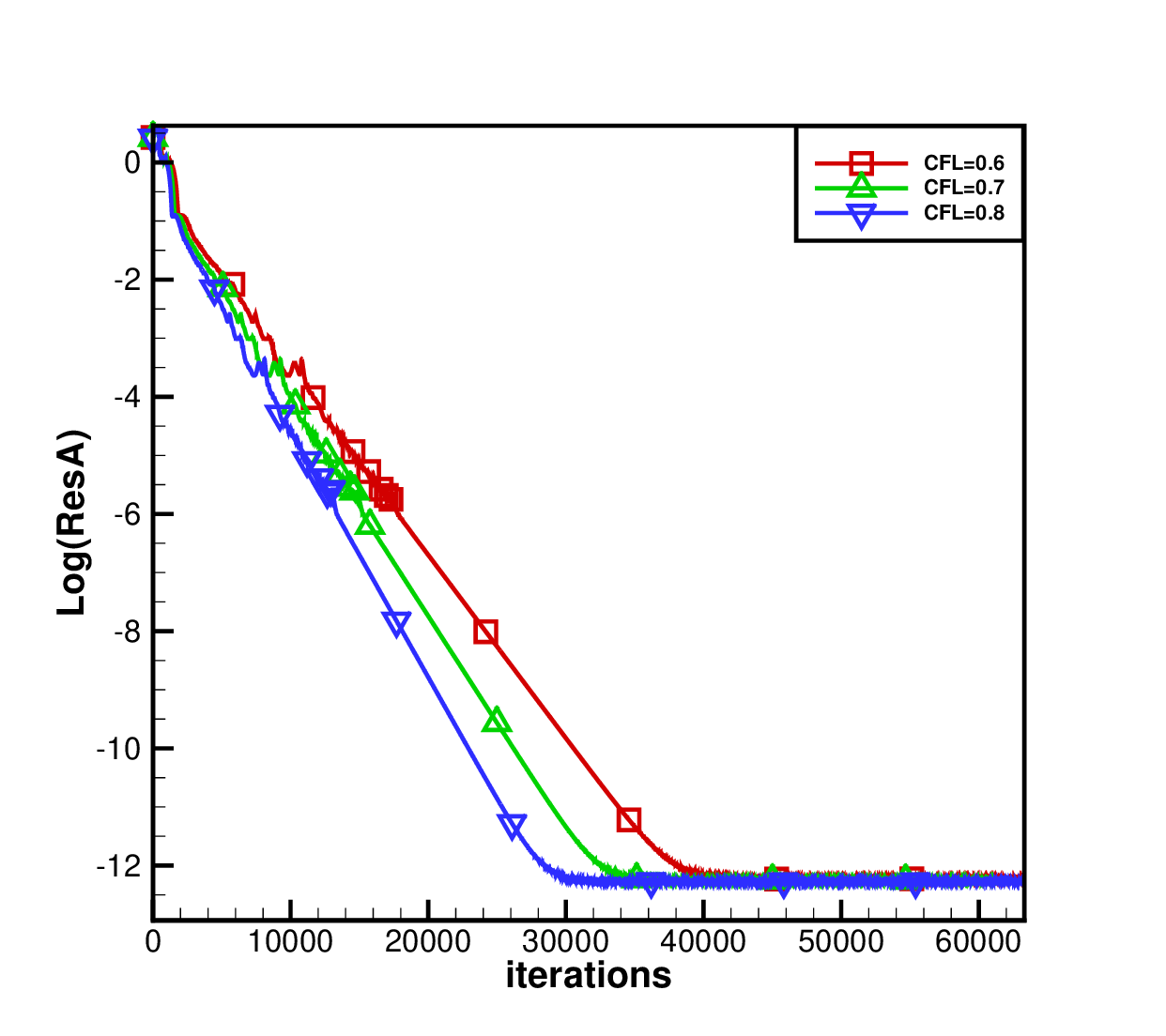}
\end{minipage}%
}
\subfigure[FS-HAUSWENO-LLF]{
\begin{minipage}[t]{0.4\linewidth}
\centering
\includegraphics[width=3.0in]{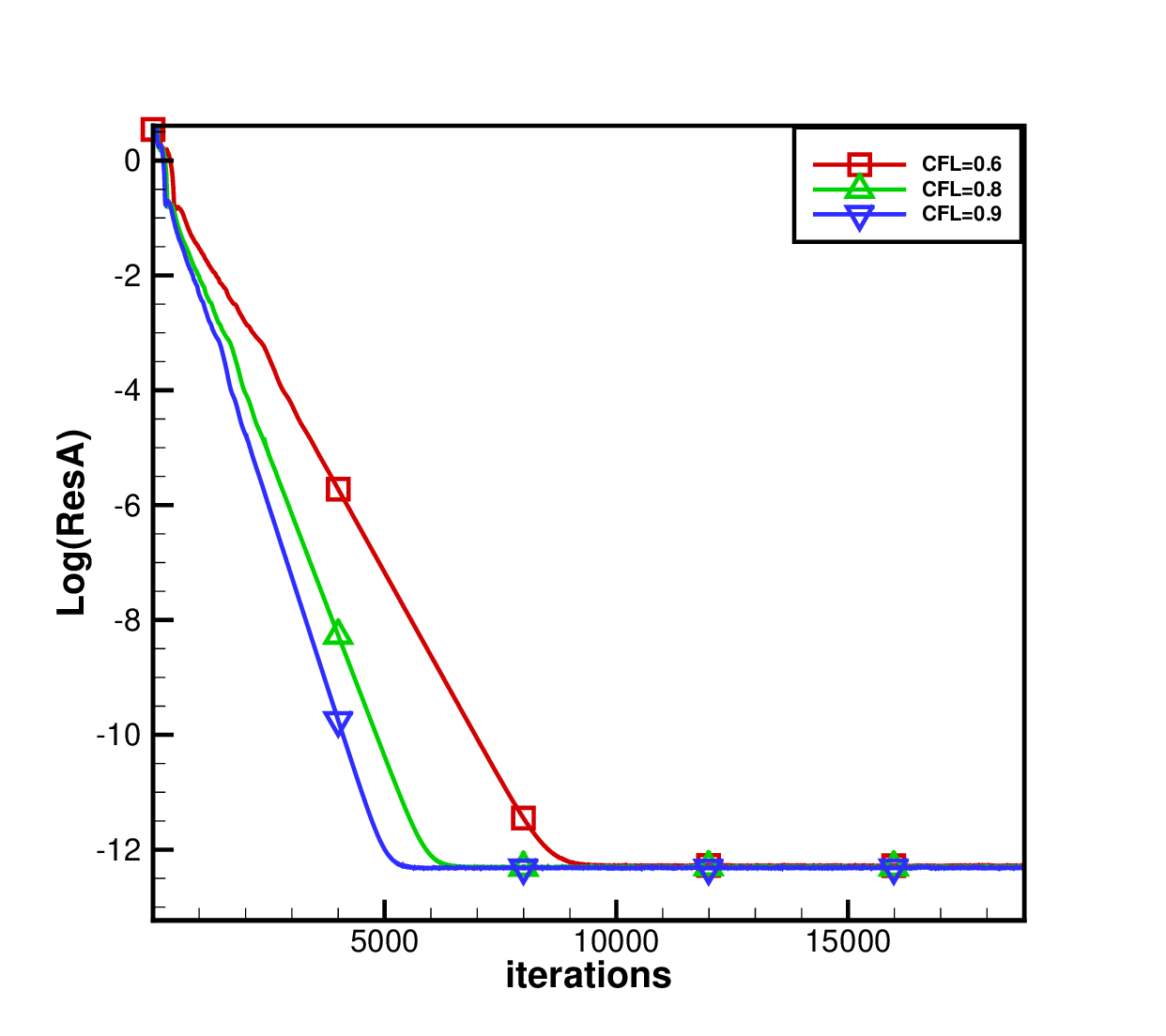}
\end{minipage}%
}%
\caption{\label{4.3}Example 4: Flow around a circular cylinder. The convergence history of the residue as a function of number of iterations for different  schemes with various CFL numbers.}
\end{figure}

In this example, we consider a supersonic flow past a curved wall which is a circular cylinder of unit diameter positioned at the origin on the $x-y$ plane. The setup of this problem is the same as in $\cite{SIRUIT}$. The computation domain has the range $[-5,0]\times[-5,5]$. The simulation is initialized with a Mach 3 flow advancing towards the cylinder from the left-hand side. Along the left inflow boundary at $x=-5$, we impose the uniform far-field data. For the top boundary at $y=5$ and the right boundary at $x=0$, constant extrapolation is applied as in \cite{S.Tan}. Due to the symmetry of this problem, we perform the simulation at the upper half of the whole physical domain, and at $y=0$, the reflection technique is implemented. The fifth-order inverse Lax-Wendroff procedure is used to handle the surface of the solid wall cylinder at the boundary as in \cite{LILW}. Via thorough numerical tests, we observe that both of the hybrid iterative schemes (FS-HAUSWENO and RK-HAUSWENO) exhibit absolute convergence when the LLF flux is used. However, absolute convergence can not be achieved when the HLLC flux is applied to them. This example shows that although the HLLC flux has lower dissipation, the LLF flux is more robust for the proposed hybrid fast sweeping AWENO method to achieve the convergence. The contour plots of the pressure variable of the converged steady-state numerical solutions obtained using these two hybrid iterative schemes (FS-HAUSWENO and RK-HAUSWENO) with the grid size $\Delta x=\Delta y=\frac{1}{32}$, and the cells where the WENO interpolation is used in the FS-HAUSWENO scheme are shown in Fig.~$\ref{4.2}$. We observe similar results for both iterative schemes and the well-captured bow shock. Similar to the previous examples, it is observed that the hybrid technique identifies the region of troubled-cells pretty well.  Residue history in terms of iterations for these two hybrid schemes under various CFL numbers is presented in Fig.~$\ref{4.3}$. The iteration residuals settle down to values less than $10^{-12}$, a level of round off errors. As in the previous examples, the absolute convergence of the proposed hybrid fast sweeping AWENO method is achieved. The number of iterations, the final time, and the total CPU time for the RK-HAUSWENO scheme, the FS-AUSWENO scheme and the FS-HAUSWENO scheme with various CFL numbers to achieve absolute convergence are reported in Table $\ref{4.1}$. In this example, significant improvement is observed on the computational efficiency of the fast sweeping method over the TVD-RK method. Approximately $65\%$ CPU time cost is saved using the FS-HAUSWENO scheme, compared to the RK-HAUSWENO scheme. Table \ref{4.1} also shows that the hybrid fast sweeping method (FS-HAUSWENO) is more efficient than the scheme without the hybrid technique (FS-AUSWENO).  These comparisons again verify that both the fast sweeping and the hybrid techniques help in enhancing the efficiency for the iterative schemes to converge.
\begin{table}
		\centering
\begin{tabular}{|c|c|c|c|}\hline
			\multicolumn{4}{|c|}{RK-HAUSWENO-LLF}\\\hline
            $\gamma:$ CFL number & iteration number & final time & CPU time \\\hline
0.6	&38127	&35.85	&1091.95 \\\hline
0.8	&28446	&35.67	&895.52 \\\hline
1.0	&not conv	&	&\\
		\end{tabular}
\begin{tabular}{|c|c|c|c|}\hline
			\multicolumn{4}{|c|}{FS-AUSWENO-LLF}\\\hline
            $\gamma:$ CFL number & iteration number & final time & CPU time \\\hline
0.6	&8544	&24.22	&647.72\\\hline
0.9	&5008	&21.28	&374.14\\\hline
1.0	&not conv	&	&\\
		\end{tabular}
\begin{tabular}{|c|c|c|c|}\hline
			\multicolumn{4}{|c|}{FS-HAUSWENO-LLF}\\\hline
            $\gamma:$ CFL number & iteration number & final time & CPU time \\\hline
0.6	&8540	&24.21	&534.05\\\hline
0.9	&5008	&21.30	&315.02\\\hline
1.0	&not conv	&	&\\\hline
		\end{tabular}
		\caption{\label{4.1}Example 4: Flow around a circular cylinder. Number of iterations, the final time and total CPU time when convergence is obtained. Convergence criterion threshold value  is $10^{-12}$. CPU time unit: second.}
	\end{table}

\bigskip
\noindent{\bf Example 5: A Mach 3 wind tunnel with a forward-facing step}

\begin{figure}
\centering
\subfigure[RK-HAUSWENO-LLF]{
\begin{minipage}[t]{0.32\linewidth}
\centering
\includegraphics[width=2.1in]{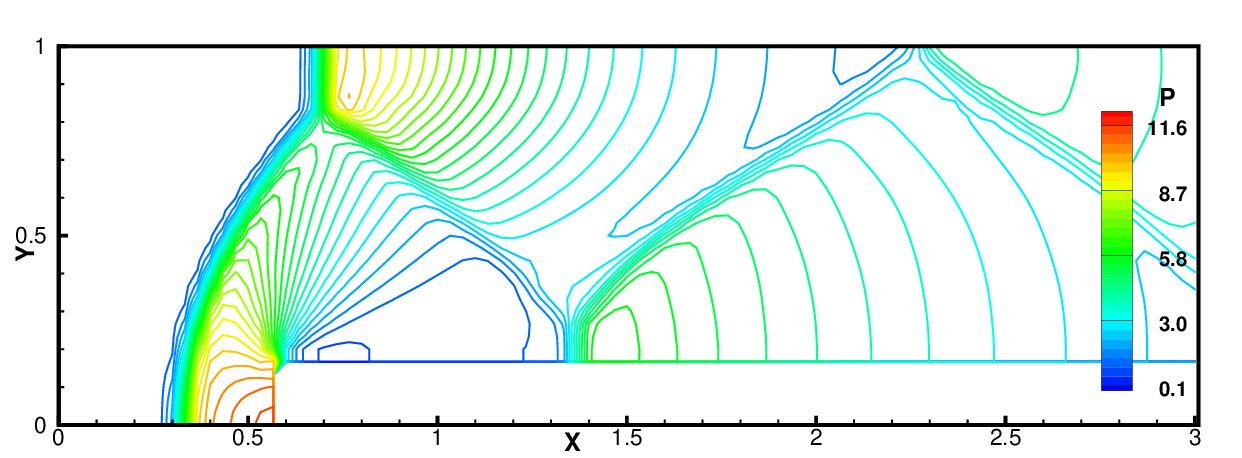}
\end{minipage}%
}%
\subfigure[FS-HAUSWENO-LLF]{
\begin{minipage}[t]{0.32\linewidth}
\centering
\includegraphics[width=2.1in]{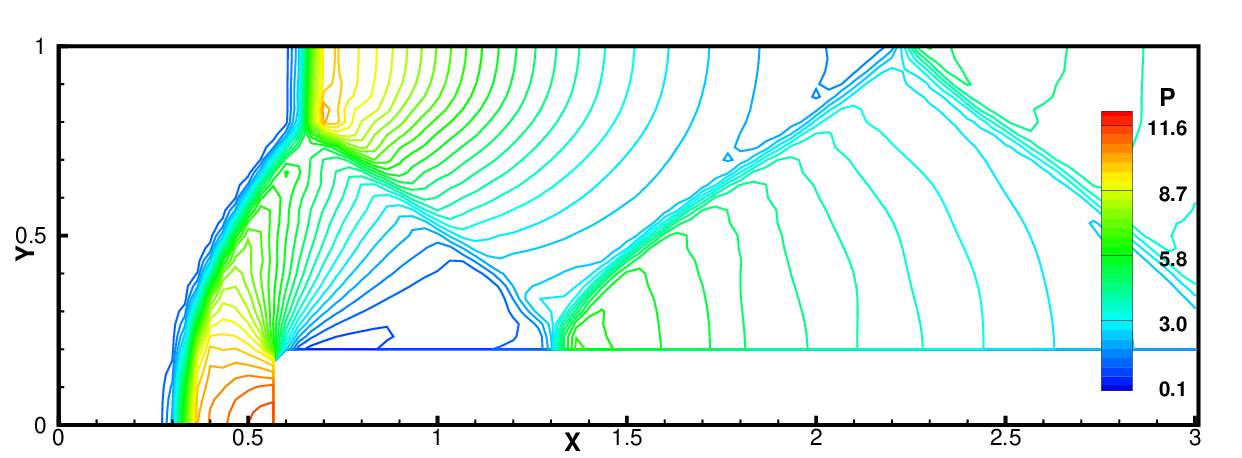}
\end{minipage}%
}%
\subfigure[FS-HAUSWENO-LLF]{
\begin{minipage}[t]{0.32\linewidth}
\centering
\includegraphics[width=2.0in]{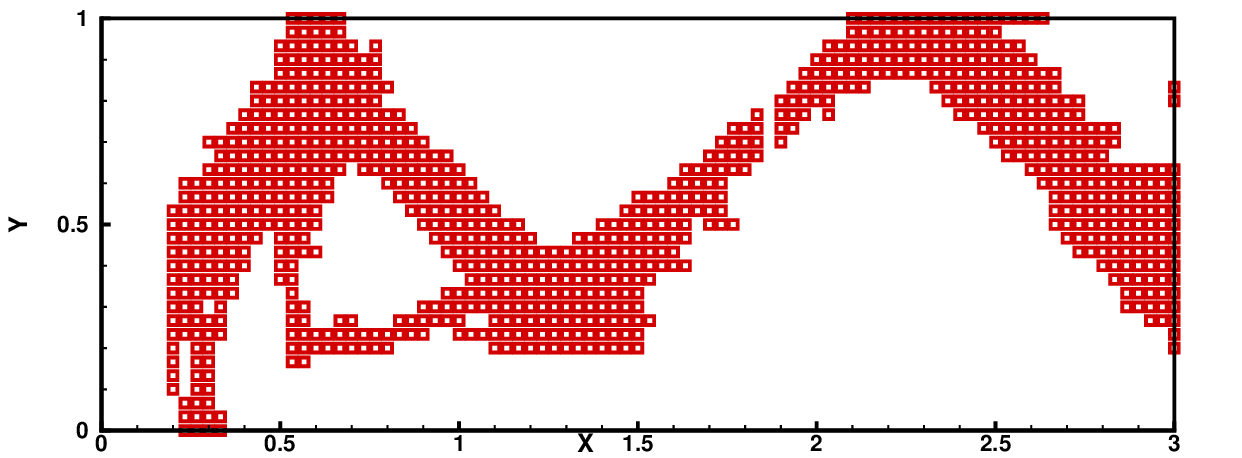}
\end{minipage}%
}
\centering
\caption{\label{5.2}Example 5: A Mach 3 wind tunnel with a forward-facing step. 30 equally spaced pressure contour from 0.1 to 12 of the converged steady states of numerical solutions by different iterative schemes and the identified troubled-cells where the WENO interpolation is used in the FS-HAUSWENO scheme.}
\end{figure}

\begin{figure}
\centering
\subfigure[RK-HAUSWENO-LLF]{
\begin{minipage}[t]{0.4\linewidth}
\centering
\includegraphics[width=3.0in]{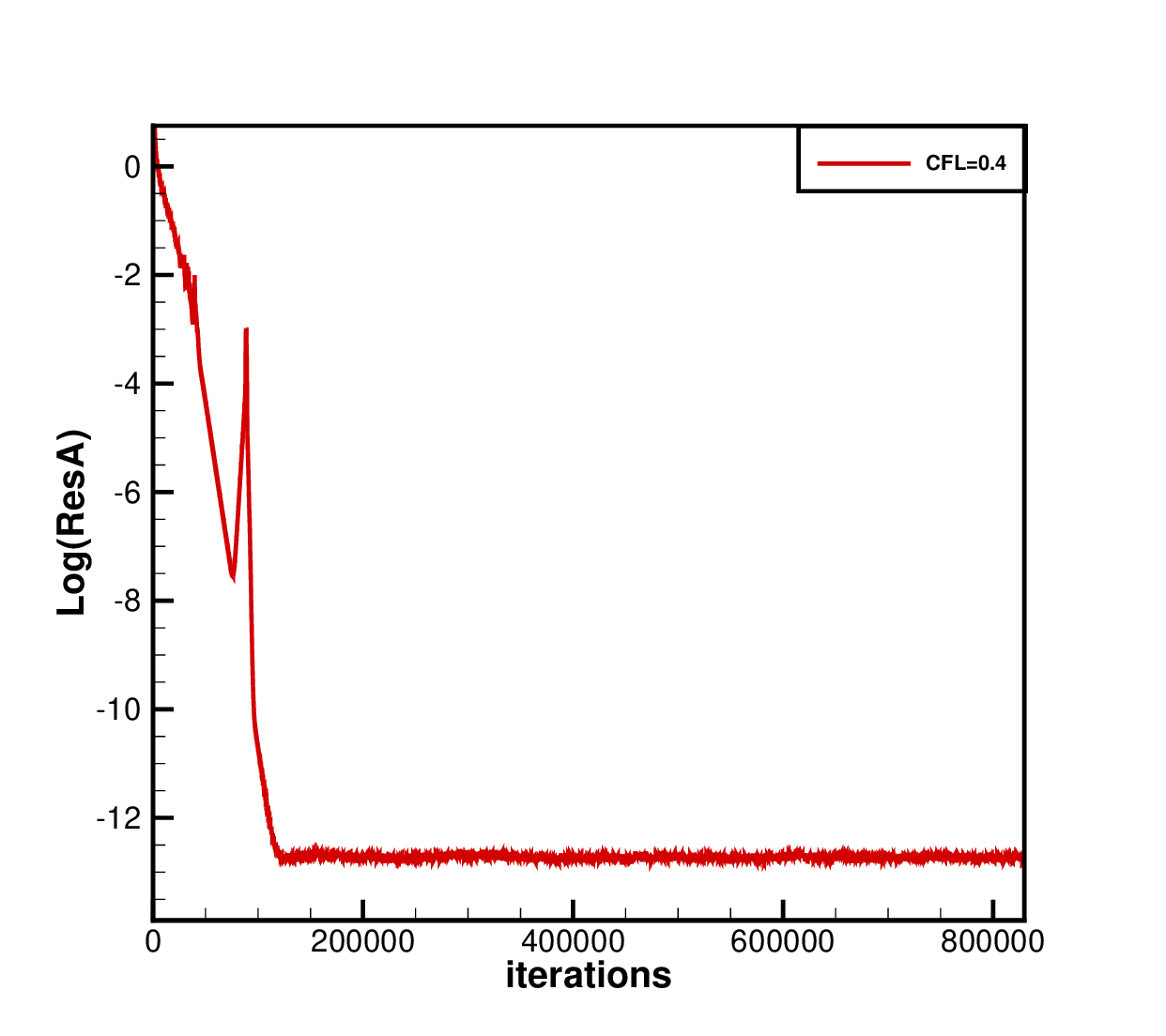}
\end{minipage}%
}%
\subfigure[FS-HAUSWENO-LLF]{
\begin{minipage}[t]{0.4\linewidth}
\centering
\includegraphics[width=3.0in]{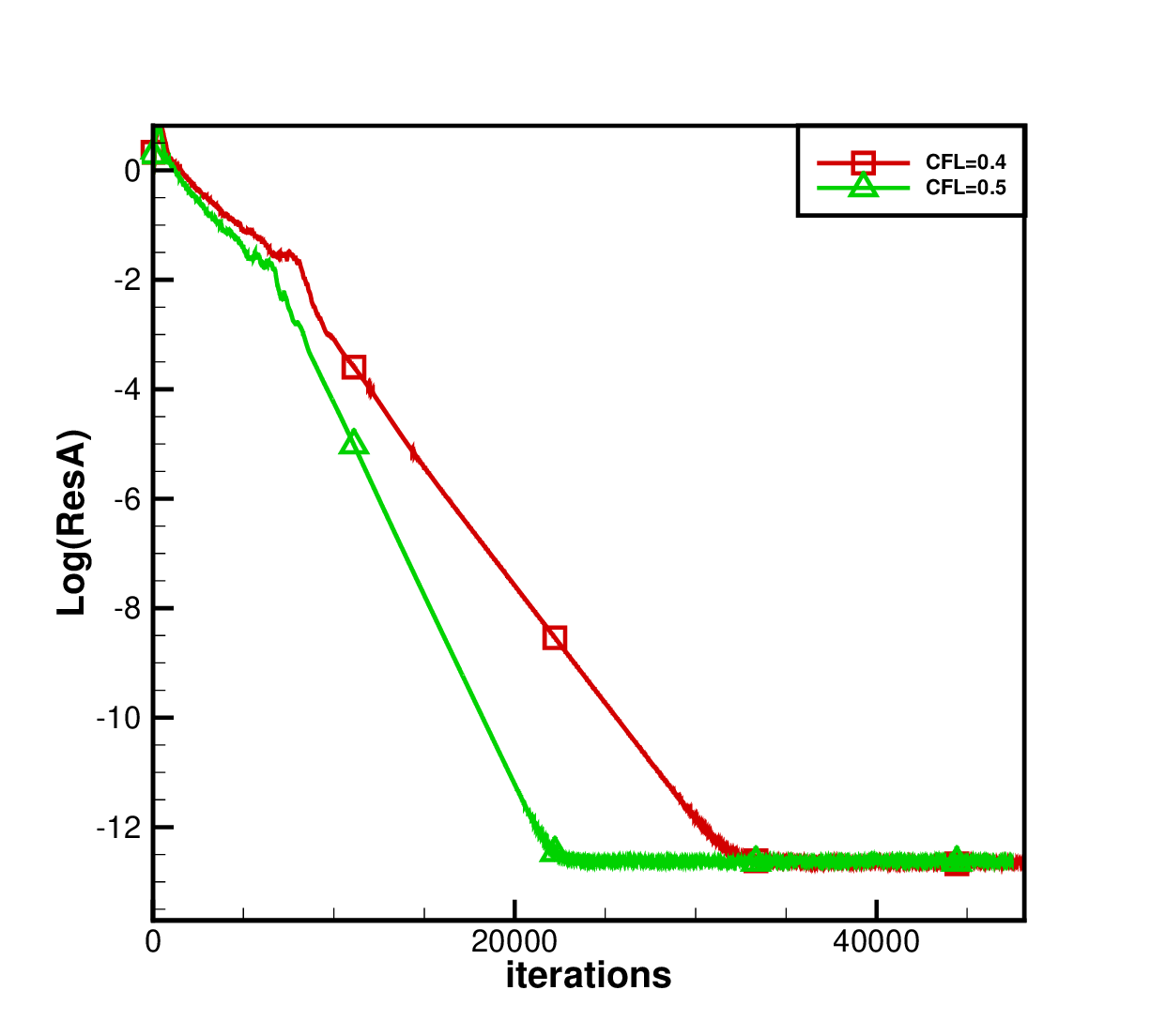}
\end{minipage}%
}%
\centering
\caption{\label{5.3}Example 5: A Mach 3 wind tunnel with a forward-facing step. The convergence history of the residue as a function of number of iterations for two hybrid schemes.}
\end{figure}

In this example, we apply the proposed hybrid fast sweeping AWENO method for solving the Mach 3 wind tunnel flow problem of Woodward and Colella \cite{taijie1}.  The setup of the problem is as following. The wind tunnel is 1 length unit wide and 3 length units long. The step is 0.2 length units high, and it is located at the bottom of the
wind tunnel from a right going Mach 3 flow. The computation grid is $90\times30$. Reflective boundary conditions are applied along the walls of the tunnel, and inflow and outflow boundary conditions are applied at the entrance and the exit. There is a singularity point at the corner of the step. We apply the same procedure as specified in \cite{taijie1, JZhu-CShu17} and use an assumption of nearly steady flow in the region near the corner to fix the singularity.
In Table \ref{5.1}, the number of iterations, the final time, and the total CPU time for the RK-HAUSWENO scheme, the FS-AUSWENO scheme and the FS-HAUSWENO scheme with various CFL numbers to achieve absolute convergence are reported. Similar to Example 4, it is observed that both of the hybrid iterative schemes (FS-HAUSWENO and RK-HAUSWENO) exhibit absolute convergence when the LLF flux is used, but absolute convergence can not be achieved when the HLLC flux is applied to them. This example also verifies that the LLF flux is more robust than the HLLC flux for the proposed hybrid fast sweeping AWENO method to achieve the convergence. Results in Table \ref{5.1} show significant improvement on the computational efficiency of the fast sweeping method over the TVD-RK method. Approximately $75\%$ CPU time cost is saved using the FS-HAUSWENO scheme with the largest CFL number permitted to converge to steady state, compared to the RK-HAUSWENO scheme. Also, the hybrid fast sweeping method (FS-HAUSWENO) is more efficient than the scheme without the hybrid technique (FS-AUSWENO).   The contour plots of the pressure variable of the converged steady-state numerical solutions obtained using these two hybrid iterative schemes (FS-HAUSWENO and RK-HAUSWENO), and the cells where the WENO interpolation is used in the FS-HAUSWENO scheme are shown in Fig.~\ref{5.2}. Similar to the previous examples, comparable results for both iterative schemes are obtained, and the hybrid technique identifies the region of troubled-cells well. The history of the residue as a function of iterations for these two hybrid schemes reported in Fig.~$\ref{5.3}$ shows that the iteration residuals settle down to values at a level of round off errors and the absolute convergence is achieved. 
\begin{table}
		\centering
\begin{tabular}{|c|c|c|c|}\hline
			\multicolumn{4}{|c|}{RK-HAUSWENO-LLF}\\\hline
            $\gamma:$ CFL number & iteration number & final time & CPU time \\\hline
0.4	&110421	&53.07	&479.13\\\hline
0.5	&not conv	&	&\\
		\end{tabular}
\begin{tabular}{|c|c|c|c|}\hline
			\multicolumn{4}{|c|}{FS-AUSWENO-LLF}\\\hline
            $\gamma:$ CFL number & iteration number & final time & CPU time \\\hline
0.4	&30396	&43.65	&208.28\\\hline
0.5	&21524	&38.76	&142.42\\\hline
0.6	&not conv	&	&\\
	\end{tabular}
\begin{tabular}{|c|c|c|c|}\hline
			\multicolumn{4}{|c|}{FS-HAUSWENO-LLF}\\\hline
            $\gamma:$ CFL number & iteration number & final time & CPU time \\\hline
0.4	&30396	&43.91	&177.30\\\hline
0.5	&21210	&38.19	&122.08\\\hline
0.6	&not conv	&	&\\\hline
		\end{tabular}
		\caption{\label{5.1}Example 5: A Mach 3 wind tunnel with a forward-facing step. Number of iterations, the final time and total CPU time when convergence is obtained. Convergence criterion threshold value is $10^{-12}$. CPU time unit: second.}
	\end{table}

\bigskip
\noindent{\bf Example 6: Supersonic flow past an NACA0012 airfoil}

\begin{figure}
\centering
\subfigure[RK-HAUSWENO-LLF]{
\begin{minipage}[t]{0.4\linewidth}
\centering
\includegraphics[width=2.6in]{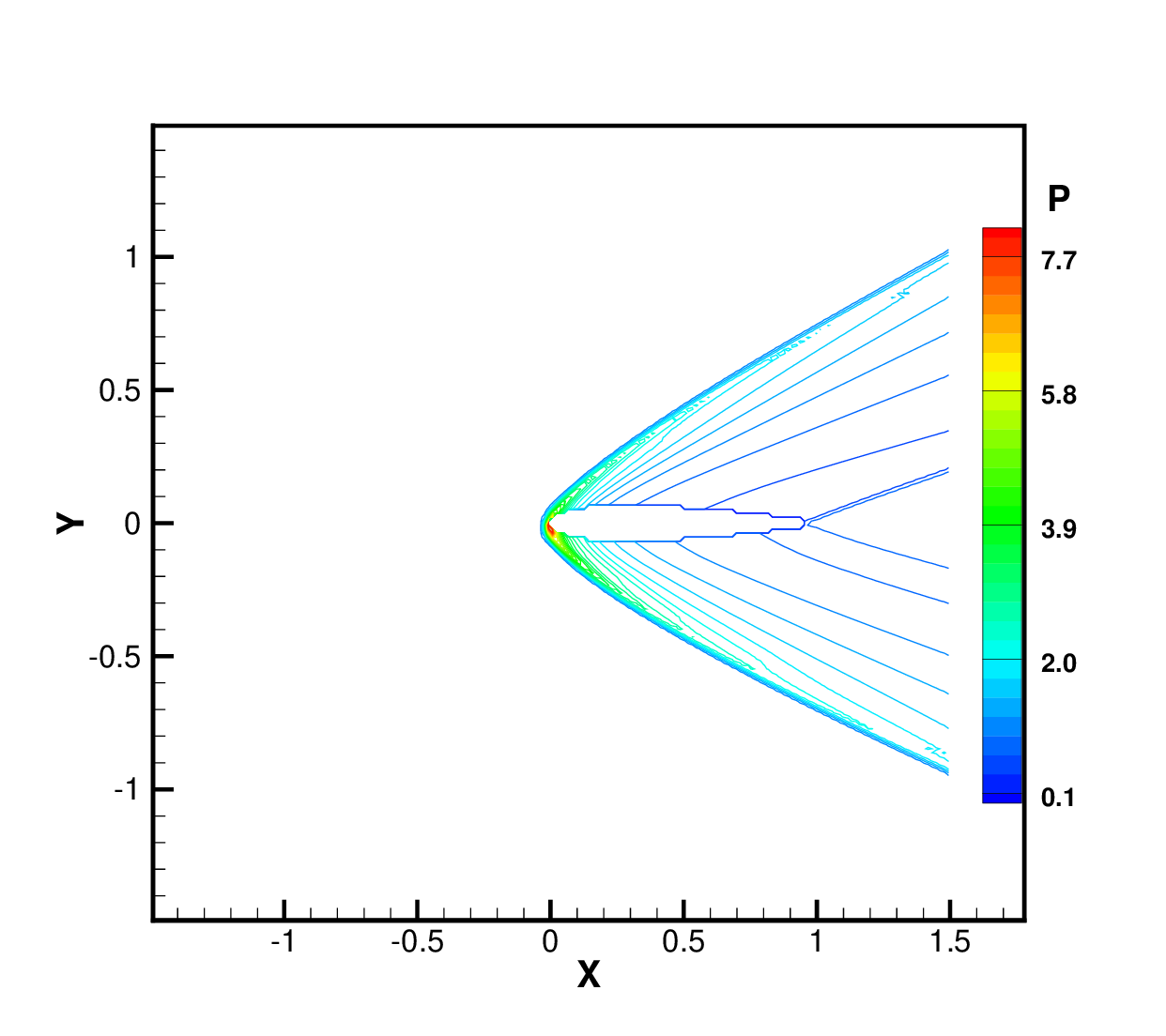}
\end{minipage}%
}
\subfigure[RK-HAUSWENO-HLLC]{
\begin{minipage}[t]{0.4\linewidth}
\centering
\includegraphics[width=2.6in]{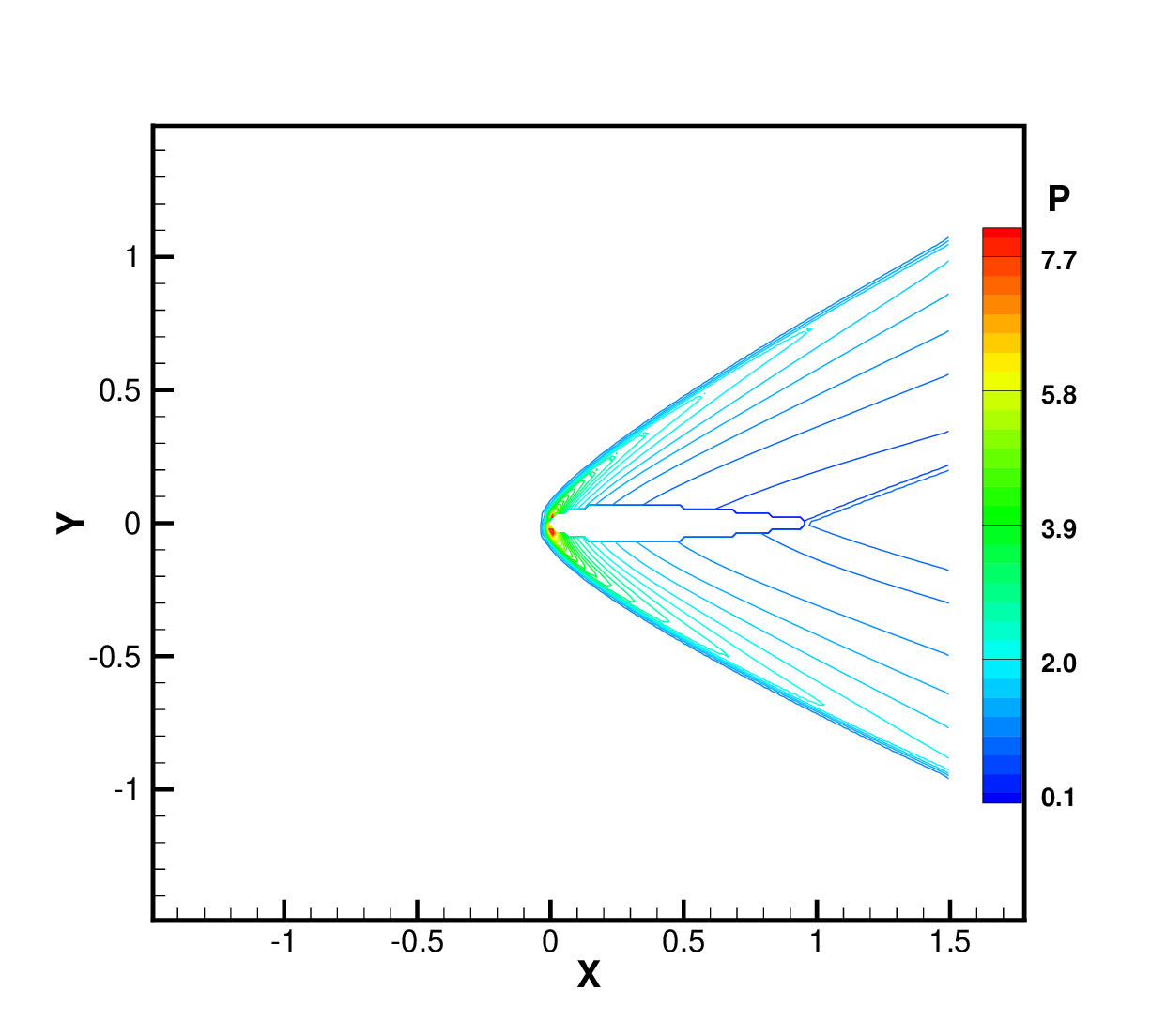}
\end{minipage}%
}%

\subfigure[FS-HAUSWENO-LLF]{
\begin{minipage}[t]{0.4\linewidth}
\centering
\includegraphics[width=2.6in]{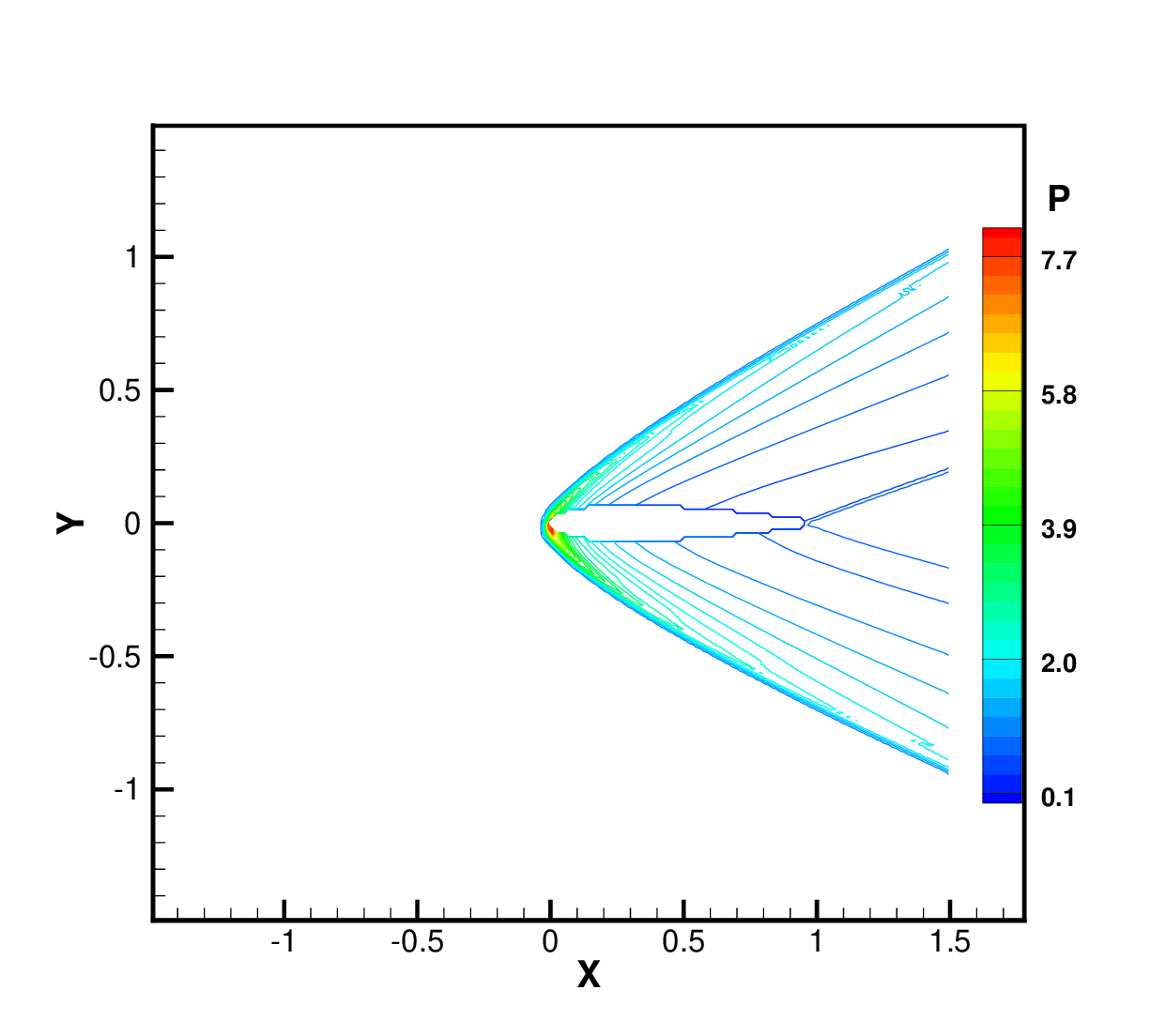}
\end{minipage}%
}
\subfigure[FS-HAUSWENO-HLLC]{
\begin{minipage}[t]{0.4\linewidth}
\centering
\includegraphics[width=2.6in]{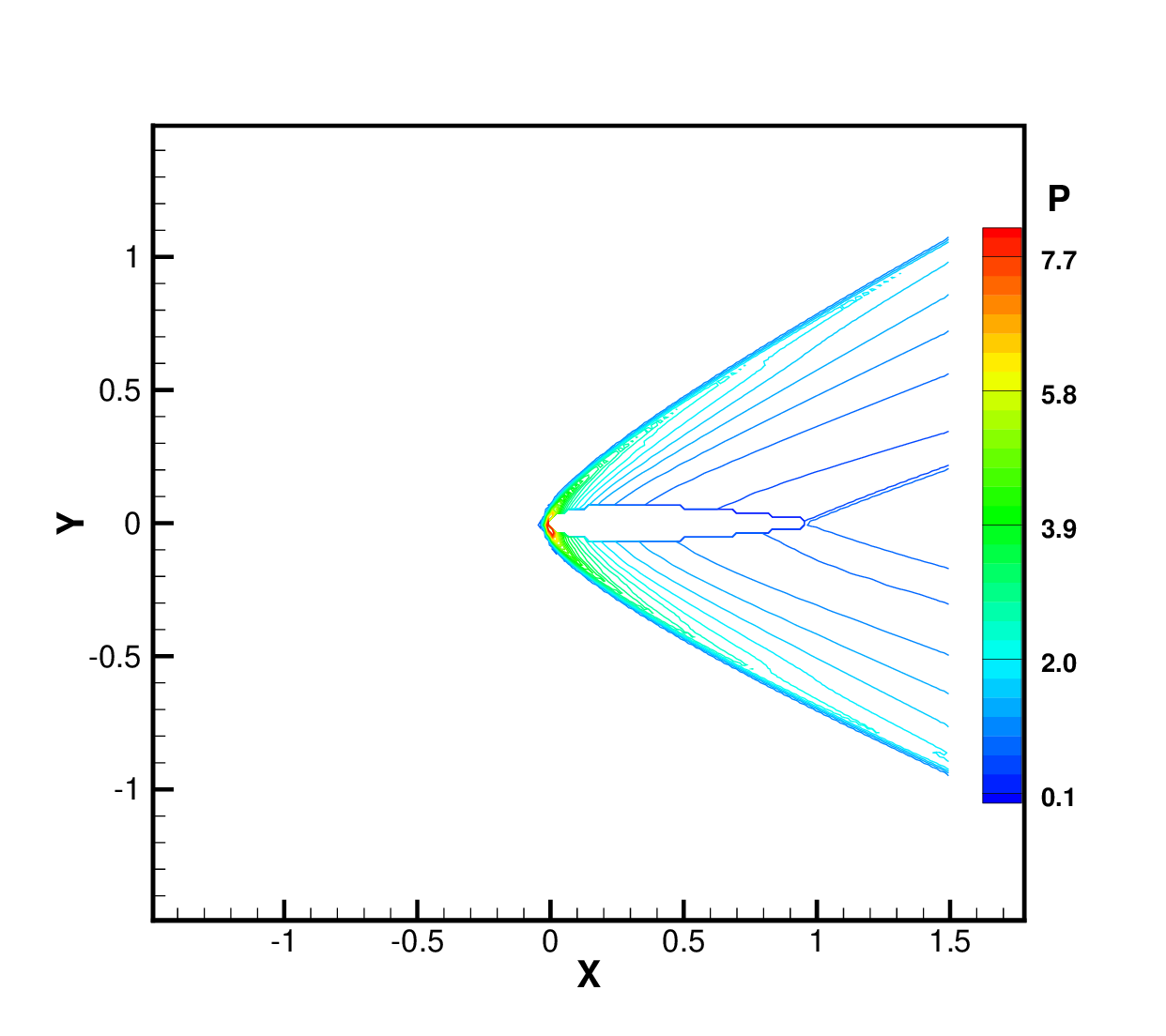}
\end{minipage}%
}%

\subfigure[FS-HAUSWENO-LLF]{
\begin{minipage}[t]{0.4\linewidth}
\centering
\includegraphics[width=2.6in]{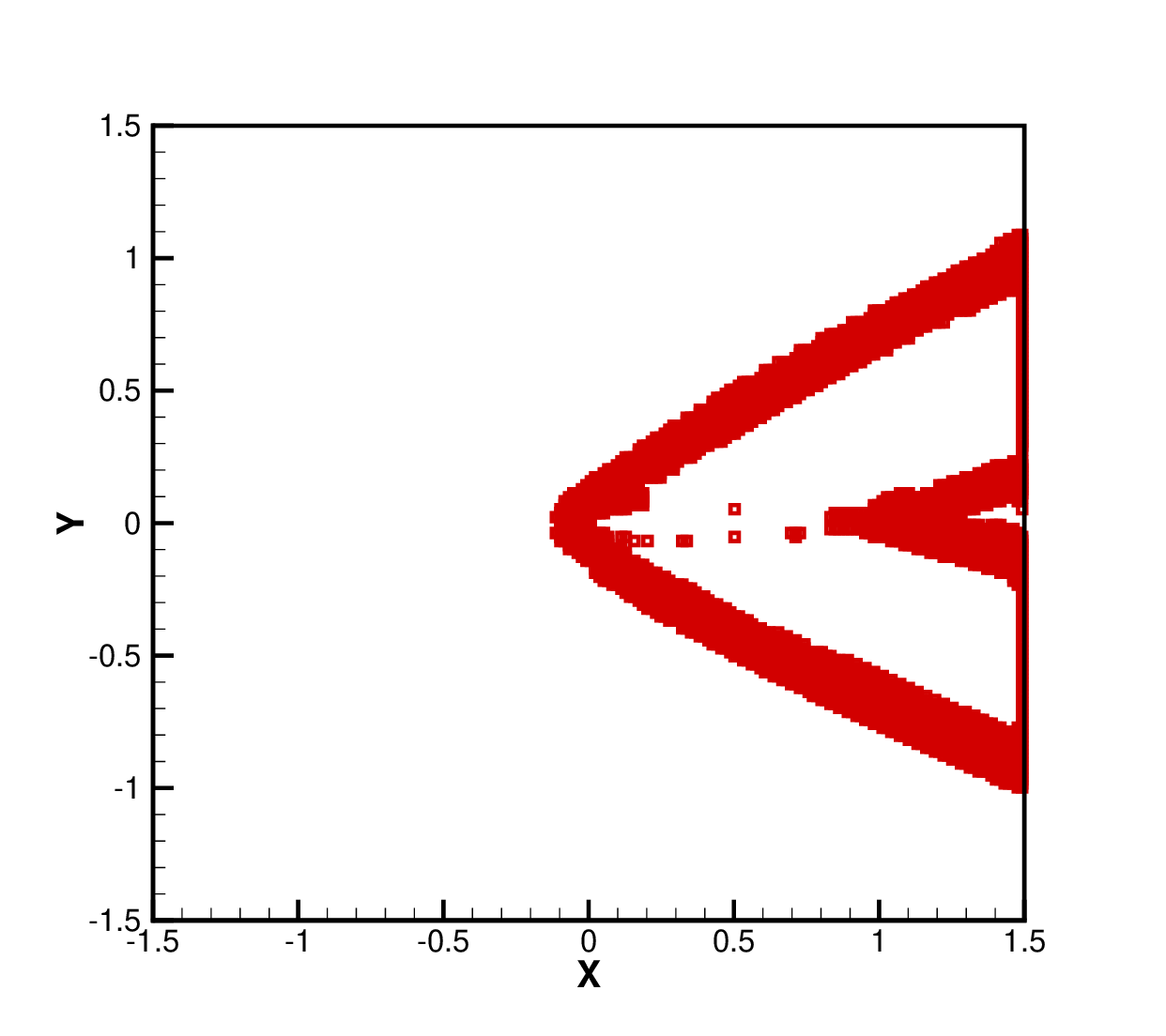}
\end{minipage}%
}
\subfigure[FS-HAUSWENO-HLLC]{
\begin{minipage}[t]{0.4\linewidth}
\centering
\includegraphics[width=2.6in]{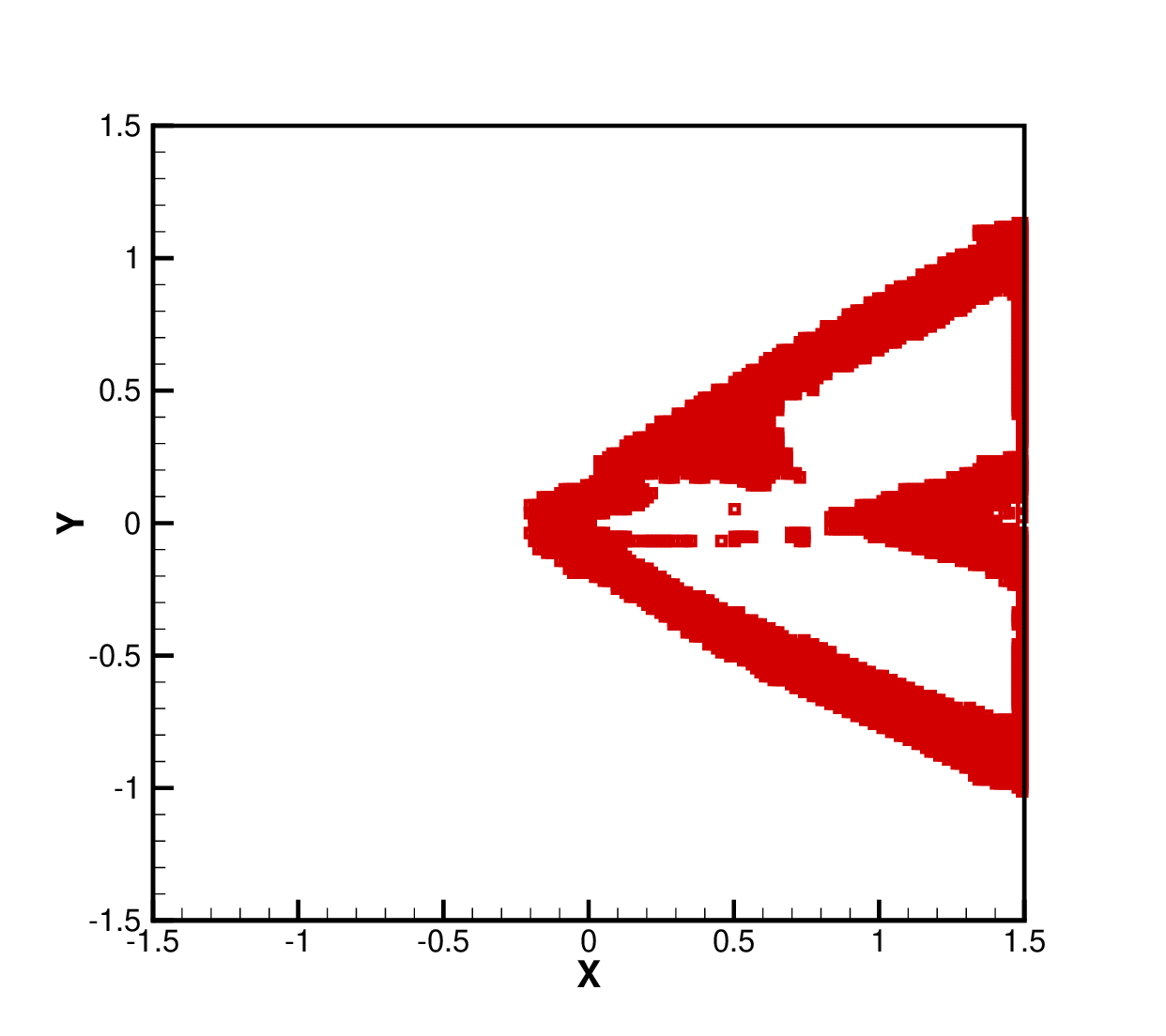}
\end{minipage}%
}%
\caption{\label{6.2}Example 6: Supersonic flow past an NACA0012 airfoil. 30 equally spaced pressure contour from 0.1 to 8.0 of the converged steady states of numerical solutions by different iterative schemes and the identified troubled-cells where the WENO interpolation is used in the FS-HAUSWENO scheme.}
\end{figure}

\begin{figure}
\centering
\subfigure[RK-HAUSWENO-LLF]{
\begin{minipage}[t]{0.4\linewidth}
\centering
\includegraphics[width=3.0in]{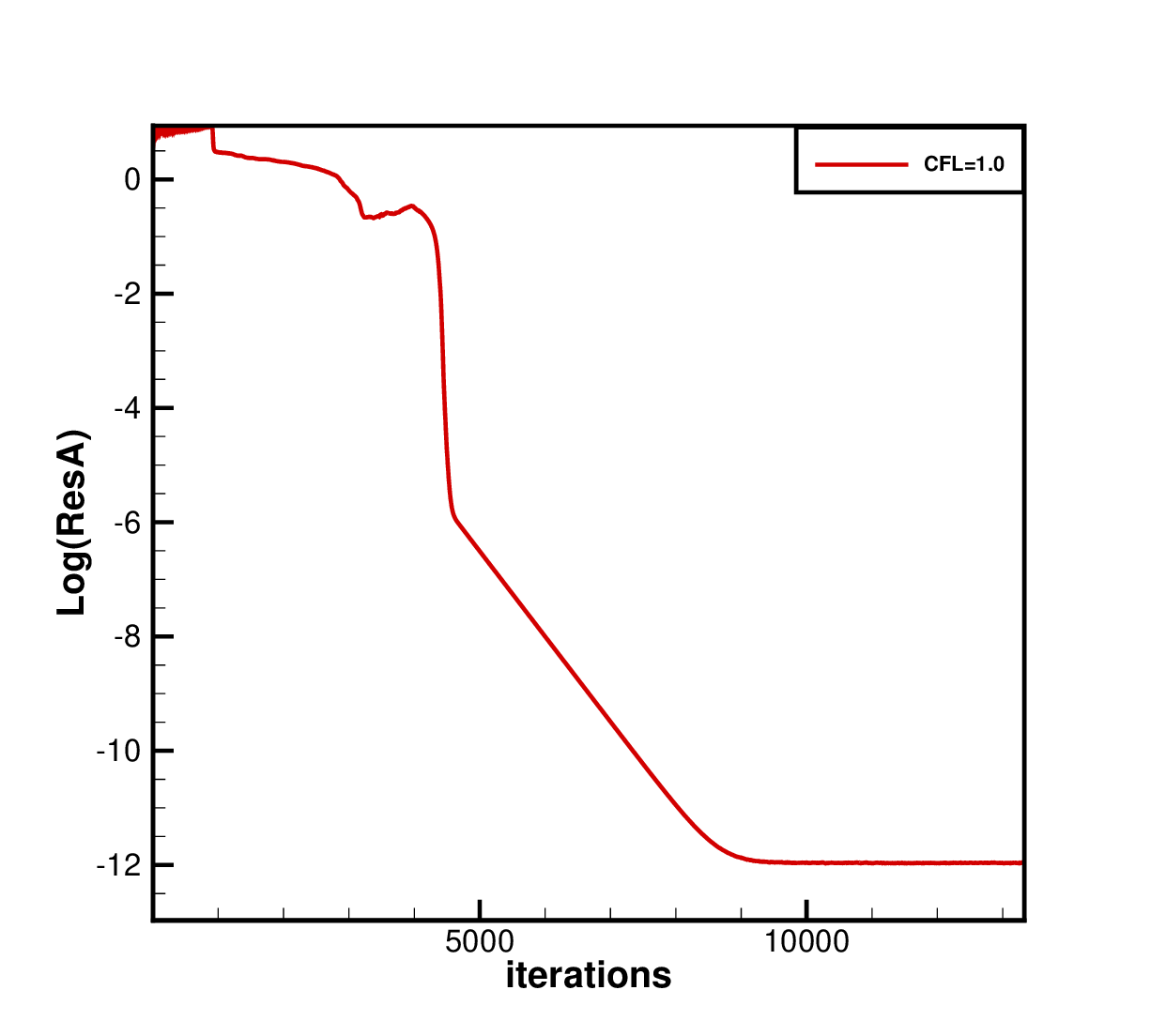}
\end{minipage}%
}
\subfigure[RK-HAUSWENO-HLLC]{
\begin{minipage}[t]{0.4\linewidth}
\centering
\includegraphics[width=3.0in]{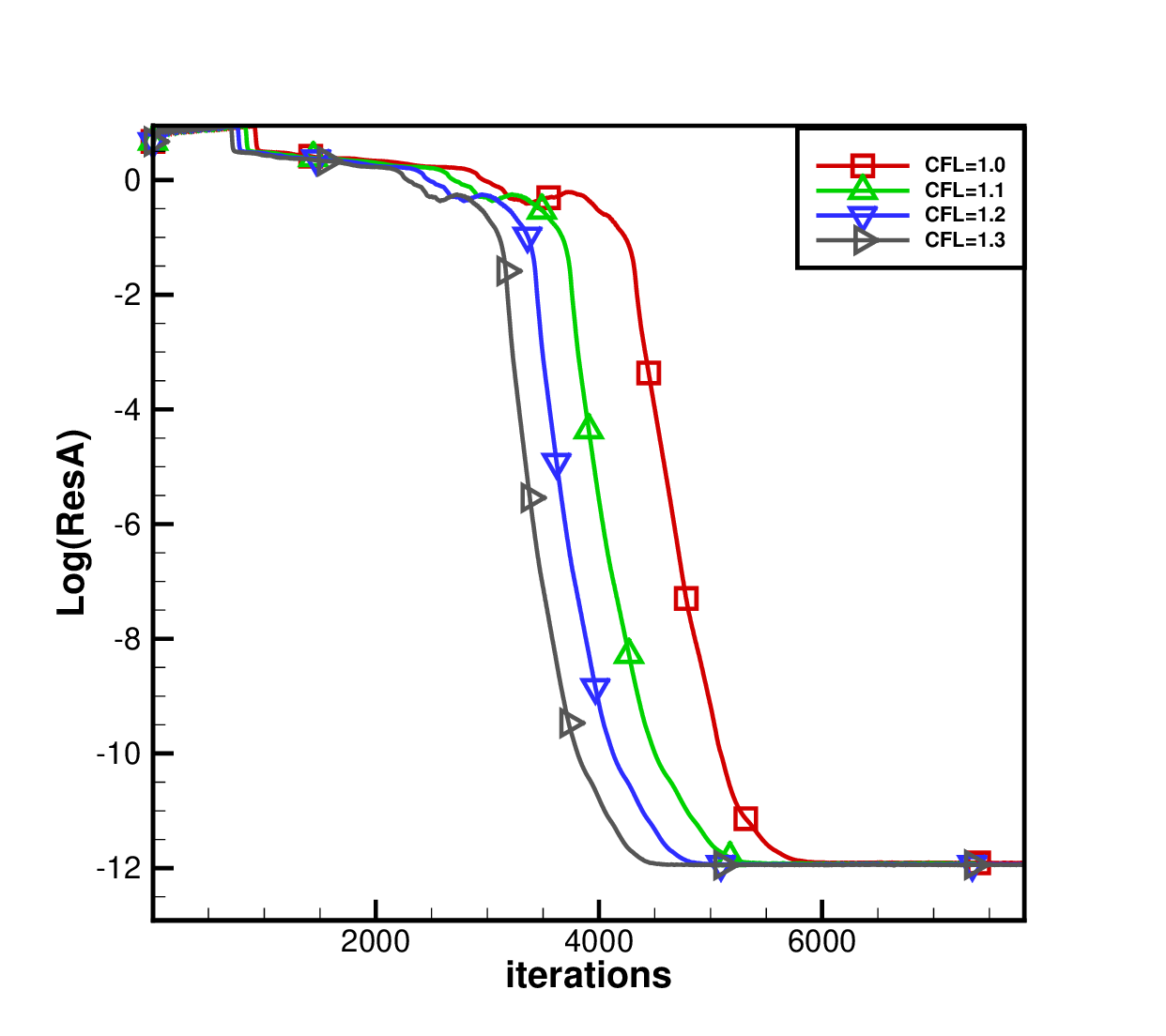}
\end{minipage}%
}%

\subfigure[FS-HAUSWENO-LLF]{
\begin{minipage}[t]{0.4\linewidth}
\centering
\includegraphics[width=3.0in]{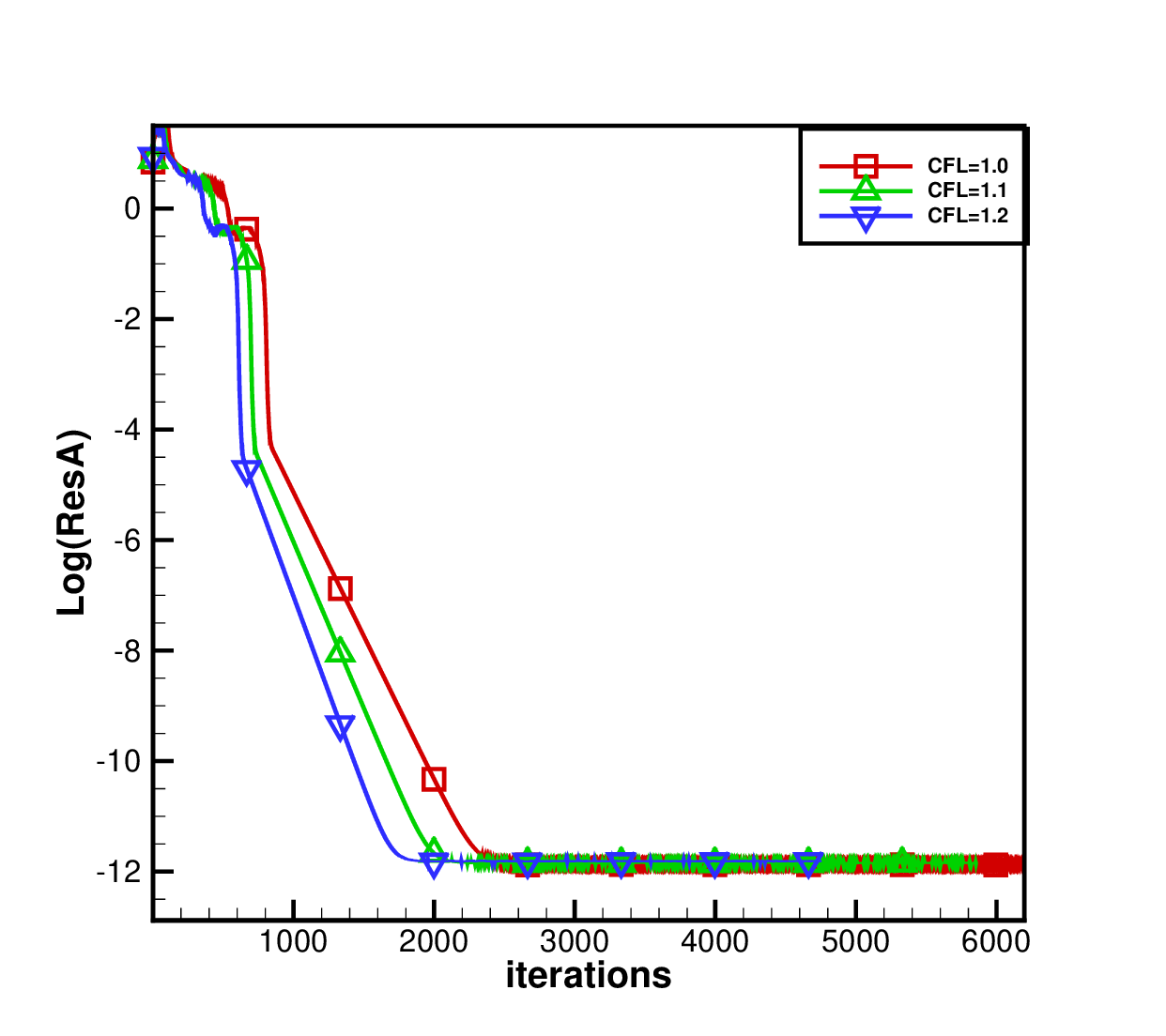}
\end{minipage}%
}
\subfigure[FS-HAUSWENO-HLLC]{
\begin{minipage}[t]{0.4\linewidth}
\centering
\includegraphics[width=3.0in]{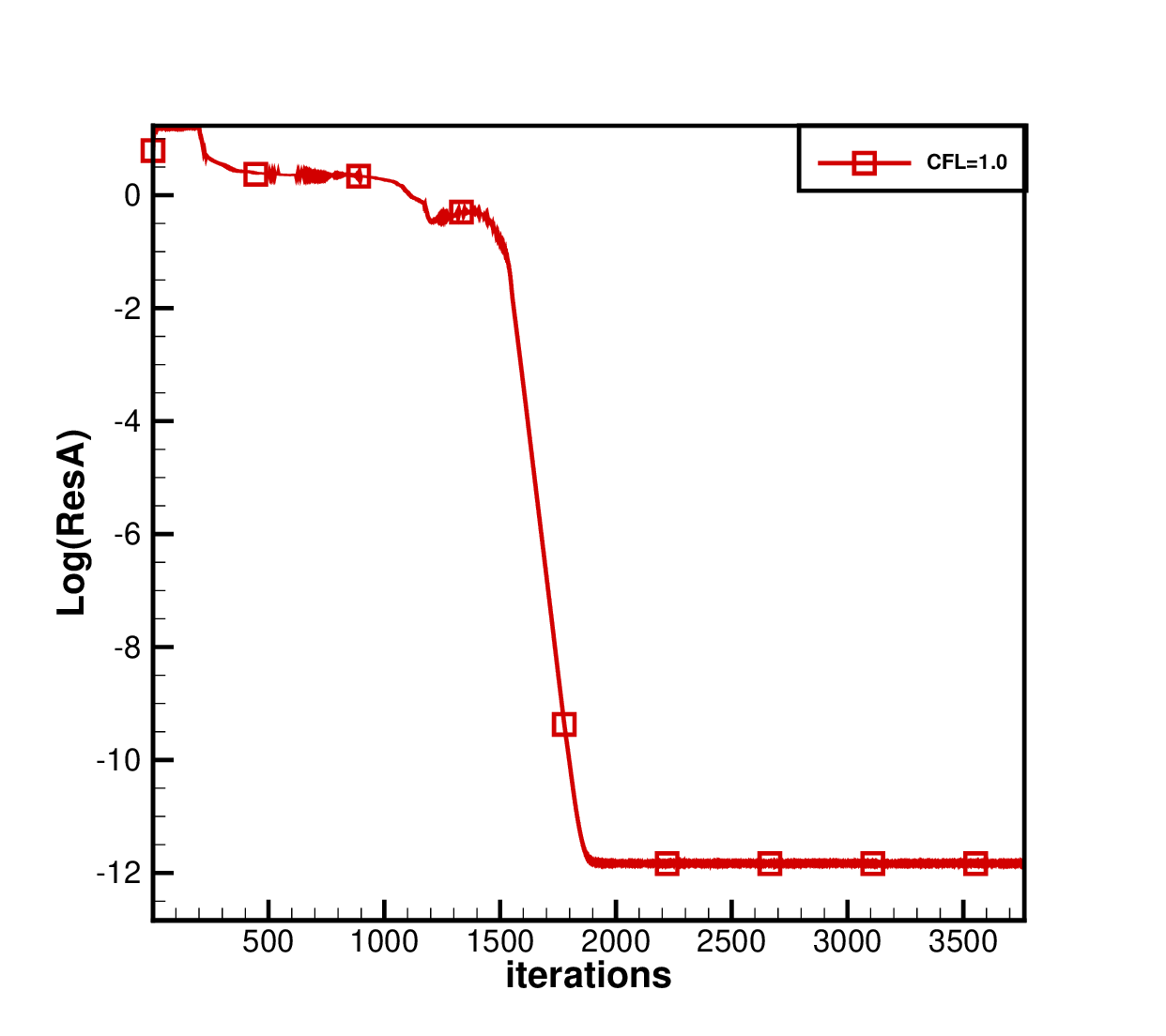}
\end{minipage}%
}%
\caption{\label{6.3}Example 6: Supersonic flow past an NACA0012 airfoil. The convergence history of the residue as a function of number of iterations for different hybrid schemes.}
\end{figure}

We consider the example of an inviscid Euler supersonic flow past a single NACA0012 airfoil configuration \cite{NACA}. The flow has a Mach number $M_{\infty}=3$ and angle of attack $\alpha=3^{\circ}$. The computational domain is $[-1.5,1.5]\times[-1.5,1.5]$, and simulations are performed on the computational grid 200 $\times$ 200. The fifth-order inverse Lax-Wendroff procedure is used to deal with the complex wall surface boundary as in \cite{LILW}. Number of iterations required to reach the convergence criterion threshold value $10^{-11}$, the final time, and the total CPU time when convergence is obtained for the RK-HAUSWENO scheme, the FS-AUSWENO scheme and the FS-HAUSWENO scheme with the LLF flux and the HLLC flux and various CFL numbers up to the largest possible ones are reported in Table \ref{6.1}. The same observations as in the previous examples are obtained, namely, the FS-HAUSWENO scheme is the most efficient one among these three iterative schemes. If the largest CFL number allowed by each scheme to converge to steady states is used, the FS-HAUSWENO scheme saves about $40\%-80\%$ CPU time of that by the RK-HAUSWENO scheme and about $20\%-30\%$ CPU time of that by the FS-AUSWENO scheme, and comparable numerical steady-state solutions are obtained. The contour plots of the pressure variable of the converged steady-state numerical solutions obtained using these two hybrid iterative schemes (FS-HAUSWENO and RK-HAUSWENO) with the LLF flux and the HLLC flux, and the cells where the WENO interpolation is used in the FS-HAUSWENO scheme are shown in Fig.~\ref{6.2}. The results show that the shock waves on the complex surface boundary are well captured, and the hybrid technique identifies the region of troubled-cells well.  Residue history in terms of iterations for these two hybrid schemes is reported in Fig.~\ref{6.3}. It shows that
the iteration residues settle down to the round off error $10^{-11} \sim 10^{-12}$ level and the absolute convergence is achieved.

\begin{table}
		\centering
\begin{tabular}{|c|c|c|c|}\hline
			\multicolumn{4}{|c|}{RK-HAUSWENO-LLF}\\\hline
            $\gamma:$ CFL number & iteration number & final time & CPU time \\\hline
1.0	&8040	&6.00	&841.92 \\\hline
1.1	&not conv	&	& \\\hline
		\end{tabular}
\begin{tabular}{|c|c|c|c|}\hline
			\multicolumn{4}{|c|}{RK-HAUSWENO-HLLC}\\\hline
            $\gamma:$ CFL number & iteration number & final time & CPU time \\\hline
1.0	&5268	&3.96	&392.59 \\\hline
1.1	&4791	&3.95	&356.08 \\\hline
1.2	&4389	&3.95	&324.95 \\\hline
1.3	&4053	&3.95	&302.91 \\\hline
1.4	&not conv	&	& \\\hline
		\end{tabular}
\begin{tabular}{|c|c|c|c|}\hline
			\multicolumn{4}{|c|}{FS-AUSWENO-LLF}\\\hline
            $\gamma:$ CFL number & iteration number & final time & CPU time \\\hline
1.0	&2227	&4.96	&289.28 \\\hline
1.1	&1913	&4.68	&245.23 \\\hline
1.2	&1649	&4.40	&220.65 \\\hline
1.3	&not conv	&	& \\\hline
		\end{tabular}
\begin{tabular}{|c|c|c|c|}\hline
			\multicolumn{4}{|c|}{FS-AUSWENO-HLLC}\\\hline
            $\gamma:$ CFL number & iteration number & final time & CPU time \\\hline
1.0	&1834	&4.18	&218.88 \\\hline
1.1	&not conv	&	& \\\hline
		\end{tabular}
\begin{tabular}{|c|c|c|c|}\hline
			\multicolumn{4}{|c|}{FS-HAUSWENO-LLF}\\\hline
            $\gamma:$ CFL number & iteration number & final time & CPU time \\\hline
1.0	&2137	&4.76	&220.44 \\\hline
1.1	&1833	&4.49	&181.52 \\\hline
1.2	&1581	&4.22	&154.52 \\\hline
1.3	&not conv	&	& \\\hline
		\end{tabular}
\begin{tabular}{|c|c|c|c|}\hline
			\multicolumn{4}{|c|}{FS-HAUSWENO-HLLC}\\\hline
            $\gamma:$ CFL number & iteration number & final time & CPU time \\\hline
1.0	&1833	&4.49	&171.91 \\\hline
1.1	&not conv	&	& \\\hline
		\end{tabular}
		\caption{\label{6.1}Example 6: Supersonic flow past an NACA0012 airfoil. Number of iterations, the final time and total CPU time when convergence is obtained. Convergence criterion threshold value is $10^{-11}$. CPU time unit: second.}
	\end{table}

\section{Concluding remarks}
\label{sec3}
\setcounter{equation}{0}
\setcounter{figure}{0}
\setcounter{table}{0}

In this paper, the absolutely convergent fixed-point fast sweeping WENO methods in \cite{LZZ} are extended to further improve their flexibility, accuracy, and efficiency for solving steady state solution of hyperbolic conservation laws. The designed fast sweeping methods incorporate the AWENO scheme with unequal-sized substencils as the local solver such that arbitrary monotone numerical fluxes can be applied easily and  absolute convergence (i.e., the residue of the fast sweeping iterations converges to machine zero / round off errors) is preserved. The explicit structure of fixed-point fast sweeping methods makes the adoption of AWENO local solver naturally. A novel hybrid technique is developed in the local solver to combine the nonlinear AWENO interpolation with the corresponding high-order linear upwind scheme for more efficient simulations. This new hybrid technique is simpler and shows a better absolutely convergence property in the AWENO scheme than its previous version. Numerical experiments on solving various examples with different schemes are carried out to study the designed method. The accuracy improvement for the new fast sweeping method using a monotone numerical flux with low dissipation such as the HLLC flux is verified. On the other hand, more dissipative flux such as the LLF flux is more robust for the fast sweeping method to achieve absolute convergence. Numerical results show improved efficiency with the hybrid technique, and significant improvement on the computational efficiency of the fast sweeping method over the popular 
third-order  TVD-RK time-marching method for steady state computations. Specifically, up to $65\%-80\%$ CPU time costs can be saved for some complex examples in this paper.

\end{document}